\documentclass[11pt,openleft,titlepage]{amsbook}

\usepackage{lscape,setspace}
\usepackage{graphicx}

\renewcommand{\baselinestretch}{1.5}

\theoremstyle{plain}
\newtheorem{theorem}{Theorem}[chapter]
\newtheorem{proposition}{Proposition}[chapter]
\newtheorem{conjecture}{Conjecture}[chapter]

\theoremstyle{definition}
\newtheorem{lem}{Lemma}[chapter]

\newtheorem{cor}{Corollary}
\newtheorem{example}{Example}[chapter]
\newtheorem{exercise}{Exercise}

\newenvironment{renumerate}%
{%
\begin{enumerate}}%
{\end{enumerate}%
}%

\newenvironment{theo}[1]%
{\begin{theorem}\label{T:#1}}%
{\end{theorem}}

\newenvironment{prop}[1]%
{\begin{proposition}\label{T:#1}}%
{\end{proposition}}

\newenvironment{conj}[1]%
{\begin{conjecture}\label{T:#1}}%
{\end{conjecture}}

\newenvironment{definition}%
{\vskip6pt%
\noindent%
{\bf Definition.}}%
{\vskip6pt}

\newenvironment{remark}%
{\vskip6pt%
\noindent%
{\it Remark.}}%
{\vskip6pt}

\newenvironment{remarks}%
{\vskip6pt%
\noindent%
{\it Remarks}. %
\begin{renumerate}}%
{\end{renumerate}\vskip6pt}

{\begin{exercise}\label{T:#1}}%
{\end{exercise}}

\newenvironment{ex}[1]%
{\begin{example}\label{T:#1}}%
{\end{example}}

\newcommand{\R}[1]{\text{${\mathbb R}^{#1}$}}
\newcommand{\RR}{\text{${\mathbb R}$}}
\newcommand{\C}{\text{$\mathbb C$}}
\newcommand{\N}{\text{$\mathbb N$}}
\newcommand{\Z}{\text{$\mathbb Z$}}
\newcommand{\Q}{\text{$\mathbb Q$}}
\renewcommand{\frak}[1]{\text{$\mathfrak{#1}$}}
\newcommand{\J}{\text{$\mathcal{J}$}}
\newcommand{\SO}[1]{\mathrm{SO}({#1})}
\newcommand{\Spin}[1]{\mathrm{Spin}({#1})}
\newcommand{\so}[1]{\mathfrak{so}({#1})}
\newcommand{\spin}[1]{\mathfrak{spin}({#1})}
\newcommand{\ga}{\text{$\alpha$}}
\newcommand{\gb}{\text{$\beta$}}
\newcommand{\gd}{\text{$\delta$}}
\newcommand{\e}{\text{$\varepsilon$}}

\newcommand{\gO}{\text{$\Omega$}}
\newcommand{\go}{\text{$\omega$}}
\newcommand{\gf}{\text{$\varphi$}}

\newcommand{\Id}{\mathrm{Id}}

\newcommand{\del}{\text{$\partial$}}

\newcommand{\delbar}{\text{$\overline{\partial}$}}
\newcommand{\tensor}{\otimes}
\newcommand{\im}{\mathrm{Im}\,}
\newcommand{\Ker}{\mathrm{Ker}\,}
\newcommand{\mc}[1]{\text{$\mathcal{#1}$}}

\newcommand{\into}{\rightarrow}
\newcommand{\noqed}{\let\qed\relax}
\newcommand{\nil}{\mathrm{nil}}
\newcommand{\idx}{\mathrm{nil}}
\newcommand{\Gg}{\mathfrak{g}}
\newcommand{\om}{\text{$\omega$}}
\newcommand{\gcs}{generalized complex structure}
\newcommand{\gacs}{generalized almost complex structure}
\newcommand{\gcss}{generalized complex structures}
\newcommand{\tgcs}{twisted generalized complex structure}
\newcommand{\gc}{generalized complex}
\newcommand{\gcm}{generalized complex manifold}
\newcommand{\tgcm}{twisted generalized complex manifold}
\newcommand{\gk}{generalized K\"ahler}
\newcommand{\tgk}{twisted generalized K\"ahler}
\newcommand{\tgks}{twisted generalized K\"ahler structure}
\newcommand{\gcy}{generalized Calabi--Yau}
\newcommand{\wrt}{with respect to}
\newcommand{\Ann}{\mathrm{Ann}}

\numberwithin{equation}{chapter}

\begin{document}

\begin{titlepage}
\vspace*{3cm}
\begin{center}
{\Huge{\bf New aspects of the $dd^c$-lemma}}\\
\vspace{7cm}

\Large{Gil R. Cavalcanti} \\[8pt] \large{New College \\ [8pt] University of Oxford}\\
\vspace{2cm}
{\large A thesis submitted for the degree of}\\
{\large {\it Doctor of Philosophy}}\\
\vskip24pt
{\large Trinity Term 2004}
\end{center}
\end{titlepage}
\newpage
\thispagestyle{empty}
~
\newpage
\pagenumbering{roman}
\vspace*{3cm}
\begin{spacing}{1}
\begin{center}
{\LARGE {\bf New aspects of the $dd^c$-lemma}}\\
\vskip12pt
Gil R. Cavalcanti\\
\vskip12pt
{\large {\bf Abstract}}
\end{center}
\vskip6pt
Generalized complex geometry has been recently introduced by Hitchin and Gualtieri and has complex and symplectic geometry as extremal cases, but other nontrivial examples were scarce. In this thesis, we produce new examples of \gcss\ on manifolds by generalizing results from symplectic and complex geometry. We produce \gcss\ on symplectic fibrations over a \gc\ base. We study in some detail different invariant \gcss\ on compact Lie groups and provide a thorough description of invariant structures on nilmanifolds, achieving a classification on 6-nilmanifolds.

Also, drawing on Gualtieri's work, where a generalization of the complex operator $d^c$ was introduced, we study implications of the generalized `$dd^c$-lemma'. Similarly to the standard $dd^c$-lemma, its generalized version induces a decomposition of the cohomology of a manifold and causes the degeneracy of the spectral sequence associated to the splitting $d = \del + \delbar$ at $E_1$. But, in contrast with the $dd^c$-lemma, its generalized version is not preserved by symplectic blow-up or blow-down (in the case of a \gcs\ induced by a symplectic structure) and does not imply formality.

Moreover, we study T-duality for principal circle bundles and how to transport invariant  twisted \gcss\ between T-dual spaces. We use T-duality to find new twisted \gk\ structures on semi-simple Lie groups and prove that no 6-nilmanifold but the torus admits invariant \gk\ structures. We also show that the $dd^c$-lemma is preserved if a \gcs\ is transported via T-duality.

Finally, we tackle the question of formality of a compact orientable $k$-connected $n$-manifold $M^n$, $n= 4k+3,4k+4$. We show that if $b_{k+1}(M)= 1$ then $M$ is formal. In particular, a 1-connected 7-manifold with $b_2=1$ is formal. We prove that if $M$ also satisfies a hard Lefschetz-like property, and $b_{k+1}(M)=2$, then $M$ is formal, while if $b_{k+1}(M)=3$, all the Massey products vanish.  We use these results to study formality of  $G_2$- and $\Spin{7}$-manifolds as well as 8-dimensional symplectic manifolds.
\end{spacing}

\newpage ~
\thispagestyle{empty}
\newpage
\thispagestyle{empty}
\chapter*{Acknowledgements}


Completing a D.Phil. is always a difficult task, and it is never a solo effort. To this end, I owe an enormous debt of gratitude to a number of my colleagues and friends, without whom this thesis would never have been completed.

Without question, the most important person in this endeavor has been my supervisor, Professor Nigel Hitchin, who has been a wonderful and inspirational mentor. It is difficult to overstate the importance of his help and advice through all stages of this quest.  

I am also in great debt to Marco Gualtieri, who introduced me to \gc\ geometry in greater detail and took long hours away from his own thesis writing to explain the main features he had just found about them. I hope the occasions for joint work will be plentiful and this will be an academically fruitful friendship.

I would also like to thank Fred Witt and Ben Hoff for constructive and amusing conversations over blackboard or tea. And from overseas, Alcibiades Rigas, for the great help to get me here in the first place, and Adriano Moura, for keeping in touch and always remembering my birthday, even if I had just forgot his the day before.

Of course, life (even during a D.Phil.) is never all about work, and I would feel a residing sense of shame if I did not mention the amazing friends who I have met along the way, and in particular,  Carlo Maresca-Von Beckh Widmanstetter, David Cleugh, Liz Hanbury, Oliver Thomas, Simon Howes and Wallace Wong. Bwana Aly-Khan Kassam deserves a special mention here. Not only for the great friend he became (and also brother-in-law-to-be), but also for the help in so many different occasions, varying from cheering up to actual down-to-earth dirty computer programing.

I also want to offer heartfelt thanks to my parents, Ana Maria and Alberto, for being there for me and for giving all the long distance support they could throughout this long journey.

I would like to thank to CAPES--MEC in Brasil for their generous support that made everything possible.

Finally a very special place, both in this acknowledgement, and in my heart, goes to Sigrid Heinsbroek. A woman who made the past two years very happy. My time in Oxford would not have been half as good without you.

\tableofcontents

\pagenumbering{arabic}
\chapter*{Introduction}\label{introduction}

The object of our study will be generalized complex manifolds. These are manifolds with a complex structure $\J$ defined in the direct sum of the tangent and cotangent bundles $TM \oplus T^*M$. \J\ is required to be orthogonal \wrt\ the natural pairing:
$$\langle X + \xi , Y +\eta\rangle = \frac{1}{2}(\xi(Y) + \eta(X)),\qquad X,Y \in TM,~~ \xi,\eta \in T^*M$$
and satisfy an integrability condition: the $+i$-eigenspace has to be closed \wrt\ the Courant bracket.

The Courant bracket can also be {\it twisted} by a closed 3-form. And if we use this twisted bracket in the integrability condition we obtain {\it twisted \gcss}.

A \gcs\ has a {\it type} at each point in $M^{2n}$. The type indicates how generic the structure is and varies from 0 (generic) to $n$. Standard (and extremal) examples of \gc\ manifolds are complex (type $n$) and symplectic (type $0$, i.e., `generic') manifolds, products of those and their deformations. Symplectic and complex geometry also provide a series of guiding properties which can be generalized to this new setting. The most relevant to us concerns a generalization of the operators \del, \delbar\ and $d^c$ from complex geometry.

The first remark towards that is the fact that a \gcs\ induces a decomposition of the complex of differential forms
$$\gO^{\bullet}(M) = \oplus_{k=-n}^n \mc{U}^k.$$
In the case of a complex structure, this decomposition is given by
$$ \mc{U}^k = \oplus_{p-q = k}\gO^{p,q}(M).$$
The corresponding decomposition for a symplectic manifold had not been encountered before.

The relevant fact about the decomposition of forms induced by a \gcs\ is that it also decomposes the exterior derivative
$$d:\mc{U}^k \into \mc{U}^{k+1} + \mc{U}^{k-1}$$
which, inspired by the complex case, allows us to define elliptic differential operators \del\ and \delbar:
$$d = \del+ \delbar;\qquad \del: \mc{U}^k \into \mc{U}^{k+1} \qquad \delbar:\mc{U}^k \into \mc{U}^{k-1}.$$
These operators agree with their homonyms from complex geometry, and allow us to define also an analogue of $d^c$ as the real operator $d^{\J} = -i(\del - \delbar)$.

\begin{definition}
A \gcm\ satisfies the {\it $dd^{\J}$-lemma} if
$$ \im d \cap \ker d^{\J} =\im d^{\J} \cap \ker d = \im dd^{\J}$$
\end{definition}

Well known results in complex geometry about the $dd^c$-lemma are that it induces a decomposition of the cohomology of the manifold according to the decomposition of forms into $\gO^{p,q}$ and, more striking, it implies formality \cite{DGMS75}. Also, in \cite{DGMS75}, Deligne {\it et al} argue that the $dd^c$-lemma is preserved by birational transformations, in particular, it is preserved by the blow up of a complex manifold along a submanifold. In the symplectic case, although the decomposition of forms into $\mc{U}^k$ had not been encountered before, the operator $d^{\J}$ induced by a symplectic structure had, and it coincides with Brylinki's {\it canonical derivative} \cite{Br88}. It is a result of Merkulov \cite{Me98} that the $dd^{\J}$-lemma in the symplectic case is equivalent to the hard Lefschetz property (we shall say just ``Lefschetz property'' for short).

\begin{definition}
A symplectic manifold $(M^{2n},\go)$ satisfies the {\it Lefschetz property} if  the maps
$$[\go^{n-k}] : H^{k}(M) \into H^{2n-k}(M)$$
are surjective for $k <n$.
\end{definition}

It was conjectured by Babenko and Taimanov \cite{BT00a} that the Lefschetz property in a symplectic manifold implies formality.

In a different vein, as \gc\ geometry unifies complex and symplectic geometry, it suggests itself as a natural setting in which to study mirror symmetry. In a recent paper \cite{BEM03}, Bouwknegt, Evslin and Mathai introduce a very simple and appealing formalism for T-duality for principal circle bundles and we use this description to study how \gcss\ transform.

We are concerned mostly with two general problems in this thesis:

\begin{renumerate}

\item The search for new examples of \gcss;
\item Implications of the $dd^{\J}$-lemma regarding decomposition of cohomology, formality and its behaviour under symplectic blow-up.
\end{renumerate}

\vskip6pt
{\it Examples of \gcm s}
\vskip6pt

Regarding the search for examples of new structures, we achieve this in different ways. Recalling the early days of symplectic geometry, when the search for nonk\"ahler symplectic manifolds was a `hot problem', we use {\it symplectic fibrations} --- Thurston's solution \cite{Th76} to the symplectic problem --- to produce examples of \gcm s which are not just deformations of products of complex and symplectic manifolds. And we also study these structures on {\it Lie groups} and on {\it nilmanifolds} (Thurston's simplest example of symplectic structure on a symplectic fibration is given by a symplectic nilmanifold).

We prove that any compact Lie group $G^{2n}$ of rank $2k$ admits invariant structures of any type from $n-k$ to $n$. The type $n$ structures coincide with the complex structures determined by Samelson \cite{Sam53} and our technique is a generalization of his.

Regarding nilmanifolds, we find some obstructions to admitting an {\it invariant} \gcs\ based on crude information about the corresponding Lie algebra: the nilpotence step and the dimension of the spaces forming the descending central series. Then we  give a classification of invariant \gcss\ on 6-nilmanifolds. This includes 5 examples that are not products of lower dimensional manifolds and do not admit either symplectic or invariant complex structures. This shows that all 6-nilmanifolds admit \gcss. We also use our results to give an 8 dimensional example that does not admit any invariant \gcs.

Moving to the second of our objectives, we answer the following questions:

\vskip6pt
{\it Does the $dd^{\J}$-lemma induce a decomposition of cohomology?}
\vskip6pt
We establish that it does. By analogy with the Hochschild homology of a Calabi-Yau manifold, we denote this decomposition of the cohomology of a generalized complex manifold $M^{2n}$ by $HH^k(M)$, where $k$ varies from $-n$ to $n$. In the case of a \gcs\ induced by a complex one, these spaces are $HH^k(M) = \oplus_{p-q=k} H^{p,q}(M)$, which agrees with the Hochschield homology of a Calabi--Yau \cite{Cal03}. One of the interesting facts about this decomposition is that generalized complex submanifolds are always represented by elements of $HH^0(M)$.

\vskip6pt
{\it What is this decomposition of forms in the symplectic case?}
\vskip6pt
We prove that the spaces $\mc{U}^k$ in the symplectic case are given by:
$$ \mc{U}^k  = e^{i \go} e^{\frac{\Lambda}{2i}} \gO^{n-k}(M).$$
This expression allows us to find the following unexpected relations between  $\del$ and $\delbar$ and $d$ and $d^{\J}$:
$$ \del(e^{i \go} e^{\frac{\Lambda}{2i}}  \ga) =  e^{i \go} e^{\frac{\Lambda}{2i}} d\ga$$
$$- 2i \delbar ( e^{i \go} e^{\frac{\Lambda}{2i}} \ga) =  e^{i \go} e^{\frac{\Lambda}{2i}} d^{\J} \ga.$$
With this, we can, for example, prove in the symplectic case that the spectral sequence associated to the splitting $d = \del + \delbar$  degenerates at $E_1$. Given the formulae above for $\del$ and \delbar, this spectral sequence is in fact equivalent to the spectral sequence for the {\it canonical complex} introduced by Brylinski \cite{Br88}.

\vskip6pt
{\it What is the equivalent `$dd^{\J}$-lemma' in the symplectic case?}
\vskip6pt
As we mentioned before, 
it is a result by  Merkulov \cite{Me98} that, in the symplectic case, {\it the $dd^{\J}$-lemma is equivalent to the Lefschetz property}. Due to some missing details in Merkulov's original proof, his result is often contested. We give proofs for the parts missing in Merkulov's paper.

\vskip6pt
{\it How does the Lefschetz property behave under symplectic blow-up?}
\vskip6pt
The first example of a simply-connected nonk\"ahler symplectic manifold was the blow-up of $\C P^n$ along the Kodaira-Thurston manifold \cite{McD84}. The resulting manifold has $b_3 =3$ and hence does not satisfy the Lefschetz property. So, in contrast with the complex case, the symplectic blow-up does not preserve the $dd^{\J}$-lemma in the ambient manifold. We prove that neither does the blow-down, since, depending on the existence of symplectic submanifolds in certain homology classes, we have the following:
\vskip6pt
\noindent
{\it a symplectic manifold that does not satisfy the Lefschetz property can be blown up to a manifold that satisfies this property.}
\vskip6pt
Furthermore we give sufficient conditions for the Lefschetz property to be preserved. Namely, under some co-dimension conditions,
\vskip6pt
\noindent
{\it if both the ambient manifold and submanifold satisfy the Lefschetz property, the blow-up will satisfy this property too.}
\vskip6pt

\vskip6pt
{\it Does the Lefschetz property imply formality?}
\vskip6pt

No. Using our results about blow-up, we are able to construct examples of nonformal symplectic manifolds satisfying the Lefschetz property. This answers in the negative the conjecture by Babenko and Taimanov \cite{BT00a} and is related to an earlier question on formality of simply-connected symplectic manifolds by Lupton and Oprea \cite{LO94}. This again puts the symplectic case in contrast with the complex case, as, in the latter case, the $dd^{\J}$-lemma implies formality.

\vskip6pt
{\it How do \gcss\ and the $dd^{\J}$-lemma behave under $T$-duality?}
\vskip6pt

We study T-duality on principal circle bundles with $NS$-fluxes turned on, following work by Bouwknegt {\it et al} \cite{BEM03,BHM03}. We prove that invariant twisted \gcss\ transform well under T-duality and further the $dd^{\J}$-lemma on one side corresponds to the same property on the other. Using T-duality and results of Fino {\it et al} on SKT structures \cite{FPS02}, we prove that no 6-nilmanifold admits invariant \gk\ structures.
\vskip6pt

In the last chapter, we leave \gc\ geometry and study the question of formality of $k$-connected $n$-manifolds, $n= 4k+3, 4k+4$, with special attention to $G_2$- and $\Spin{7}$-manifolds.  The relevance comes from the fact that other Riemannian manifolds with special holonomy are formal (the question is open only for $G_2$, $\Spin{7}$ and quaternionic K\"ahler  manifolds) and $G_2$- and $\Spin{7}$-manifolds have a property that resembles the Lefschetz property for symplectic manifolds:
$$\int_{M} [\gf]  a^2 < 0, a \in H^2(M)\backslash \{0\},$$
where $\gf$ is the closed $(n-4)$-form determining the $G_2$- or $\Spin{7}$-structure on $M^n$, and hence
\begin{equation}\label{E:equationintheintroduction1}
[\gf]:H^2(M) \stackrel{\cong}{\into} H^{n-2}(M)
\end{equation}
is an isomorphism.
In this chapter we stretch the bounds of Miller's theorem \cite{Mi79} stating that {\it $k$-connected $n$-manifolds are formal if $n \leq 4k+2$}. We prove that the same claim holds in dimensions $4k+3$ and $4k+4$ as long as $b_{k+1}=1$. This is relevant to us as our focus is on simply-connected 7-manifolds. We can also prove formality of any $k$-connected $n$-manifold $M^n$, $n=4k+3, 4k+4$, with $b_{k+1}=2$ for which there is a class $\gf \in H^{n-2k-2}(M)$ such that
\begin{equation}\label{E:equationintheintroduction}
[\gf]:H^{k+1}(M) \stackrel{\cong}{\into} H^{n-k-1}(M).
\end{equation}
Further, if $n=4k+3$ and $b_{k+1}=3$, we can prove that all Massey products vanish. We finish with examples inspired by $G_2$- and $\Spin{7}$-manifolds as well as by symplectic manifolds.

This thesis is organized as follows. In the first Chapter we introduce some basic concepts regarding differential graded algebras which will be needed later. This includes Sullivan minimal models, the definition of formality and the proof that complex manifolds satifying the $dd^c$-lemma are formal \cite{DGMS75}, which makes relevant the question of whether \gcm s satisfying the $dd^{\J}$-lemma are formal. We also include the concept of $s$-formality as introduced in \cite{FM02}, which will reappear in the last chapter. And finally we also establish that, for a formal manifold, the cohomology of the exterior derivative $d$ twisted by a 3-form $H$ is isomorphic to the $H$-cohomology.

In Chapter \ref{gcss}, we introduce the central object of this thesis: \gcm s.  Most of the material of this chapter is already in Gualtieri's thesis \cite{Gu03}. There, we recall from \cite{Gu03} that a \gcs\ induces a decomposition of forms into $\mc{U}^k$, as well as the decomposition of the exterior derivative $d = \del + \delbar$. I am able to introduce new examples through symplectic fibrations and principal torus bundles, including Lie groups. We present Gualtieri's deformation theorem for local deformations of a \gcs. We also introduce \gc\ submanifolds and \gk\ manifolds, as defined in \cite{Gu03}.

Chapter \ref{nilmanifolds} gives examples of \gcss\ on nilmanifolds as well as constraints to the existence of certain types of structures on these spaces.

We move on to the $dd^{\J}$-lemma and its implications in Chapters \ref{ddj-lemma} to \ref{blowup}. In Chapter \ref{ddj-lemma}, we prove that if a \gcm\ satisfies the $dd^{\J}$-lemma, its cohomology splits into $HH^k$, according to the decomposition of forms into $\mc{U}^k$. Conversely, the $dd^{\J}$-lemma is itself equivalent to this decomposition of cohomology and the degeneracy of the canonical spectral sequence at $E_1$. We also prove that a \gc\ submanifold is always represented by an element of $HH^0$.

In the next chapter we deal with the symplectic case. We present Yan's results about the representation of $\frak{sl}(2,\C)$ in the exterior algebra of a symplectic manifold \cite{Ya96} and give the decomposition of differential forms into $\mc{U}^k$ for a symplectic manifold. We also use Yan's results to give his proof that the Lefschetz property is equivalent to the existence of symplectic harmonic representatives in each cohomology class and then complete Merkulov's proof of the equivalence of the Lefschetz property and the $dd^{\J}$-lemma \cite{Me98}.

In Chapter \ref{blowup} we study how the Lefschetz property behaves under symplectic blow-up and show that in some cases it is possible to blow up a manifold that does not have this property and obtain one that does. One the other hand, we prove that if both ambient and embedded manifolds satisfy the Lefschetz property and the codimension is high enough, so will the blow-up. As an application, we use blow-up to give the first examples of nonformal symplectic manifolds satisfying the Lefschetz property, or, in the language of \gc\ geometry, a nonformal \gcm\ which satisfies the $dd^{\J}$-lemma.

In Chapter \ref{tduality} we study T-duality for principal circle bundles. There we show that twisted \gcss\ can be transported to the T-dual and that all the structures defined from a \gcs\ on one side can be transported to the T-dual. 

In the last Chapter we study the question of formality on $k$-connected $n$-manifolds, $n=4k+3, 4k+4$. Using the concept of $s$-formality, we improve slightly Miller's bounds for formality of $k$-connected spaces and give some results for manifolds satisfying \eqref{E:equationintheintroduction}.

\chapter{Minimal Models, Formality and Massey Products}\label{minimalmodels}

The main objects of this chapter belong to rational homotopy theory. This theory begins with the discovery by Sullivan in the 60's that not only the abelian groups $H_i(X,\Z)$ and $\pi_i(X)$, the homology and homotopy groups, of a simply connected topological space $X$ can be rationalized, but also the space itself can be rationalized (to $X_{\Q}$ with $\pi_i(X_{\Q}) = \pi_i(X) \tensor \Q$ and  $H_i(X_{\Q}) = H_i(X) \tensor \Q$) as well as maps between spaces $f:X \into Y$ (to $f_{\Q}: X_{\Q} \into Y_{\Q}$). Rational homotopy theory is the study of properties that depend only on the rational type of the space.

The basic object in this theory are differential graded algebras (DGAs, for short) and Sullivan \cite{Su78} introduces a subclass of those which are minimal in some sense and can be used to model the DGA of  PL-forms on $X$. These `minimal models' can be constructed from the algebra of PL-forms on $X$, contain all the rational homotopic information about $X$ and are homotopy invariant as two minimal models are homotopic if and only if they are isomorphic.

From the computational point of view, when dealing with manifolds, a very important point is that we have de Rham-like theorems. These theorems imply that once we tensor everything with $\RR$, we can just work with differential forms instead of PL-forms and the minimal model obtained from differential forms is just the minimal model obtained from PL-forms tensored with \RR. Hence the minimal model for the differential forms on a simply connected manifold $M$ contains all the rational homotopic information of $M$.

Sullivan's theory can also be adapted to manifolds with nontrivial $\pi_1$. Notably, if $\pi_1$ is nilpotent we still have minimal models and the theory is very similar. The case of a general fundamental group is trickier and is still the subject of research \cite{GHT00}.

One particular property of minimal models will concern us frequently in this thesis: formality.
For manifolds with this property, the minimal model for the manifold is the same as the minimal model for its cohomology algebra and hence, by Sullivan's theorem, all the information about the rational homotopy can be extracted from the cohomology algebra or {\it the rational homotopy type is a formal consequence of the cohomology}. Although we have no use for the particular rational homotopy type of a manifold (and hence Sullivan's theorem does not play an important role in our investigation), formality has further implications that make it an interesting property {\it per se}.

Two properties of formal manifolds will be important to us:
\begin{renumerate}
\item All their Massey products vanish.\\
Indeed, Massey products carry topological information and are constructed from the complex of differential forms, having no analogue construction using just the cohomology algebra. Hence the existence of nontrivial Massey products indicates that there is more to the rational homotopy type of the manifold than can be seen from the cohomology algebra.

\item Their $d_H$-cohomology is isomorphic to the $H$-cohomology, where $H$ is any closed 3-form.\\
Here $d_H = d + H\wedge$ is the exterior derivative $d$ twisted by the 3-form $H$. Its cohomology is the twisted cohomology. While $[H]$ acts on the standard cohomology in the obvious way. The difference between these two cohomologies is given by some high order Massey products, hence, if these vanish, the two cohomologies coincide.
\end{renumerate}
And the main result making this topic so relevant to this thesis: Complex manifolds satisfying the $dd^{c}$-lemma are formal (Deligne {\it et al} \cite{DGMS75}).

Most, if not all, of the results we present here are standard and we refer to \cite{GM81,FHT01} for a thorough discussion on these and related topics. However, references may be hard to find for results on twisted cohomology --- I didn't find any. As some may not be familiar with this language, we go through this theory quickly. We start with the definitions of minimal models and formality in Section \ref{minimalmodelssection} where we also give Deligne's proof of the theorem quoted above. In Section \ref{construction} we sketch the proof that every compact manifold has a minimal model and from there give an alternative definition of formality. Then we introduce the concept of $s$-formality as defined by Fern\'andez and Mu\~noz \cite{FM02} and state their theorem about when $s$-formality implies formality. In Section \ref{masseyproducts}, we introduce Massey products and argue that in a formal manifold, these vanish. Finally in the last section, which is the only with some chance of originality in this chapter, we prove that for a formal manifold the twisted cohomology and the $H$-cohomology are (noncanonically) isomorphic.

\section{Minimal Models}\label{minimalmodelssection}

We start introducing the central object of the theory:

\begin{definition}
A {\it differential graded algebra}, or {\it DGA} for short, is an $\N$ graded vector space $\mc{A}^{\bullet}$ over a field $k$, endowed with a product and a differential $d$ satisfying:
\begin{enumerate}
\item The product maps $\mc{A}^{i}\times \mc{A}^{j}$ to $\mc{A}^{i+j}$ and is graded commutative:
$$a \cdot b = (-1)^{ij} b\cdot a;$$
\item The differential has degree 1: $d:\mc{A}^{k} \into \mc{A}^{k+1}$;
\item $d^2 =0$;
\item The differential is a derivation: for $a \in \mc{A}^{i}$ and $b \in \mc{A}^{j}$
$$d(a\cdot b) = da \cdot b + (-1)^i a \cdot db.$$
\end{enumerate}
\end{definition}

We will only work with the field of the real numbers $\RR$. A nontrivial example of DGA is the complex of differential forms on a manifold equipped with the exterior derivative. Another example is the cohomology of a manifold, with $d=0$ as differential. The cohomology of a DGA is defined in the standard way:

\begin{definition}
The cohomology of a DGA $(\mc{A},d)$ is the quotient
$$H^{\bullet}(\mc{A}) = \frac{\ker d:\mc{A}^{\bullet} \into \mc{A}^{\bullet+1}}{\im d:\mc{A}^{\bullet-1} \into \mc{A}^{\bullet}}.$$
\end{definition} 

Sometimes, given a DGA, $(\mc{A},d)$, one can construct another differential  graded algebra that captures all the information about the differential and which is minimal in the following sense.

\begin{definition}
A DGA $(\mc{M},d)$ is {\it minimal} if it is free as a DGA (i.e. polynomial in even degree and skew symmetric in odd degree) and has generators $e_1, e_2 \dots, e_n, \dots$ such that
\begin{enumerate}
\item The degree of the generators form a weakly increasing sequence of positive numbers;
\item There are finitely many generators in each degree;
\item The differential satisfies $de_i \in \wedge\{e_1, \dots, e_{i-1}\}$.
\end{enumerate}
A {\it minimal model} for a differential graded algebra $(\mc{A},d)$ is a minimal DGA, $\mc{M}$, together with a map of differential graded algebras, $\rho: \mc{M} \into \mc{A}$, inducing isomorphism in cohomology. If $(\mc{A},d)$ is the complex of differential forms on a manifold $M$ we also refer to the minimal model of $\mc{A}$ as the minimal model of $M$.
\end{definition}

Every simply-connected compact manifold has a minimal model which is unique up to isomorphism and later we will see how it can be constructed. It is a result of Sullivan \cite{Su78}, that the rational homotopy groups of a simply-connected manifold, $\pi_n(M) \tensor \Q$, are isomorphic to $\mbox{span}\{e_i\mid deg(e_i) =n\}$, where $e_i$ are the generators of the minimal model. This means that all the information about the rational homotopy type of a simply-connected manifold is contained in its minimal model.

As the cohomology algebra of a manifold with real coefficients is also a DGA we can also construct its minimal model. This minimal model is not necessarily the same as the minimal model for the manifold and in general it carries less information.

\begin{definition}\label{formalitypage}
A manifold is {\it formal} if the minimal models for $M$ and for its cohomology algebra are isomorphic, or equivalently, there is a map of differential graded algebras $\psi:\mc{M} \into H^{\bullet}(M,\RR)$ inducing isomorphism in cohomology, where $\mc{M}$ is the minimal model of $M$.
\end{definition}

Hence formality means that the rational homotopy type of $M$ can be obtained from its cohomology algebra.

\begin{ex}{symmetric spaces}
The most standard examples of formal spaces are compact symmetric spaces. The key feature of those is that they have a natural metric for which the product of harmonic forms is harmonic. Hence the map $\phi: H^{\bullet}(M) \into \gO^{\bullet}(M)$ which associates to each cohomology class its unique harmonic representative is actually a map of differential graded algebras inducing an isomorphism in cohomology. If we let $\psi:\mc{M} \into H^{\bullet}(M)$ be the minimal model for $H^{\bullet}(M)$, then, by letting $\rho = \phi \circ \psi: \mc{M} \into \gO(M)$, we obtain a map of DGAs from the minimal algebra $\mc{M}$ inducing an isomorphim in cohomology. Therefore $(\mc{M},\rho)$ is the minimal model for $M$, which is formal.
\end{ex}

This example implies that spheres, projective spaces, Grassmannians and Lie groups, amongst other spaces, are formal. 
More surprising, we have the following result.

\begin{theo}{formality1}
{\em (Deligne} et al {\em \cite{DGMS75}):} A compact complex manifold satisfying the $dd^c$-lemma is formal.
\end{theo}
\begin{proof} Let $M$ be a manifold satisfying the $dd^c$-lemma and let $(\gO,d)$ be the algebra of
differential forms on $M$ with differential $d$,
$(\gO_c,d)$ the algebra of $d^c$-closed forms also with
differential $d$ and $(H_{d^c}(M),d)$ the cohomology of
$M$ with respect to $d^c$, again with differential $d$. Then we
have maps
$$ (\gO,d) \stackrel{i}{\longleftarrow} (\gO_{c},d)
\stackrel{\pi}{\longrightarrow} (H_{d^c}(M), d),$$
and, as we are going to see, these maps induce an isomorphism on cohomology
and the differential of $(H_{d^c}(M), d)$ is zero.

\begin{renumerate}
\item $i^*$ is surjective:

Let \ga\ be a closed form and set $\beta = d^c\ga$. Then $d\beta = d
d^c \ga = - d^c d \ga = 0$, so $\beta$ satisfies the conditions of the
$dd^c$-lemma, hence $\beta = d^c d \gamma$. Let $\tilde{\ga} = \ga -
d\gamma$, then $d^c\tilde{\ga} = d^c\ga - d^c d \gamma = \beta - \beta
=0$, so $[\ga] \in \im(i^*)$.

\item $i^*$ is injective:

If $i^* \ga$ is exact, then \ga\ is $d^c$-closed and exact, hence by
the $dd^c$-lemma $\ga = dd^c\beta$, so \ga\ is the derivative of a
$d^c$-closed form and hence its cohomology class in $\gO_c$ is
also zero.

\item The differential of $(H_{d^c}(M), d)$ is zero:

Let \ga\ be $d^c$-closed, then $d\ga$ is exact and $d^c$ closed so, by
the $dd^c$-lemma, $d\ga = d^cd \beta$ and so it is zero in $d^c$-cohomology.

\item $\pi^*$ is onto:

Let \ga\ be $d^c$-closed. Then, as above,  $d\ga  = dd^c \beta$. Let
$\tilde{\ga} = \ga - d^c \beta$, and so $d\tilde{\ga} = 0$ and
$[\tilde{\ga}]_{d^c} = [\ga]_{d^c}$, so $\pi^*([\tilde{\ga}]_d)= [\ga]_{d^c}$.

\item $\pi^*$ is injective:

Let \ga\ be closed and $d^c$-exact, then the $dd^c$-lemma implies
that \ga\ is exact and hence $[\ga] = 0$ in $\mc{E}_c$.
\end{renumerate}

With this, we have seen that $(H_{d^c},d)$ is isomorphic to the
cohomology algebra $H^{\bullet}(M)$, since its differential is zero and all the
maps induce isomorphisms on cohomology. Now let $\tilde{\rho}: \mc{M} \into
\gO_c$ be the minimal model for $(\gO_c,d)$, then $\rho = i \circ
\tilde{\rho}: \mc{M} \into \gO$ is the minimal model for the complex
of differential forms and $\psi = \pi \circ \tilde{\rho}: \mc{M} \into
H_{d^c}$ induces an isomorphism on cohomology, so $M$ is formal.
\end{proof}

\begin{remark}
A vital ingredient in this proof is that $d^c$ satisfies the Leibniz rule $d^c(a\wedge b) = d^ca \wedge b + (-1)^{|a|}a \wedge d^cb$. This implies that the complex $\gO_c$ of $d^c$-closed forms is actually a graded algebra.
\end{remark}

\section{Construction of Minimal Models and $s$-formality}\label{construction}

As mentioned before, every compact simply-connected manifold has a minimal model. The proof of this fact is done in a constructive way and the fundamental tool is a Hirsch extension.

If $\mc{A}$ is a differential algebra,
$V$ a finite dimensional vector space and $d:V \into
\ker(d)\cap (\mc{A}^{n+1})$ a linear map, we can form the {\it Hirsch extension}
$\mc{A}\tensor_d (\wedge V)_n$ as follows: the elements of
$\mc{A}\tensor_d (\wedge V)_n$ are elements of the tensor product of
$\mc{A}$ with the free algebra generated by $V$ in degree $n$, the
differential is defined by linearity plus the product rule
$d(a\tensor v) = da \tensor v + (-1)^{\deg(a)} a \wedge dv \tensor
1$. We denote 
$1\tensor v$ by $v$ and $a \tensor 1$ by $a$ and then $a \tensor v$
becomes simply $a \wedge v$.

In doing such an extension one is, amongst other things,

\begin{renumerate}
\item killing cohomology in degree $n+1$ since the elements in the image
of $d$  which were closed in $\mc{A}$ but not exact and have
become exact in $\mc{A}\tensor_d (\wedge V)_n$;

\item creating cohomology in degree $n$, namely, the new classes are
given by the kernel of $d:V \into \mc{A}^{n+1}/\im~d:\mc{A}^{n} \into \mc{A}^{n+1}$.
\end{renumerate}

Therefore, we can always split $V = B^n +Z^n$, where the map $d:B^n \into \mc{A}^{n+1}/\im~d:\mc{A}^{n} \into \mc{A}^{n+1}$ is an injection and $Z_n = \ker~d:V \into \mc{A}^{n+1}/\im~d:\mc{A}^{n} \into \mc{A}^{n+1}$, although this splitting is not canonical.

\begin{theo}{existence of minimal models}
Every compact simply-connected manifold has a minimal model.
\end{theo}
\begin{proof}
We start with the free DGA $\mc{M}_2$ generated by $H^2(M)$ in degree 2 and with vanishing differential. We define $\rho_2: \mc{M}_2 \into \gO^{\bullet}(M)$ by choosing arbitrary representatives for the cohomology classes in $H^2(M)$ and extend it to higher symmetric powers of $H^2(M)$ so that it is a map of algebras.

Not all elements in $Sym^2(H_2(M)) = \mc{M}_2^4$ represent nonvanishing cohomology classes, though. And also $\mc{M}^2_3 = \{0\}$. The next stage in the construction of the minimal model fixes these two particular problems. Let
$$V^3 = (\ker Sym^2(H^2(M)) \into H^4(M)) \oplus H^3(M),$$
with a map $d: V^3 \into \mc{M}_2^4$ defined by
$$ d\tilde{a} =a, \mbox{ for } \tilde{a} \in \ker Sym^2(H^2(M)\into H^4(M));$$
$$ da = 0, \mbox{ for } a \in H^3(M),$$
where we used $\tilde{a}$ to indicate the element in $\ker Sym^2(H^2(M)) \into H^4(M) < V^3$ which corresponds to $a \in Sym^2(H^2(M)) <\mc{M}_2^4$ when these two spaces are identified.
This allows us to make the Hirsch extension of $\mc{M}_2$ by $V^3$, which we denote by $\mc{M}_3$. We extend $\rho_2$ to $\rho_3:\mc{M}_3 \into \gO^{\bullet}{M}$ in the following way: If $\tilde{a} \in \ker Sym^2(H^2(M)) \into H^4(M)$, then $\rho_2(a)$ was defined in the first step and is exact, say $\rho_2(a) = db$, then considering $\tilde{a} \in V_3$ we define
$$\rho_3(\tilde{a}) = b$$
And we also define $\rho_3$ in $H^3(M)$ by choosing a representative for each comology class. Extend $\rho_3$ to the whole $\mc{M}_3$ by requiring it to be a map of algebras.

With these choices, $\rho_3$ induces an injection $H^4(\mc{M}_3) \into H^4(M)$ and isomorphisms in cohomology in degrees 2 and 3. The next step in the construction consists in taking another Hirsch extension, this time by
$$V^4 = (\ker H^5(\mc{M}_3) \into H^5(M)) \oplus H^4(M)/\rho(H^4(\mc{M}_3)).$$
Again the first set of generators kills the degree 5 cohomology classes in $\mc{M}_3$ that do not exist in $H^5(M)$, assuring that $H^5(\mc{M}_4)$ injects in $H^5(M)$ and the second introduces the degree 4 cohomology classes that are missing in $H^4(\mc{M}_3)$ when compared to $H^4(M)$, guaranteeing $H^4(\mc{M}_4) \cong H^4(M)$.

This procedure gives us a family of free DGAs, each obtained by a Hirsch extension from the previous one:
$$ \mc{M}_2 < \mc{M}_3 < \cdots \mc{M}_n < \cdots,$$
and maps $\rho_n:\mc{M}_n \into \gO^{\bullet}(M)$ such that:
\begin{renumerate}
\item $\rho_n|{\mc{M}_k} = \rho_k$ for $k< n$, hence induce a map in the limit DGA $\rho: \mc{M} \into \gO^{\bullet}(M)$;
\item The induced maps in cohomology are isomophisms $\rho^*: H^k(\mc{M}_n) \cong H^k(M)$, for $k < n$ and therefore, in the limit, $H^k(\mc{M}) \cong H^k(M)$ for all $k$.
\end{renumerate}
Thus the limit $\mc{M}$ of the $\mc{M}_k$ is the minimal model for $M$.
\end{proof}
We should remark that if the fundamental group of a manifold is nilpotent, then a similar construction of a minimal model can still be carried out \cite{GM81}. There is also a concept of minimal model for non-simply-connected spaces, but we will have no need to consider that.

The main reason to present the proof of this theorem is that it is a constructive one and sometimes properties about the minimal model are more easily obtained using this point of view. For example we have an alternative description of formality (see \cite{FHT01,GM81}, for a proof).

\begin{theo}{alternativeformal}
Let $V^i$ be the spaces introduced in degree $i$ when constructing the minimal model of a manifold $M$. $M$ is formal if and only if there is a splitting $V^i = B^i \oplus Z^i$, such that
\begin{enumerate}
\item $d:B^i \into \mc{M}^{i+1}$ is an injection;
\item $d(Z^i) = 0$;
\item If an element of the ideal $\mc{I}(\oplus_{i \in \N} B^i) < \mc{M}$ is closed, then it is exact.
\end{enumerate}
\end{theo}

This characterization of formality allows one to consider weaker versions. Notably, Fern\'andez and Mu\~noz introduced in \cite{FM02} the useful concept of $s$-formality.

\begin{definition}
A manifold is $s$-formal if there is a choice of splitting $V^i = Z^i \oplus B^i$ satisfying (1) and (2) above and such that every closed element in the ideal $\mc{I}(\oplus_{i \leq s}B^i) < \mc{M}_s$ is exact in $\mc{M}$.
\end{definition}

Clearly this is a weaker concept, in general, but it is also obvious that if an $n$-dimensional manifold is $n$-formal, it is formal. The following surprising result of Fern\'andez and Mu\~noz shows that sometimes the weaker condition of $s$-formality implies formality.

\begin{theo}{fernandezmunoz}{\em (Fern\'andez and Mu\~noz \cite{FM02})}: If an orientable compact manifold $M^n$ is $s$-formal for $s \geq n/2 -1$, then $M$ is formal.
\end{theo}

The good point of this theorem is that if one wants to prove formality by constructing the minimal model and finding the splitting as in Theorem \ref{T:alternativeformal}, it is not necessary to determine the full minimal model, but only the beginning of it. This is most useful in low dimensions and when one has some information about the behaviour of the Betti numbers, as we will see in Chapter \ref{7manifolds}.

\section{Massey Products}\label{masseyproducts}

So far we have two ways to tell whether a manifold is formal: from the definition of formality and from a (partial) construction of the minimal model and the Theorem of Fern\'andez and Mu\~noz. Now we determine an easy way to test if a manifold is {\it not} formal by considering Massey products.

We will start with the triple products. The ingredients are  $a_{12},a_{23},a_{34} \in
\mc{A}$ three closed elements such that $a_{12}a_{23}$ and $a_{23}a_{34}$ are exact. Then, denoting $\bar{a}= (-1)^{|a|}a$, we define
\begin{equation}\tag{$*$}
\begin{cases}
\overline{a_{12}}a_{23} &= da_{13}\\
\overline{a_{23}}a_{34} &= da_{24}.
\end{cases}
\end{equation}
In this case, one can consider the element
$\overline{a_{13}} a_{34} + \overline{a_{12}} a_{24}$. By the choice of $a_{13}$ and $a_{24}$ this element is closed, hence it represents a cohomology class. Observe,
however, that $a_{13}$ and $a_{24}$ are not well defined and we can change
them by any closed element, hence the expression above does not define a
unique cohomology class but instead an element in the
quotient $H^{\bullet}(\mc{A})/\mathcal{I}([a_1],[a_3])$.

\begin{definition}
The {\it triple Massey product} or just {\it triple product}
$\langle [a_{12}], [a_{23}], [a_{34}]\rangle$, 
of the cohomology classes $[a_{12}]$, $[a_{23}]$ and $[a_{34}]$
with $[a_{12}] [a_{23}]=[a_{23}] [a_{34}]=0$
 is the coset
$$\langle [a_{12}],[a_{23}],[a_{34}]\rangle = [\overline{a_{12}}
  a_{24} + \overline{a_{13}} a_{34}] + ([a_{12}], [a_{34}])
  \in H^{\bullet}(\mc{A})/\mathcal{I}([a_{12}],[a_{34}]),$$
where $a_{13}$ and $a_{24}$ are defined by $(*)$. 
\end{definition}

Now if the triple products $\langle [a_{12}], [a_{23}], [a_{34}]\rangle$ and $\langle [a_{23}], [a_{34}], [a_{45}] \rangle$ vanish simultaneously, i.e., are represented by exact elements for a fixed set of  choices, then one can consider
$$\overline{a_{12}}a_{25} + \overline{a_{13}}a_{35} + \overline{a_{14}}a_{45},$$
where $d a_{14} = \langle a_{12},a_{23},a_{34}\rangle$ and similarly for $a_{25}$. This element is always closed, and again depends on choices, hence is defined on a quotient space of $H^{\bullet}(\mc{A})$. This is the {\it quadruple Massey product} $\langle [a_{12}], [a_{23}], [a_{34}], [a_{45}]\rangle$.

In general, before defining an $n$-Massey product $\langle [a_{12}], \cdots,[a_{nn+1}]\rangle$, one needs vanishing lower order products:
$$da_{ij} = \sum_{i<k<j} \overline{a_{ik}} a_{kj},$$
and the {\it n-product} is represented by
$$\sum_{1<k<n+1} \overline{a_{1k}}a_{kn+1}.$$

\begin{remark}
Although in order to compute Massey products we have to make choices of elements in $\mc{A}$ representing the cohomology classes involved, the Massey product itself is independent of those choices, as one can check that if one of the initial elements is exact (and the product defined) then the product vanishes.
\end{remark}

The importance of Massey products for this work comes from two observations
\begin{renumerate}
\item If $\rho:\mc{M} \into \mc{A}$ is the minimal model for $\mc{A}$ and we have a Massey product $\langle v_{12}, \cdots, v_{n-1n} \rangle$ in $\mc{A}$, we can use the same cohomology classes to get a Massey product in $\mc{M}$. As $\rho$ induces an isomorphism in cohomology, the Massey product vanishes in $\mc{A}$ if and only if it vanishes in $\mc{M}$. In some sense, the minimal model is the natural place in which to define Massey products;
\item Using the notation of Therorem \ref{T:alternativeformal}, there is a decomposition of $\mc{M}$, as a vector space, into
$$ \mc{M} = \mbox{span}(\oplus Z^i) \oplus \mc{I}(\oplus B^i);$$
where every element in the first summand is closed. If an element $a$ is exact, say $a= db$, we can split $b$ according to the decomposition above into $b_1 + b_2$, but the element $b_1 \in  \mbox{span}(\oplus Z^i)$ is closed, so in fact we have that every exact element satisfies:
$$a = d b_2 \in d(\mc{I}(\oplus B^i)).$$
\end{renumerate}

Now assume that $M$ is formal. Then any Massey product gives rise to a product in the minimal model from ({\it i\,}) and we can always make choices so that it lies in $\mc{I}(\oplus B^i)$, from ({\it ii\,}), but, by formality, this implies that the Massey product vanishes.

\begin{theo}{vanishing massey products and formality}
If a compact manifold has nonvanishing Massey products, it is not formal.
\end{theo}

Therefore, later in this thesis, when we are interested in formality and have a particular manifold in mind, we will prove that it {\it fails} to be formal by providing nontrivial Massey products.

\section{Twisted Cohomology and Formality}\label{twistedcohomology}

Another appearance that formality will make in this thesis concerns the relationship between the cohomology of the twisted differential operator $d_H = d + H \wedge$ --- also called the {\it twisted cohomology} --- where $H$ is a closed 3-form, and the $H$-cohomology.

\begin{definition}
The {\it $H$-cohomology} of a manifold is given by
$$\frac{\ker \{[H] : H^{\bullet}(M) \into H^{\bullet+1}(M)\}}{\im \{[H] : H^{\bullet-1}(M) \into H^{\bullet}(M)\}},~ \bullet \in \Z_2.$$
\end{definition}

The problem we want to solve is to describe the twisted cohomology, in terms of the standard $d$-cohomology and the cohomology class $[H]$.

The differential forms are $\Z_2$-graded --- $\gO^{ev/od}$ --- And $d_H$ has degree 1 with respect to this grading. The usual grading of forms by degree gives us a filtration of $\gO^{\bullet}(M)$ by
$$F^p = \oplus_{k\geq p} \gO^p,$$
and $d_H: F^p \into F^p$, preserves the filtration. Therefore this gives us a spectral sequence converging to the twisted cohomology. Denoting by
$$F^p\gO^q = F^p \cap \gO^q, q \in \Z_2,$$
the term $E_2^{p,q}$, $q \in \Z_2$, is given by $E^{p,1} = \{0\}$ and 
$$E_2^{p,0} = \frac{\{\ga \in F^p\gO^p : d_H \ga \in F^{p+2}\gO^{p+1} \}}{\{d\gb : \gb \in  F^{p-1}\gO^{p-1} \} \oplus \{ \ga \in F^{p+1}\gO^{p}\}}.$$
Decomposing \ga\ by degrees, $\ga = \ga_p + \ga_{p+2} +\dots$, the above is equivalent to
$$E_2^{p,0} =   \frac{\{\ga \in F^p\gO^p : d \ga_p =0 \}}{\{\ga \in F^p\gO^p : \ga_p = d\gb \}} \cong H^{p}(M),$$
i.e., the $E_2$ term is just the ordinary cohomology. For this spectral sequence, $E_3 \cong E_2$ and the next nontrivial term is $E_4$, which is
$$E_4^{p,0} = \frac{\{a \in H^{p}(M) : a H =0 \}}{\{b H : b\in H^{p-3}(M) \}}.$$
As the spectral sequence does not stop here necessarily, we conclude that the twisted cohomology is a quotient of the $H$-cohomology.

\begin{ex}{liegroups dh-cohomology}
Recall that a compact semi-simple Lie group has the rational homotopy type of a product of odd spheres of dimension greater than or equal to 3 \cite{Sam52,Bor55}. This implies that for any nonzero twist $H$, the $H$-cohomology is trivial. Hence, in this case, the spectral sequence for the $d_H$-cohomology degenerates at $E_4$ and the $d_H$-cohomology is trivial.
\end{ex}

There are other cases where the $d_H$  cohomology and the $H$-cohomology coincide. And formality is the key.

\begin{ex}{H-Massey product}
Let $M$ be a formal manifold, $H$ a closed 3-form and $a$ a closed form such that $H \wedge a$ is exact. Then one can consider the triple product $\langle H, H, a\rangle$, which vanishes by formality.

As we can compute the Massey product in the minimal model we can take $H \wedge H= d0$ and $- H \wedge a = d a_1$, with $a_1 \in \mc{I}(\oplus B_i)$. Therefore, the vanishing of the product is equivalent to
$$ \langle H, H, a \rangle = -H \wedge a_1 = d a_2.$$
And again we can assume that $a_2 \in  \mc{I}(\oplus B_i) $ 

Hence one can go one step further and consider the Massey product $\langle H, H, H, a \rangle$ and so on. Also, as we are working in the minimal model, we can actually always choose $\langle H, \cdots, H \rangle = d 0$ and still get vanishing Massey products:
$$ da_k = \langle\overset{k}{\overbrace{H, \cdots, H}}, a\rangle,$$
and
$$ d a_{k+1} = \langle\overset{k+1}{\overbrace{H, \cdots, H}}, a\rangle= - H \wedge a_k.$$
Therefore
$$d\left(\sum a_k\right) = -H \wedge \sum a_k,$$
Showing that, in a formal manifold, if $[H \wedge a]=0$, we can create a form of mixed degree, \ga, such that $d_H \ga = d\ga + H \wedge \ga =0$.
\end{ex}

\begin{theo}{degeneracyatE4}
If a manifold $M$ is formal, then the twisted cohomology is isomorphic to the $H$-cohomology.
\end{theo}
\begin{proof}
Let $a$ be a nontrivial $H$-cohomology class, and $a_0$ be a form representing such a class. Then, from Example \ref{T:H-Massey product} we can generate a $d_H$-closed form from $a_0$: $b = a_0 + a_1 + a_2 + \cdots$. It is clear that if $b$ was $d_H$-exact, $a_0$ would be the trivial class in the $H$-cohomology.

By choosing a basis for the $H$-cohomology, we can use the argument above to construct an injective linear map from the $H$-cohomology to the $d_H$-cohomology. As the latter has at most the same dimension of the former, they must be isomorphic.
\end{proof}


\chapter{Generalized Complex Geometry}\label{gcss}

The concept of generalized complex manifold was introduced by Hitchin \cite{Hi03} and studied by Gualtieri in his thesis \cite{Gu03}. It consists of a complex structure, not on the tangent bundle, but on the direct sum of the tangent and cotangent bundles, orthogonal with respect to the natural inner product given by evaluation of forms on vectors:
$$\langle X+ \xi, Y + \eta \rangle = \frac{1}{2}(\xi(Y) + \eta(X)),$$
and satisfying an integrability condition.

 The basic examples of such manifolds are complex and symplectic manifolds. Also, as early examples show \cite{Hi03}, in some cases it is possible to deform complex to symplectic structures on a manifold by a path of \gcss.

Many of the properties shared by complex and symplectic structures are present in the generalized complex case, which provides a unified way to study these seemingly distinct cases and establishes relations between some already known properties. Also, by having complex and symplectic geometry as subcases, generalized complex geometry seems to be a natural environment to study {\it mirror symmetry}. Despite having been introduced only recently, generalized complex geometry is already being used to describe mirror symmetry at different levels in a growing number of papers \cite{BeB04,FMT03,GMPT04,LMTZ04} and we will return to this topic at a later stage in this thesis.

The point of this chapter is to introduce the central object of our study  and its basic properties: generalized complex manifolds. All the results here presented are already in Gualtieri's thesis \cite{Gu03}, except for Sections \ref{fibrations} and \ref{lie groups} and Example \ref{T:fibrations}, which are my own contribution. Here, we aim to introduce a series of concepts that parallel more familiar ones in complex geometry, as the decomposition of forms into $\Omega^{p,q}$, the differential operators $\del,\delbar$ and $d^c$, submanifolds and K\"ahler manifolds. We will also present results about how to construct nontrivial examples through symplectic fibrations and the deformations of known structures.

\section{Linear Algebra of a Generalized Complex Structure}\label{linear algebra of gcs}

For any vector space $V^n$ there is a natural symmetric nondegenerate bilinear  form on $V\oplus V^*$ of signature $(n,n)$:
$$\langle X + \xi, Y + \eta  \rangle = \frac{1}{2}(\xi(Y) + \eta(X)).$$
A \gcs\ on $V$ is a complex structure \J\ on $V\oplus V^*$ orthogonal with respect to the natural pairing. Since $\J^2 = - \Id$ on $V\oplus V^*$, it splits $(V \oplus V^*)\tensor \C$ as a direct sum of $\pm i$-eigenspaces, $L$ and $L'$, and conjugation swaps them: $\overline{L} = L'$. Further, as \J\ is orthogonal, we obtain that for $v,w \in L$,
$$\langle v,w \rangle = \langle \J v, \J w\rangle = \langle iv, iw \rangle = -\langle v, w\rangle,$$
and hence $L$ (and $\overline{L}$) is isotropic with respect to the natural pairing.

Conversely, prescribing such an $L$ as the $i$-eigenspace determines a unique \gcs\ on $V$, therefore we could have defined equally well a \gcs\ on a vector space $V^n$ as being an $n$-complex-dimensional subspace $L < (V\oplus V^*)\tensor \C$ such that:
\begin{itemize}
\item $L\cap \overline{L} = \{0\}$;
\item $L$ is isotropic with respect to the natural pairing.
\end{itemize}

\begin{ex}{complex1}
({\bf Complex Structures}): If $V$ has a complex structure, it induces a \gcs\ on $V$ by letting $L=\Lambda^{0,1}V \oplus \Lambda^{1,0}V^*$ be the $i$-eigenspace. It is clear that $L$ satisfies both conditions above. Using the natural splitting of $V\oplus V^*$ we can represent this \gcs\ in the matrix form:
$$\J = \begin{pmatrix}
-J & 0\\
0 & J^*
\end{pmatrix}.$$
\end{ex}

\begin{ex}{symplectic}
({\bf Symplectic Structures}): A symplectic form \go\ on a vector space $V$ also induces a generalized complex structure on $V$ by letting $L =\{X - i X \lfloor \go : X \in V \}$. The nondegeneracy of \go\ implies that $L\cap \overline{L}= \{0\}$ and skew symmetry implies that $L$ is isotropic. Again $\J$ can be represented in the matrix form as:
$$\J = \begin{pmatrix}
 0& -\go^{-1}\\
\go & 0
\end{pmatrix}.$$
\end{ex}

\begin{ex}{B-field}
A real 2-form $B$ (also called a $B$-field) acts naturally on $V\oplus V^*$ by
$$X +  \xi \mapsto X + \xi - X\lfloor B.$$
If $V$ is endowed with a \gcs, \J\ whose $+i$-eigenspace is $L$, we can consider its image under the action of a $B$-field: $L_B = \{X+ \xi - X\lfloor B : X + \xi \in L\}$. Since $B$ is real, $L_B \cap \overline{L_B} = (\Id -B) L\cap \overline{L} = \{0\}$. Again, skew symmetry implies that $L_B$ is isotropic. In matrix form, if $\J_B$ is $\J$ transformed by $B$, we have
$$\J_B = \begin{pmatrix}
1 & 0\\
-B & 1
\end{pmatrix} \J
\begin{pmatrix}
1 & 0\\
B & 1
\end{pmatrix}.$$
\end{ex}

\begin{ex}{beta-field}
An element $\beta \in \wedge^2 V$ also acts on $V \oplus V^*$ in a similar fashion:
$$X +  \xi \mapsto X + \xi + \xi\lfloor \beta.$$
An argument similar to the one above shows that the $\beta$-transform of a \gcs\ is still a \gcs.
\end{ex}

\subsection{Mukai Pairing and Pure Forms}
One more characterization of a \gcs\ on a vector space $V^n$ can be obtained from an interpretation of forms as spinors.

The Clifford algebra of $V \oplus V^*$ is defined using the natural form $\langle \cdot, \cdot \rangle$, i.e., for $v \in V \oplus V^* \subset Cl(V \oplus V^*)$ we have $v^2 = \langle v,v \rangle$.

Since $V$ and $V^*$ are maximal isotropics, their exterior algebras are subalgebras of $Cl(V\oplus V^*)$. In particular, $\wedge^n V^*$ is a distinguished line in the Clifford algebra and generates a left ideal ($\wedge^{\bullet}V^*$) with the property that any of its elements can be written as $s\sigma$, with $\sigma \in \wedge^nV^*$ and $s \in \wedge^{\bullet}V$. This, in turn, determines an action of the Clifford algebra on $\wedge^{\bullet}V^*$ by
$$(X + \xi)\cdot \gf = X\lfloor\gf + \xi \wedge \gf.$$

If we let \ga\ be the antiautomorphism of $Cl(V\oplus V^*)$ defined on decomposables by
\begin{equation}\label{E:main antiautomorphism}
 \ga(v_1 \tensor \cdots \tensor v_k) = v_k \tensor \cdots \tensor v_1,
\end{equation}
then we have the following bilinear form on $\wedge^{\bullet}V^* \subset Cl(V\oplus V^*)$:
$$ (\xi_1, \xi_2) \mapsto (\ga(\xi_1)\wedge \xi_2)_{top},$$
where $top$ indicates taking the top degree component on the form. If we decompose $\xi_i$ by degree: $\xi_i = \sum \xi_i^j$, with $\deg(\xi_i^j)= j$, the above can be rewritten, in an $n$-dimensional space, as
\begin{equation}\label{E:mukai pairing}
(\xi_1,\xi_2) = \sum_j (-1)^j (\xi_1^{2j}\wedge \xi_2^{n-2j} + \xi_1^{2j+1} \wedge \xi_2^{n-2j-1})
\end{equation}

 This bilinear form coincides in cohomology with the {\it Mukai pairing} \cite{Mu84}, where it was introduced in a $K$-theoretical framework.

Now, given a form $\rho \in \wedge^{\bullet}V^* \tensor \C$ (of possibly mixed degree) we can consider its Clifford annihilator
$$ L_{\rho} = \{v \in (V\oplus V^*)\tensor \C : v \cdot \rho = 0\}.$$
It is clear that $\overline{L_{\rho}} = L_{\overline{\rho}}$. Also, for $v \in L_{\rho}$,
$$0 = v^2 \cdot \rho = \langle v, v\rangle \rho,$$
thus $L_{\rho}$ is always isotropic.
\begin{definition}
An element $\rho \in \wedge^{\bullet}V^*$ is a {\it pure form} if $L_{\rho}$ is maximal, i.e., $\dim_{\C}L_{\rho} = \dim_{\RR}V$.
\end{definition}

Gualtieri shows in his thesis \cite{Gu03} that any pure form is of the form $c e^{B+i\go} \gO$, where $B$ and \go\ are real 2-forms and $\gO$ is a decomposable complex form. The relation between the Mukai pairing and \gcss\ comes in the following:

\begin{prop}{chevalley}{\em (Chevalley \cite{Ch96})}
Let $\rho$ and $\tau$ be pure forms. Then $L_{\rho}\cap L_{\tau} = \{0\}$ if and only if $(\rho,\tau) \neq 0$.
\end{prop}

Therefore a pure form $\rho = e^{B+i \omega}\gO$ determines a \gcs\ if and only if $(\rho,\overline{\rho})  = \gO \wedge \overline{\gO} \wedge \go^{n-k}\neq 0$, where $k$ is the degree of $\gO$ and $V$ is $2n$-dimensional. This also shows that there is no \gcs\ on odd dimensional spaces. And with this we have established that the following are equivalent definitions:
\begin{definition}
A {\it generalized complex structure} on a vector space $V^{2n}$ is
\begin{itemize}
\item A complex structure on $V \oplus V^*$ orthogonal with respect to the natural pairing;
\item A maximal isotropic subspace $L < (V\oplus V^*)\tensor \C$ such that $L \cap \overline{L} = \{0\}$;
\item A line in $\wedge^{\bullet}V^*\tensor\C$ generated by a form $e^{B+i \go}\gO$, such that $\gO$ is a decomposable complex form of degree, say, $k$, $B$ and $\go$ are real 2-forms and $\gO\wedge \overline{\gO}\wedge \go^{n-k} \neq 0$.
\end{itemize}
The degree of the form $\gO$ is the {\it type} of the \gcs. 
\end{definition}

\vskip6pt
{\sc Examples \ref{T:complex1} -- \ref{T:beta-field} revised}: The line in $\wedge^{\bullet}V^*\tensor\C$ that gives the \gcs\ for a complex structure is $\wedge^{n,0}V^*$, while the line for a symplectic structure $\go$ is generated by $e^{i\go}$. If $\rho$ is a pure form determining a \gcs, its $B$-field transform will be associated to the form $e^B \wedge \rho$ and its $\beta$-field transform will be $e^{\beta} \lfloor \rho$.
\vskip6pt

\subsection{The Decomposition of Forms}\label{decomposition of forms}
Using the same argument (used before with $V$ and $V^*$) to the maximal isotropics $L$ and $\overline{L}$ determining a \gcs\ on $V$, we see that $Cl((V\oplus V^*)\tensor \C) \cong Cl(L \oplus \overline{L})$ acts on $\wedge^{2n} L$ and the left ideal generated is the subalgebra $\wedge^{\bullet} L$. The choice of a pure form $\rho$ for the \gcs\ and a volume element $\sigma \in \wedge^{2n} L$ gives an isomorphism of Clifford modules: 
$$\phi:\wedge^{\bullet} L \into \wedge^{\bullet}V^*\tensor \C; \qquad
\phi(s \cdot \sigma) = s \cdot \rho.$$

The decomposition of $\wedge^{\bullet} L$ by degree gives rise to a new decomposition of $\wedge^{\bullet} V^*\tensor \C$ and the Mukai pairing on $\wedge^{\bullet}V^*\tensor \C$ corresponds to the same pairing on $\wedge^{\bullet}L$. But in $\wedge^{\bullet} L$ the Mukai pairing is nondegenerate in
$$\wedge^k L \times \wedge^{2n-k} L \into \wedge^{2n}L$$
and vanishes in $\wedge^k L \times \wedge^{l} L$ for all other values of $l$. Therefore the same is true for the induced decomposition on forms.

\begin{prop}{mukai pairing on U_k}
Letting $U^k = \wedge^{n-k} \overline{L} \cdot \rho < \wedge^{\bullet}V^*\tensor \C$, the Mukai pairing on $U^k \times U^l$ is trivial unless $l=-k$, in which case it is nondegenerate.
\end{prop}

\begin{remark}
The elements of $U^k$ are even/odd forms, according to the parity of $k$, $n$ and the type of the structure. For example, if $n$ and the type are even, the elements of $U^k$ will be even if and only if $k$ is even.
\end{remark}

\begin{ex}{complex decomposition}
In the complex case, we take $\rho \in \wedge^{n,0}V^*\backslash\{0\}$ to be a form for the induced \gcs. Then, from Example \ref{T:complex1}, we have that $\overline{L} = \wedge^{1,0}V\oplus \wedge^{0,1}V^*$, so
$$ U^k = \oplus_{p-q=k} \wedge^{p,q}V^*.$$
Then, in this case, one can see Proposition \ref{T:mukai pairing on U_k} as a  consequence of the fact that the top degree part of the exterior product vanishes on $\wedge^{p,q}V^*\times \wedge^{p',q'}V^*$, unless $p+p'=q+q'=n$, in which case it is a nondegenerate pairing.
\end{ex}
\begin{ex}{B-field decomposition}
If a \gcs\ induces a decomposition of the differential forms into the spaces $U^k$, then the $B$-field transform of this structure will induce a decomposition into $U^k_B= e^B \wedge U^k$. Indeed, by Example \ref{T:B-field} revised, $U_{B}^n = e^{B} \wedge U^n$, and
$$U^k_B = (\Id - B)(\overline{L}) \cdot U^{k+1}_{B}.$$
The desired expression can be obtained by induction.
\end{ex}

\subsection{The Actions of \mc{J} on Forms}

Recall that the group $\Spin{n,n}$ sits inside $Cl(V \oplus V^*)$ as
$$\Spin{n,n} = \{v_{1} \cdots v_{2k} : v_i \cdot v_i = \pm 1; k \in \N\}.$$
And $\Spin{n,n}$ is a double cover of $\SO{n,n}$:
$$
\gf:\Spin{n,n} \into \SO{n,n}; \qquad \gf(v)X = v \cdot X \cdot \ga(v),$$
where \ga\ is the main antiautomorphism of the Clifford algebra as defined in \eqref{E:main antiautomorphism}.

This map identifies the Lie algebras $\spin{n,n} \cong \so{n,n} \cong \wedge^2V \oplus \wedge^2V^* \oplus \mbox{End}(V)$:
$$\spin{n,n} \into \so{n,n}; \qquad d\gf(v)(X) = [v,X]= v \cdot X - X \cdot v$$

 But, as the exterior algebra of $V^*$ is naturally the space of spinors, each element in $\spin{n,n}$ acts naturally on $\wedge^{\bullet}V^*$.

\begin{ex}{lie algebra action of B}
Let $B = \sum b_{ij}e^i \wedge e^j \in \wedge^2V^* <\so{n,n}$ be a 2-form. As an element of $\so{n,n}$, $B$ acts on $V\oplus V^*$ via
$$X + \xi \mapsto X\lfloor B.$$
Then the corresponding element in $\spin{n,n}$ inducing the same action on $V\oplus V^*$ is given by $\sum b_{ij} e^j e^i$, since, in $\so{n,n}$, we have
$$ e^i \wedge e^j: e_k \mapsto \delta_{ik}e^j- \delta_{jk}e^i.$$
And, in $\spin{n,n}$,
$$d\gf(e^je^i)e_k =  (e^j e^i) \cdot e_k - e_k \cdot (e^j e^i) = e^j \cdot (e^i\cdot e_k) - (e_k \cdot e^j)\cdot e^i = \delta_{ik}e^j- \delta_{jk}e^i.$$
Finally, the spinorial action of $B$ on a form \gf\ is given by
$$\sum b_{ij} e^j e^i \cdot \gf = - B \wedge \gf.$$
\end{ex}

\begin{ex}{lie algrebra action of beta}
Similarly, for $\beta = \sum\beta^{ij}e_i\wedge e_j \in \wedge^2 V < \so{n,n}$, its action is given by
$$\beta \cdot(X +\xi)= i_{\xi}\beta.$$
And the corresponding element in $\spin{n,n}$ with the same action is $\sum \beta^{ij} e_{j}e_i$. The action of this element on a form \gf\ is given by
$$\beta \cdot \gf = \beta\lfloor \gf.$$
\end{ex}

\begin{ex}{endomorphisms}
Finally, an element of $A = \sum A_i^j e^i \tensor e_j \in End(V) < \so{n,n}$ acts on $V \oplus V^*$ via
$$A(X +\xi) = A(X) + A^*(\xi).$$
The element of $\spin{n,n}$ with the same action is $\frac{1}{2}\sum A^j_i(e_je^i - e^i e_j)$. And the action of this element on a form \gf\ is given by:
\begin{align*}
A \cdot \gf &= \frac{1}{2}\sum A^j_i(e_je^i - e^i e_j)\cdot \gf\\
&= \frac{1}{2}\sum A^j_i (e_j\lfloor (e^i\wedge\gf) - e^i \wedge (e_j\lfloor \gf))\\
& = \frac{1}{2} \sum_i A^i_i \gf - \sum_{i,j} A^j_i e^i\wedge(e_i \lfloor \gf)\\
& = -A^*\gf + \frac{1}{2}\mbox{Tr} A\gf,
\end{align*}
where $A^*\gf$ is the Lie algebra adjoint of $A$ action of \gf\ via
$$A^*\gf(v_1, \cdot v_p) = \sum_i \gf(v_1, \cdots, Av_i, \cdots, v_p).$$
\end{ex}

The reason for introducing this Lie algebra action of $\spin{n,n}$ on forms is  because $\J \in \spin{n,n}$, hence we can compute its action on forms.

\begin{ex}{symplectic lie algebra action}
In the case of a \gcs\ induced by a symplectic one, we have that $\J$ is the sum of a 2-form, \go, and a bivector, $-\go^{-1}$, hence its Lie algebra action on a form \gf\ is
\begin{equation}\label{E:symplectic J}
\J \gf =(-\go \wedge\ -\,\go^{-1}\lfloor)\gf
\end{equation}
\end{ex}

\begin{ex}{complex lie algebra action}
If $\J$ is a \gcs\ on $V$ induced by a complex structure $J$, then its Lie algebra action is the one corresponding to the traceless endomorphism $-J$ (see Example \ref{T:complex1}). Therefore, when acting on a $p,q$-form \gf:
$$ \J \cdot \gf = J^*\gf = i(p-q) \gf.$$
\end{ex}

From this example and Example \ref{T:complex decomposition} it is clear that in the case of a \gcs\ induced by a complex one, the subspaces $U^k < \wedge^{\bullet}V^*$ are the $ik$-eigenspaces of the action of $\J$. This is general.

\begin{prop}{U_k as eigenspaces} The spaces $U^k$ are the $ik$-eigenspaces of the Lie algebra action of $\J$.
\end{prop}
\begin{proof}
Recall from Subsection \ref{decomposition of forms} that the choice of a volume form  $\sigma \in \wedge^{2n}L$ and a pure form $\rho \in  \wedge^{\bullet}V^*$ determining the \gcs\ gives an isomorphism of Clifford modules:
$$\phi:\wedge^{\bullet}L \into \wedge^{\bullet}V^*\tensor \C; \qquad \phi(s \cdot \sigma) = \phi(s \cdot \rho).$$
And the spaces $U^k$ are defined as $U^k = \phi(\wedge^{n+k}L)$. Further, $\J$ acts on $L^* \cong \overline{L}$ as multiplication by $-i$. Hence, by Example \ref{T:endomorphisms}, its Lie algebra action on $\gamma \in \wedge^{n+k}L$ is:
$$
\J\cdot \gamma = -J^*\gamma + \frac{1}{2}\mbox{Tr}J \gamma = i(n+k)\gamma - \frac{1}{2}2n \gamma = ik \gamma.
$$
As $\phi$ is an isomorphism of Clifford modules, for $\alpha \in U^k$ we have $\phi^{-1}\ga \in \wedge^{n+k}L$ and
$$\J \cdot \ga =  \J\cdot \phi(\phi^{-1}\ga) = \phi(\J \phi^{-1} \cdot \ga) = ik\phi(\phi^{-1}\ga) = ik \ga.$$
\end{proof}

And we finish this section with a useful lemma.

\begin{lem}\label{T:maximal real in U_0}
Assume that $W < V \oplus V^*$ is a maximal isotropic subspace invariant under \J. Let $\rho$ be a form whose Clifford annihilator is $W$. Then, $\rho \in U^0$.
\end{lem}
\begin{proof}
Split $v \in W$ into its $L$ and $\overline{L}$ components: $v = v^{1,0} + v^{0,1}$ and $\rho$ into its $U^k$ components: $\rho = \sum \rho_k$. Then, as $W$ is invariant under $J$, we get that
$$v\cdot \rho = (v^{1,0} + v^{0,1})\cdot \rho=0 \qquad \mbox{ and } \qquad Jv \cdot \rho = i(v^{1,0}-v^{0,1})\cdot \rho = 0,$$
which, together with the decomposition of $\rho$, implies that $v^{1,0}\cdot \rho_k = v^{0,1} \cdot \rho_k =0$.
Since $\J$ act on $U^k$ as multiplication by $ik$, we get that $v$ also annihilates $\J\cdot \rho$. But as $W$ is a maximal isotropic, $\J\cdot \rho$ has to be a multiple of $\rho$ and therefore $\rho$ is an eigenvetor of $\J$, say $\rho \in U^k$. But $\rho$ is also real, as $W$ is real, and hence $\rho = \overline{\rho} \in U^{-k}$, implying that $\rho \in U^{k} \cap U^{-k}$ and hence $k =0$.
\end{proof}

\subsection{Metrics on $V\oplus V^*$}\label{metric}

One last topic we have to cover in the linear algebra of the geometry of $T\oplus T^*$ is the metric. The natural pairing is indefinite with signature $(n,n)$ on a $2n$-dimensional space. A metric compatible with this pairing on $V\oplus V^*$ corresponds to an orthogonal 
selfadjoint operator $G:V\oplus V^* \into V \oplus V^*$ such that
$$ \langle Gv,v \rangle >0 \qquad \forall v \neq 0.$$
As $G$ is self adjoint, $G^* = G$ and the orthogonality implies $G^* = G^{-1}$. Therefore $G^2=1$ and $V \oplus V^*$ splits as an orthogonal sum of $\pm 1$-eigenspaces $C_{\pm} < V\oplus V^*$. As $G$ is positive definite, the natural pairing is $\pm$-definite in $C_{\pm}$ and the choice of a pair of such spaces clearly gives us a metric back:
$$\langle \langle V, W \rangle \rangle = \langle V_+,W_+ \rangle - \langle V_-, W_- \rangle,$$
where $V_{\pm}$ and $W_{\pm}$ are the components of $V$ and $W$ in $C_{\pm}$. Therefore a metric is equivalent to a choice of orthogonal spaces $C_{\pm}$ where the natural pairing is definite.

Since $V$ is maximal isotropic, any such $C_+$ can be written as the graph of an element in $\tensor^2 V^*$. More precisely, using the splitting $\tensor^2 V^* = Sym^2 V^* + \wedge^2 V^*$ of a 2-tensor into its symmetric and skew-symmetric parts, we can write $C_+$ as the graph of $b+g$, where $g$ is a symmetric 2-form and $b$ is skew:
$$ C_+ = \{ X + b(X, \cdot)+ g(X,\cdot) | X \in V\}.$$

The fact that the natural pairing is positive definite on $C_+$ places restrictions on $g$. Indeed,
$$ g(X,X) = \langle X + b(X, \cdot)+ g(X,\cdot),X + b(X, \cdot)+ g(X,\cdot) \rangle  >0 \qquad \mbox{ if } X \neq 0.$$
Hence $g$ is a metric on $V$. Further, $C_-$, the orthogonal complement of $C_+$, is also a graph of $b_- + g_-$. But using orthogonality we can determine $g_-$ and $b_-$:
\begin{align*}
 0 & = \langle X + b(X, \cdot)+ g(X,\cdot), Y + b_-(Y, \cdot)+ g_-(Y,\cdot)\rangle \\
&= b(X,Y) + b_-(Y,X) + g(X,Y) + g_-(Y,X),
\end{align*}
which holds for all $X,Y \in V$ if and only if $b_- = b$ and $g_- = -g$ and hence $C_-$ is the graph of $b-g$.

This means that a metric on $V\oplus V^*$ compatible with the natural pairing is equivalent to a choice of metric on $V$ and 2-form $b$.

\section{Generalized Complex Manifolds}
Based on the work of previous section, the following is a natural definition
\begin{definition}
A {\it generalized almost complex structure} on a manifold $M^{2n}$ is one of the following equivalent structures:
\begin{itemize}
\item A section \J\ of $\mathrm{End}(TM\oplus T^*M)$ orthogonal with respect to the natural pairing and satisfying $\J^2 = - \Id$; 
\item An $n$-complex dimensional subbundle $L < (TM \oplus T^*M)\tensor \C$ isotropic with respect to the natural pairing and satisfying $L \cap \overline{L}=\{0\}$;
\item A line bundle in $\wedge^{\bullet}T^*M\tensor\C$ generated locally by a form of the form $e^{B+i \go}\gO$, such that $\gO$ is a decomposable complex form, $B$ and $\go$ are real 2-forms and $\gO\wedge \overline{\gO}\wedge \go^{n-k} \neq 0$ in the points where $\deg(\gO)=k$.
\end{itemize}
A point is {\it regular} for the \gacs\ if it has a neighbourhood where the type of the structure is constant. The line bundle in $\wedge^{\bullet}T^*M\tensor\C$ determining a \gcs\ is the {\it canonical line bundle}.
\end{definition}

To state the integrability condition, we recall that an almost complex structure is integrable if the $+i$-eigenspace is closed under the Lie bracket. In our case, there is no Lie bracket on $TM \oplus T^*M$, but we recollect that the Lie bracket satisfies the following identity when acting on a form \gf\ (see \cite{KoNo63}, Chapter 1, Proposition 3.10):
$$2 [v_1,v_2]\lfloor \gf =
 v_1\wedge v_2\lfloor d\gf + d(v_1 \wedge v_2\lfloor\gf) + 2 v_1 \lfloor d(v_2\lfloor \gf) - 2 v_2\lfloor d(v_1\lfloor \gf).$$

Now we observe that this formula gives a natural extension of the Lie bracket to a bracket on $TM \oplus T^*M$, acting on forms via the Clifford action:
\begin{equation}\label{E:Courant bracket1}
2[v_1,v_2]\cdot\gf =
 v_1\wedge v_2\cdot d\gf + d(v_1\wedge v_2\cdot\gf)) + 2v_1\cdot d(v_2\cdot \gf) - 2 v_2\cdot d(v_1\cdot \gf).
\end{equation}
Spelling it out we obtain (see Gualtieri \cite{Gu03}, Lemma 4.24):
\begin{equation}\label{E:Courant bracket2}
[X+\xi,Y + \eta] = [X,Y] + X\lfloor d\eta - Y \lfloor d\xi +\frac{1}{2}d(X\lfloor \eta - Y \lfloor \xi).
\end{equation}
This is the Courant bracket introduced in \cite{Co90,CW86} in the context of Dirac structures, of which \gcss\ are the complex analogue.

With this, it is natural to define that a \gacs\ is {\it integrable} if the $+i$-eigenbundle $L$ is closed under the Courant bracket. Of course, this can also be translated to an integrability condition for the forms defining a \gcs. Indeed, if $\rho$ has $L$ as its Clifford annihilator, then $L$ is involutive with respect to the Courant bracket if and only if, for all $v_1,v_2$ sections of $L$ , $[v_1,v_2]\cdot \rho =0$. Using \eqref{E:Courant bracket1}, this is equivalent to
$$ v_1\cdot (v_2 \cdot d\rho)=0,$$
as all the other terms vanish.
If we let $\mc{U}^k$ be the space of sections of the bundle $U^k$, defined in Proposition \ref{T:mukai pairing on U_k},
this is only the case if $d\rho \in \mc{U}^{n-1} = \overline{L}\cdot \mc{U}^n$. Therefore we have the following:
\begin{definition}
A \gacs\ is {\it integrable} if one of the following equivalent properties holds:
\begin{itemize}
\item The $+i$-eigenbundle is involutive \wrt\ the Courant bracket;
\item The exterior derivative satisfies
$$ d :\mc{U}^n \into \mc{U}^{n-1},$$
i.e., for a nonvanishing local section \gf\ of $\mc{U}^n$, $d\gf = (X+\xi)\cdot\gf$, for some $X+ \xi \in (TM\oplus T^*M)\tensor \C$.
\end{itemize}
A generalized complex  manifold for which the canonical bundle admits a nonvanishing closed section is called {\it generalized Calabi--Yau}.
\end{definition}

\begin{ex}{complex2}
An almost complex structure on a manifold $M$ induces a \gacs\ with $+i$-eigenspace $T^{0,1}M \oplus T^{*1,0}M$. If this \gacs\ is integrable, then $T^{0,1}M$ has to be closed \wrt\ the Lie bracket and hence the almost complex structure was actually a complex structure. Conversely, any complex structure gives rise to an integrable \gcs. A complex structure will give rise to a \gcy\ structure only if the canonical bundle $\wedge^{n,0}T^*M$ admits a nonvanishing closed section, i.e., if it is holomorphically trivial.
\end{ex}

\begin{ex}{symplectic2}
If $M$ has a nondegenerate 2-form \go, then the induced \gacs\ will be integrable if for some $X +\xi$ we have
$$de^{i\go} = (X + \xi)\cdot e^{i\go}.$$
The degree 1 part gives that $X\lfloor \go + \xi=0$ and the degree 3 part, that $d\go=0$ and hence $M$ is a symplectic manifold. Also, $e^{i\go}$ is a globally defined closed section of the canonical bundle, hence symplectic manifolds are generalized Calabi--Yau manifolds.
\end{ex}

\begin{ex}{B-field2}
A {\it real closed 2-form} $B$, also called a $B$-field, acts on $TM\oplus T^*M$ and maps \gcss\ to \gcss. Indeed, the action of any closed 2-form (real or complex) is a symmetry of the bracket and these, together with diffeomorphisms of the manifold, form the group of orthogonal symmetries of the Courant bracket (see \cite{Gu03}, Proposition 3.24). As \gcss\ are preserved only by the action of $B$-fields and diffeomorphisms these give the natural concept of equivalence for these structures.
\end{ex}

Similarly to the complex case, the integrability condition places restrictions on $d(\mc{U}^k)$ for every $k$ and hence allows us to define operators $\del$ and $\delbar$. This is the object of the following:

\begin{theo}{U_k to U_k+1}{\em (Gualtieri \cite{Gu03}, Theorem 4.3)} A \gacs\ is integrable if and only if $d:\mc{U}^k \into \mc{U}^{k+1}\oplus \mc{U}^{k-1}$.
\end{theo}
\begin{proof}
It is clear that if $d:\mc{U}^k \into \mc{U}^{k+1}\oplus \mc{U}^{k-1}$ for all $k$, then, in particular, for $k=n$ we get the integrability condition
$$d:\mc{U}^n \into \mc{U}^{n-1},$$
as $\mc{U}^{n+1} = \{0\}$.

In order to obtain the converse we shall first prove that
\begin{equation}\label{E:useless equation of this theorem}
d:\mc{U}^{k} \into \oplus_{j \geq k-1}  \mc{U}^{j}.
\end{equation}

This is done is by induction on $k$, $k$ starting at $n$ and going downwards.
The first step is just the integrability condition. Assuming that the claim is true for $k' > k$, let $v_1$, $v_2$ be sections of $L$ and $\gf \in \mc{U}^k$, then by equation \eqref{E:Courant bracket1} we have
$$ v_1\wedge v_2\cdot d\gf  = 2 [v_1,v_2]\cdot \gf - d(v_1 \wedge v_2\cdot\gf)) -2 v_1 \cdot d(v_2\cdot \gf) + 2 v_2\cdot d(v_1\cdot \gf)$$
but the right hand side is, by inductive hypothesis, in $\oplus_{j \geq k+1} \mc{U}^{j}.$
Therefore, so is $v_2\wedge v_1 \cdot d\gf$. Since $v_1$ and $v_2$ are sections of $L$, we conclude that $d\gf \in \oplus_{j\geq k-1} \mc{U}^{j}$.

In order to finish the proof of the converse, we remark that conjugation swaps $\mc{U}^{\pm k}$, but preserves $d$, as it is a real operator, i.e., $ d\gf = \overline{d\overline{\gf}}$. Thus, for $\gf \in \mc{U}^{k}$, $\overline{\gf} \in \mc{U}^{-k}$ and, using \eqref{E:useless equation of this theorem},
$$d\gf \in \oplus_{j\geq k-1} \mc{U}^{j} \qquad d\gf = \overline{d\overline{\gf}} \in  \oplus_{j \leq k+1} \mc{U}^{j}.$$
Showing that $d\gf \in \mc{U}^{k-1} \oplus \mc{U}^k \oplus \mc{U}^{k+1}$. Finally, if $\gf \in \mc{U}^k$ is an even/odd form, then all the elements in $\mc{U}^{k}$ are even/odd whereas the elements of $\mc{U}^{k-1} \oplus \mc{U}^{k+1}$ are odd/even. As $d$ has degree 1 \wrt\ the normal grading of $\gO^{\bullet}$, $d\gf$ is odd/even and hence has no $\mc{U}^k$ component.
\end{proof}
\begin{definition}
Let $M$ be a generalized complex manifold. We define the operators
$$\del : \mc{U}^k \into \mc{U}^{k+1};$$
$$\delbar:\mc{U}^k \into \mc{U}^{k-1};$$
as the projections of $d$ onto each of these factors. Also we define $d^{\J} = -i(\del-\delbar)$.
\end{definition}

Similarly to $d^c$, we can find other expressions for $d^{\J}$ based on the action of the generalized complex structure on forms. One can easily check that if we consider the Lie group action of $\J$, i.e, $\J$ acts on $\mc{U}^k$ as multiplication by $i^k$, then
$$d^{\J} = \J^{-1} d \J.$$
And if one considers the Lie algebra action, then $$d^{\J} = [d,\J].$$

As a consequence of Example \ref{T:complex2} it is clear that in the complex case, \del\ and \delbar\ are just the standard operators denoted by the same symbols and $d^{\J}=d^c$. In the symplectic case, $d^{\J}$ corresponds to Kosul's {\it canonical homology derivative} \cite{Kos85} introduced in the context of Poisson manifolds and studied by Brylinski \cite{Br88}, Mathieu \cite{Ma90},Yan \cite{Ya96}, Merkulov \cite{Me98} and others \cite{FIL94,IRTU01} in the symplectic setting. We will study the symplectic case in detail in Chapters \ref{symplectic_ddbar} and \ref{blowup}.

\begin{ex}{B-field4}
If a \gcs\ induces a splitting of $\wedge^{\bullet}T^*$ into subspaces $U^k$, then, according to Example \ref{T:B-field decomposition}, a $B$-field transform of this structure will induce a decomposition into $e^B U^k$. As $B$ is a closed form, for $v \in \mc{U}^k$ we have:
$$d(e^B v) = e^B d v = e^B \del v + e^B \delbar v \in e^B \mc{U}^{k+1} + e^B \mc{U}^{k-1},$$
hence the corresponding operators for the $B$-field transform, $\del_B$ and $\delbar_B$, are given by
$$\del_B = e^B \del e^{-B};\qquad \delbar_B = e^{B} \delbar e^{-B}.$$
\end{ex}

\section{Twisted Generalized Complex Structures}

The idea of a twisted generalized complex structure comes from the observation by \v{S}evera and Weinstein \cite{SW01} that the Courant bracket can be twisted by a closed 3-form $H$ (also called an NS-flux):
$$[X+\xi, Y+ \eta]_H = [X+\xi,Y+\eta]+ Y\lfloor (X \lfloor H).$$
In the same way that the Courant bracket is the bracket induced by $d$ (see equation \ref{E:Courant bracket1}), the twisted bracket can be described as the bracket induced by the twisted differential $d_H = d+ H$:
\begin{equation}\label{E:twisted bracket}
2[v_1,v_2]_H\cdot\gf =
 v_1\wedge v_2\cdot d_H\gf + d_H(v_1\wedge v_2\cdot\gf)) + 2v_1\cdot d_H(v_2\cdot \gf) - 2 v_2\cdot d_H(v_1\cdot \gf).
\end{equation}
Hence, if we use this twisted bracket in the integrability condition we obtain that $L$ is integrable (i.e., closed under the twisted barcket) if and only if $[v_1,v_2]_{H} \cdot \rho=0$, which, similarly to the nontwisted case, is equivalent to $d_H\rho \in \mc{U}^{n-1}$. Although one can twist a \gcs\ by any 3-form, we shall restrict $H$ to be a {\it real}\/ 3-form.
\begin{definition}
An {\it $H$-twisted generalized complex structure}, or just {\it twisted \gcs}, when the twisting 3-form is clear,  is a \gacs\ such that one of the following equivalent properties hold:
\begin{itemize}
\item The $+i$-eigenspace is closed under the twisted Courant bracket;
\item The twisted exterior derivative satisfies
$$ d_H :\mc{U}^n \into \mc{U}^{n-1},$$
i.e., for a nonvanishing local section \gf\ of $\mc{U}^n$, $d\gf = (X+\xi - H)\cdot\gf$, for some $X+ \xi \in (TM\oplus T^*M)\tensor \C$.
\end{itemize}
\end{definition}
It is clear that  Theorem \ref{T:U_k to U_k+1} also holds in this context, with the exterior derivative $d$ replaced by $d_H$, hence we can still define the operators \del\ and \delbar\ by $d_{H} = \del+\delbar$ and $d^{\J} = -i(\del - \delbar)$. Also, the transform of a twisted complex structure by a nonclosed $B$-field is a twisted complex structure with twist $H + dB$. Hence the existence of such structures on a manifold only depends on the cohomology class of $H$.

\section{Symplectic Fibrations}\label{fibrations}

The differential form description of a \gcs\ furnishes also a very pictorial one around {\it regular points}. Indeed, in this case one can choose locally a closed form $\rho$ defining the structure. Then the integrability condition tells us that $\gO\wedge\overline{\gO}$ is a real closed form and therefore its annihilator $Ann(\gO\wedge\overline{\gO}) <TM$ is an integrable distribution. The algebraic condition $\gO\wedge \overline{\gO} \wedge \go^{n-k}\neq 0$ implies that $\go$ is nondegenerate on the leaves and the integrability condition $\gO\wedge d\go =0$, that $\go$ is closed when restricted to these leaves. Therefore, around a regular point, the \gcs\ furnishes a natural symplectic foliation, and further, the space of leaves has a natural complex structure given by $\gO$.

This suggests that symplectic fibrations should be a way to construct nontrivial examples of \gcss\ which are neither complex nor symplectic nor deformations of those, at least, a priori. Next we see that Thurston's argument for symplectic fibrations \cite{Th76,MS00} can also be used in the generalized complex setting.

\begin{theo}{Thurston}
Let $\pi:X \into B$ be a symplectic fibration over a compact twisted generalized complex base with (possibly zero) twist $H$ and  compact fiber $(F,\go)$. Assume that there is $a \in H^2(X)$ which restricts to the cohomology class of $\go$ on each fiber. Then $X$ admits a twisted \gcs\ whose twist is $\pi^*H$.
\end{theo}
\begin{proof}
Let $\tau_0$ be a closed 2-form representing the cohomology class $a \in H^2(X)$ and $U_{\ga}$ be a cover of $B$ by contractible open sets where the fibration and the canonical line bundle are trivial. Let $\rho_{\ga}$ be a form determining the \gcs\ over $U_{\ga}$. Also, let $\go_{\ga}$ be the pull back of $\go$ via the projection $U_{\ga} \times F \into F$. Then, in $U_{\ga} \times F$, there is a 1-form $\lambda_{\ga}$ such that
$$\go_{\ga} - \tau_0 = d\lambda_{\ga}.$$
Choosing a partition of unity $\psi_{\ga}:B \into [0,1]$ subordinate to the cover $(U_{\ga})$ we define
$$\tau = \tau_0 + d(\sum_{\ga} (\psi_{\ga}\circ \pi) \lambda_{\ga}).$$
This form is clearly in the same cohomology class as $\tau_0$ and further $d(\psi_{\ga}\circ \pi)$ vanishes on vectors tangent to the fibers, hence the restriction of $\tau$ to the fiber over $b\in B$ is
\begin{align*}
i^* \tau &= i^*\tau_0 +  \sum_{\ga}(\psi_{\ga}\circ \pi) i^*(d\lambda_{\ga})\\
&=i^*\tau_0 +  \sum_{\ga}(\psi_{\ga}\circ \pi) i^*(\go_{\ga}- \tau_0))\\
& = \sum_{\ga}(\psi_{\ga}\circ \pi) i^*\go_{\ga} = \go.
\end{align*}
We now claim that for \e\ small enough the forms $ e^{i \e\tau} \wedge \pi^*\rho_{\ga}$ patch together to give a \gcs\ on the total space. Indeed, from the above, $\tau$ is nondegenerate on the vertical subspaces $\Ker \pi_*$ and therefore it determines a field of horizontal subspaces
$$ Hor_x = \{X \in T_xM : \tau(X,Y) = 0, \forall Y \in T_xF\}.$$
The subspace $Hor_x$ is a complement to $T_x F$ and isomorphic to $T_{\pi(x)}B$ via $\pi_*$. Also, denoting by $(\cdot, \cdot)_B$ the Mukai pairing on $B$, $(\rho_{\ga},\overline{\rho_{\ga}})_B \neq 0$ and hence pulls back to a volume form on $Hor$. Therefore, for \e\ small enough,
$$(e^{i \e\tau} \wedge\pi^*\rho_{\ga}, e^{-i \e\tau} \wedge \pi^*\overline{\rho}_{\ga})_X  = (\e\tau)^{\dim(F)} \wedge \pi^*(\rho_{\ga},\overline{\rho_{\ga}})_B\neq 0,$$
and $e^{i \e\tau} \wedge\pi^*\rho_{\ga}$ is of the right algebraic type. Finally, from the integrability condition for $\rho_{\ga}$, there are $X_{\ga},\xi_\ga$ such that
$$d\rho_{\ga} = (X_{\ga} + \xi_{\ga} - H) \rho_{\ga}.$$
If we let $X_{\ga}^{hor}\in Hor$ be the horizontal vector projecting down to  $X_{\ga}$, then the following holds on $U_{\ga} \times F$:
$$d (e^{i\e\tau}\wedge\pi^* \rho_{\ga}) = (X_{\ga}^{hor} + \pi^*\xi_{\ga} - X_{\ga}^{hor}\lfloor i \e \tau - \pi^*H)e^{i\e\tau}\wedge\pi^*\rho_\ga,$$
showing that the induced \gcs\ is integrable.

On the overlap of $U_{\ga}$ and $U_{\beta}$, there is a function $g_{\ga}^{\beta}$ such that $\rho_{\ga} = g_{\alpha}^{\beta}\rho_{\beta}$, therefore the same applies to $e^{i\e\tau}\pi^*\rho_{\ga}$ and $e^{i\e\tau}\pi^*\rho_{\beta}$ on the overlap $U_{\ga} \cap U_{\beta} \times F$, hence these local forms patch together to give a line bundle.
\end{proof}

Several cases for when the conditions of the theorem are fulfilled have been studied for symplectic manifolds and many times purely topological conditions on the base and on the fiber are enough to ensure that the hypotheses hold. We give the following two examples adapted from McDuff and Salamon \cite{MS00}, Chapter 6.

\begin{theo}{surface bundles}
Let $\pi:X \into M$ be a symplectic fibration over a compact generalized complex base with fiber $(F,\go)$. If the first Chern class of $TF$ is a nonzero multiple of $[\go]$, then the conditions of Theorem \ref{T:Thurston} hold. In particular, any oriented surface bundle can be given a symplectic fibration structure and, if the fibers are not tori, the total space has a \gcs.
\end{theo}

\begin{theo}{simply connected}
A symplectic fibration with compact and 1-connected fiber and a compact twisted generalized complex base admits a twisted \gcs.
\end{theo}

These results also show that the case of torus bundles is `more interesting', in the sense that more structure has to come into play in order to determine whether there is a \gcs\ on the total space. In Chapter \ref{nilmanifolds} we will study \gcss\ on nilmanifolds, which are torus bundles over other nilmanifolds. In many cases we find these manifolds possess \gcss, though not necessarily compatible with the most obvious fibration.

\section{Further Examples: Lie Groups}\label{lie groups}

Of course, one can have a \gcs\ on the total space even if the conditions of the Theorem \ref{T:Thurston} fail. One such case is the `more interesting' one of principal torus bundles. As we have mentioned, the techniques from the previous section fail in this case and further structure has to come into play.

\begin{ex}{principal torus bundle}
Let $M^{2n}$ be a complex manifold satisfying the $dd^c$-lemma and let $X$ be the total space of a principal 2-torus bundle over $M$. Seeing $X$ as a double circle bundle, assume that the Chern classes of each of these are in $H^{1,1}(M)$. We claim that $X$ admits \gcss\ of type $2n+1$ (complex) and $2n$. Indeed, let $\gO$ be a local trivialization of $\wedge^{n,0}TM$ and $u_1, u_2$ be connection forms of the torus bundle. Then, as $du_i \wedge \Omega =0$, the form $e^{iu_1\wedge u_2}\gO$ determines locally a \gcs\ on $X$ of type $n$ (even though there may be no $a \in H^2(X)$ which restricts to $u_1\wedge u_2$ on the fibers) and the form $(u_1 + i u_2)\wedge \gO$ determines locally a type $n+1$ structure. It is clear that in either case these forms patch together yielding the canonical line bundles $e^{iu_1\wedge u_2}\wedge^{n,0}TM$ and $(u_1 + i u_2)\wedge^{n,0}TM$ respectively.

One concrete example where this happens is the product of two odd spheres $S^{2n_1 +1} \times S^{2n_2 +1}$ as those can be expressed as principal 2-torus bundles over $\C P^{n_1} \times \C P^{n_2}$. The type $n_1+n_2+1$ structure (complex) coincides with the Calabi-Eckmann complex structure on these manifolds and here we have established that $S^{2n_1 +1} \times S^{2n_2 +1}$ also admits a type $n_1 + n_2$ structure.

Of course we have used strongly the fact that the base has a complex structure. In the case, say, of a symplectic base these results wouldn't hold. For example, the case of principal 2-torus bundles over surfaces with genus greater than $1$ has been studied in the symplectic setting by Etg\"u \cite{Et01} and Walczak \cite{Walc04} and, unless the bundle is flat, there are no symplectic structures on $X$.
\end{ex}

\begin{ex}{liegroups as torus bundles}
The previous example can be extended to any principal $2k$-torus bundle, $E^{2n}$, whose Chern classes of the induced $2k$ circle bundles are of type $(1,1)$. Any choice of invariant generalized complex structure on the torus gives rise to a generalized complex structure on the total space, still using the complex structure in the base. This shows that the total space of such a bundle admits \gcss\ of any type between $n-k$ and $n$. The most typical examples of torus bundles with such a property are compact even dimensional Lie groups, which, after a choice of maximal torus, can be seen as principal torus bundles over the corresponding (K\"ahler) Flag manifold. Hence a compact Lie group $G^{2n}$ of rank $2k$ admits \gcss\ of any type between $n-k$ and $n$.
\end{ex}

This same result about \gcss\ on compact Lie groups can be obtained from a Lie algebra point of view, which is useful to prove further properties of these structures. The following argument is a generalization of the one used by Samelson to prove that every even dimensional Lie group admits left/right invariant complex structures \cite{Sam53}.

\begin{theo}{invariant structures on lie groups}
A compact Lie group $G^{2n}$ of rank $2k$ admits left invariant \gcss\ of any type between $n-k$ and $n$. Further these structures can be chosen to be bi-invariant under the action of a maximal torus.
\end{theo}
\begin{proof}
Choose a maximal torus $T^{2k} < G$ with Lie algebra $\frak{t} < \frak{g}$. Let $\{ \ga: \ga \in \frak{t}^* \}$ be the corresponding root system and $e_{\ga} \in \frak{g}\tensor \C$ be eigenvectors of $\ga$:
$$ [t, e_{\ga}] = \ga(t) e_{\ga}, \mbox{ for } t \in \frak{t},$$
so that we have a splitting $ \frak{g}\tensor \C = \frak{t}\tensor \C \oplus \mbox{span}\{e_{\ga}\}$.  This splitting of $\frak{g}\tensor \C$ allows us to identify $\frak{t}^*$ with the annihilator of $\mbox{span}\{e_{\ga}\}$ and hence we can see $\frak{t}\oplus \frak{t}^*$ as a subspace of $\frak{g}\oplus \frak{g}^*$.

If we let $L_{\frak{t}} < (\frak{t}\oplus \frak{t}^*)\tensor \C$ be the $+i$-eigenspace of an invariant \gcs\ on the maximal torus, fix a set of positive roots and denote the duals by upper indexes, then we claim that
$$ L = \mbox{span}\{e_{-\ga},e^{\ga},L_{\frak{t}}:\ga >0\}$$
induces a \gcs\ on the Lie group of type $n-k + type(\J_{\frak{t}})$, i.e., $n-k$ plus the type of the \gcs\ determined by $L_{\frak{t}}$ in $T$.

Indeed, it is easy to see that $L \cap \overline{L} = \{0\}$, $\dim_{\C} L = 2n$ and the type is $n-k + type(\J_{\frak{t}})$. To prove that $L$ is closed under the Courant bracket, we use expression \eqref{E:Courant bracket2} for the Courant bracket, which in this case becomes
$$[X+\xi,Y + \eta] = [X,Y] + X\lfloor d\eta - Y \lfloor d\xi.$$
With this, if $v_i \in \frak{t}$, $w \in \frak{t}^*$ and $e_{\ga}$ and $e^{\ga}$ are as above, we have $[\frak{g}^*,\frak{g}^*] =0$ and also
\begin{center}
\noindent
\begin{tabular}{l l l l l l}
$[v_1,v_2] $&$=0$  &$ [v_1,w] $&$= 0 $& $[v_1,e_{\ga}]$ &$ = \ga(v_1) e_{\ga}$ \\
$ [v_1,e^{\ga}] $&=$ -\ga(v_1)e^{\ga}$&$ [w, {e_\ga}]$&$= \lambda_1 e^{-\ga}$& $[e_{\ga},e_{\gb}]$ & $= \lambda_2 e_{\ga+ \gb}$\\
&&$[e_{\ga},e^{\gb}]$ & $=\lambda_3 e^{\gb-\ga}, $& & 
\end{tabular}
\end{center}
where the $\lambda_i$ are constants. With the Courant bracket established, one can check by inspection that $L$ is involutive.

If $\rho_{\frak{t}}$ is the form in $\wedge^{\bullet}\frak{t}^*$ determining $L_{\frak{t}}$, then any invariant differential form determining $L$ is a multiple of
\begin{equation}\label{E:pure form for lie groups}
\rho = \rho_{\frak{t}} \wedge_{\ga>0}e^{\ga}.
\end{equation}

Bi-invariance of \J\  under the action of the maximal torus is a consequence of the following fact
$$ \J[v, t] = [\J v,t], \mbox{ for } t \in \frak{t},$$
which can also be proved by inspection, as the bracket has already been determined.
\end{proof}

As the canonical bundle of any invariant \gcs\ on a Lie group is trivial (trivialized by an invariant form), one might wonder when these structures are actually \gcy, i.e., when the canonical bundle admits a nonvanishing closed section. The following answers this for the structures established in Theorem \ref{T:invariant structures on lie groups}.

\begin{theo}{gy and lie groups}
If $G$ is not abelian, none of the \gcss\ determined in Theorem \ref{T:invariant structures on lie groups} is determined by an invariant closed form, i.e., the structure is not invariant \gcy.
\end{theo}
\begin{proof}
We shall start by showing that the form from equation \eqref{E:pure form for lie groups} is not closed. And, in order to do so, we shall first establish that
$$(d \theta) \wedge_{\ga>0} e^{\ga} = 0, \qquad \mbox{for all } \theta \in \frak{t}^*\tensor \C.$$
Indeed, if $v \in \frak{t}\tensor \C$ and $e_{\ga}$ is as above, then, from $d\theta (X,Y) = -\theta([X,Y])$, we get $d\theta(v,\cdot) = 0$  and
$$ d\theta(e_{\ga},e_{\gb}) = \lambda \theta(e_{\ga+\gb}) =0 \qquad \mbox{if } \gb \neq -\ga,$$
where $\lambda$ is a constant. Hence $dt \in \mbox{span}\{e^{\ga}\wedge e^{-\ga}:\ga >0\}$ and the claim follows.

Now, for $\rho$ defined in \eqref{E:pure form for lie groups}, the lower order term is
$$ \rho_0 = \wedge_{\ga>0}e^{\ga}\wedge_{1\leq j \leq j_0}\theta^j,$$
where $\theta^j \in \frak{t}^*\tensor \C$. If $\rho$ determines a \gcy\ structure, this form has to be closed and we can extended the $\theta^j$ so as to form a basis of $(1,0)$ vectors in $\frak{t}^*\tensor \C$ for some complex structure. From our initial claim, as $\rho_0$ is closed, so is 
$$\rho_0 \wedge_{j_0 < j \leq k}\theta^j = \wedge_{\ga>0}e^{\ga}\wedge_{1\leq j \leq k}\theta^j.$$
If we let $t \in \frak{t}$ be an element determining the set of positive roots (i.e. $\ga >0$ if and only if $\ga(t)>0$) and split it into its $(1,0)$ and $(0,1)$ parts and let $\theta_j$ be the $1,0$ vectors dual to $\theta^j$ then
\begin{align*}
0&=\mbox{Re}(d(\rho_0 \wedge_{j_0 < j \leq k}\theta^j)(\wedge_{\ga >0}e_{\ga}\wedge_{1 \leq j \leq k}\theta_j\wedge t^{0,1})) =\\
&= \mbox{Re}(-\sum_{\ga>0} \ga(t^{0,1})\wedge_{\ga>0} e^{\ga} \wedge_{1 < j \leq k}\theta^j (\wedge_{\ga >0}e_{\ga}\wedge_{1 \leq j \leq k}\theta_j\wedge t^{0,1}))\\
&= -\mbox{Re}(\sum_{\ga>0} \ga(t^{0,1})) = - \frac{1}{2}\sum_{\ga >0} \ga(t) < 0,
\end{align*}
which is a contradiction, hence $\rho$ can not be closed.
\end{proof}

Other structures that will be important on semi-simple Lie groups are twisted \gcss. This is because these will give us examples of twisted \gk\ manifolds (see Example \ref{T:liegroups}). One question that might be raised is whether compact semi-simple Lie groups can be twisted generalized Calabi--Yau, but the answer is no.

\begin{theo}{liegroups3}
No compact semi-simple Lie group admits a twisted \gcy\ structure with nontrivial twist.
\end{theo}
{\sc Proof.}
From Example \ref{T:liegroups dh-cohomology}, for a semi-simple Lie group the $d_H$-cohomology is trivial for any nonvanishing class  $H \in H^3(G)$. Therefore, if there was a globally defined $d_H$-closed form $\rho$ defining the structure, this form would be $d_H$-exact. Then, the claim is a consequence of the following lemma.
\vskip6pt
\begin{lem}\label{T:rho cannot be exact}
If a $d_H$-closed form $\rho$ defines a \gcy\ structure on a compact manifold, the $d_H$-cohomology class of $\rho$ is nonzero.
\end{lem}
\begin{proof}
Using \eqref{E:mukai pairing}, one can easily show that the Mukai pairing of two $d_H$-exact forms is $d$-exact and hence can not be a nontrivial element in the top degree cohomology. On the other hand, from Proposition \ref{T:chevalley},
$$(\rho,\overline{\rho})  \neq 0,$$
is a volume form, so the $d_H$-cohomology class of $\rho$ is nonzero.
\end{proof}

Having solved the twisted case and given that none of the \gcss\ obtained in Example \ref{T:invariant structures on lie groups} is \gcy, we are led to the following.

\begin{conjecture}
Compact semi-simple Lie groups do not admit (invariant) \gcy\ structures.
\end{conjecture}

While answering the conjecture in the invariant setting may be easy using Cartan's classification of semi-simple Lie algebras (and probably the structures above already account for all possible invariant \gcss), in the noninvariant case this is a very difficult question. Even the simple case of $S^3 \times S^3$ is still an open problem and the best results in this direction come from Lu and Tian \cite{LuTi94}, showing that the connected sum of $m$ copies of $S^3 \times S^3$ have holomorphically trivial canonical bundle for $m >1$. 

\section{Deformations of Generalized Complex Structures}

In this section we state part of Gualtieri's deformation theorem for \gcss. The space of infinitesimal deformations is naturally a subspace of the space of sections of $\wedge^2\overline{L}$ and we want to know which sections of $\wedge^2\overline{L}$ give rise to deformations of the \gcs. We shall not discuss when such deformations are trivial and instead refer to Gualtieri's thesis. Before we can state the theorem we need one extra background lemma.

Initially we observe that the Courant bracket does not satisfy the Jacobi identity. Indeed, its Jacobiator
$$\mathrm{Jac}(A,B,C)= [[A,B],C] + [[B,C],A]+[[C,A],B]$$
is given by the following expression:
\begin{lem}(Liu, Weinstein and Xu \cite{LWX97}):
The Jacobiator of the Courant bracket is given by
$$\mathrm{Jac}(A,B,C) = \frac{1}{3}d(\langle [A,B],C \rangle+\langle [B,C],A \rangle+\langle [C,A],B \rangle),$$
where, $\langle \cdot, \cdot \rangle$ is the natural pairing on $TM\oplus T^*M$.
In particular, if a subbundle $L$ determines a \gcs, then the restriction of the Courant bracket to $L$ is a Lie bracket.
\end{lem}

Since the Courant bracket is a Lie bracket in $L$ it induces an exterior derivative in $\wedge^{\bullet}L^*$: $d_L$. Also, as $L$ is isotropic, the natural pairing gives a natural isomorphism  $L^* \cong \overline{L}$ and hence $d_L$ is a derivative of sections of $\wedge^{\bullet}\overline{L}$. Drawing on a result of Liu, Weinstein and Xu \cite{LWX97}, Gualtieri established the following deformation theorem.

\begin{theo}{deformation}
{\em (Gualtieri \cite{Gu03}, Theorem 5.4):} An element $\e \in \wedge^2 \overline{L}$ gives rise to a deformation of \gcss\ if and only if $\e$ is small enough and satisfies the Maurer--Cartan equation
$$d_L \e + \frac{1}{2}[\e,\e] =0.$$
The deformed complex structure is given by
$$ L^{\e} = (\Id + \e)L \qquad \overline{L}^{\e} = (\Id + \overline{\e})\overline{L},$$
or, on forms,
$$ \rho_{\e} = e^{\e} \cdot \rho.$$
In the complex case, a bivector $\e\in \wedge^{2,0}TM < \wedge^2 \overline{L}$ gives rise to a deformation only if each of the summands vanish, i.e., $\delbar\e=0$ (\e\ is holomorphic) and $[\e,\e]=0$ (\e\ is Poisson).
\end{theo}

\begin{ex}{surfaces}
Any holomorphic bivector \e\ on a complex surface $M$ is also Poisson, as $[\e,\e] \in \wedge^{3,0}TM = \{0\}$, and hence gives rise to a deformation of \gcss. The deformed \gcs\ will be symplectic outside the divisor representing $c_1(M)$ where the bivector vanishes. At the points where $\e=0$ the deformed structure agrees with the original complex structure. 
\end{ex}

\begin{ex}{HK}
Let $M^{4n}$ be a hyperk\"ahler manifold with K\"ahler forms $\go_I, \go_J$ and $\go_K$. According to the K\"ahler structure $(\go_I, I)$, $(\go_J + i\go_K)$ is a closed holomorphic 2-form  and $(\go_J + i\go_K)^{2n}$ is a holomorphic volume form. Therefore these generate a holomorphic Poisson bivector $\Lambda \in \wedge^{2,0}TM$ by
$$\Lambda \cdot (\go_J + i \go_K)^{2n} = 2n(\go_J + i \go_K)^{2n-1}.$$
The deformation of the complex structure $I$ by $t \Lambda$ is given by
$$e^{t \Lambda}(\go_J + i \go_K)^{2n} = t^{2n} e^{\frac{\go_J + i \go_K}{t}}.$$
which interpolates between the complex structure $I$ and the $B$-field transform of the symplectic structure $\go_K$ as $t$ varies from 0 to 1.

We also remark that although in this case we could explicitly find the path connecting the deformed structures, this is peculiar to the complex case and in general it is hard to find a path interpolating between two given \gcss: indeed they may have different topological invariants.
\end{ex}

\begin{ex}{fibrations}
Let $X \into B$ be a symplectic fibration over a twisted generalized complex base $B$ satisfying the hypothesis of Theorem \ref{T:Thurston}. Then any deformation of the \gcs\ on $B$ gives rise to a deformation of the \gcs\ on $X$. Indeed, as seen in the proof of the theorem, for a given \gcs\ with canonical line bundle $K$, $e^{i\e\tau}\wedge \pi^*K$ is the  canonical line bundle of the fibration. If $\tilde{K}$ is the canonical line bundle of a nearby \gcs, we can still use the same \e\ and $e^{i\e\tau}\wedge \pi^*\tilde{K}$ will determine a deformed structure on $X$.
\end{ex}

\section{Submanifolds}\label{submanifolds}

In some sense, a generalized complex structure is actually a twisted \gcs\ with zero twist, as even though the 3-form is zero, it may be thought to exist and be part of the structure. This is most evident when one considers submanifolds. At first, the natural definition seems to be that a submanifold $M \hookrightarrow N$ is a generalized complex submanifold if
$$\tau_0 = TM \oplus T^0M = \{X + \xi \in TM \oplus T^*N : \xi|_{TM}=0\}$$
is invariant under \J. The problem with this definition is that, as remarked in Example \ref{T:B-field2}, the natural concept of equivalence for generalized complex manifolds is given by the actions of diffeomorphisms and $B$-fields. Therefore a  meaningful definition of submanifold has to be invariant under those actions and the above definition does not support $B$-field actions.

Gualtieri modifies the above definition of generalized complex submanifold by considering the twisted and nontwisted cases all at once.

\begin{definition}
A {\it submanifold with 2-form} $(M,F) \hookrightarrow (N,H)$ is a submanifold of a manifold $N$ with 3-form if $dF = H|_M$. If $(N,H,\J)$ is a twisted generalized complex manifold, then a {\it generalized complex submanifold} is a submanifold with 2-form $(M,F)$ such that the space
$$\tau_F = (Id + F) \cdot TM\oplus T^0M = \{ X + \xi \in TM \oplus T^*N : \xi|_{TM}=X\lfloor F\}$$
is invariant under \J.
\end{definition}

This way we see that if $(M,F) \hookrightarrow (N,H,\J)$ is a generalized complex submanifold, then $(M,F-B)$ will be a submanifold of the $B$-field transform of $(N,\J)$ which shows this definition is well behaved under $B$-field tranforms.

\begin{ex}{complex submanifold} (Gualtieri \cite{Gu03}, Example 7.7): A generalized complex submanifold $(M,F) \hookrightarrow (N,0)$ of a complex manifold is the same thing as a complex submanifold, as in this case $\J$ maps $TM$ to itself. Since $M$ is complex, the conormal bundle of $M$ is also invariant under $\J$ and therefore the invariance of $\tau_F$ under the complex structure implies that $F \in \Omega^{1,1}M$. Hence a generalized complex submanifold of a complex manifold is just a complex submanifold with a closed $(1,1)$-form.
\end{ex}

\begin{ex}{A-branes} (Gualtieri \cite{Gu03}, Example 7.8):
The symplectic case is more interesting and, as Gualtieri observed, \gc\ submanifolds coincide with coisotropic A-branes introduced by Kapustin and Orlov \cite{KaOr03}. We claim that if $(M,F)$ is a \gc\ submanifold of a symplectic manifold $(N,0)$, then $M$ is a coisotropic submanifold. Also, both $F$ and, obviously, $\go|_{M}$ are annihilated by the distribution
$$TM^{\go} = \{X \in TN|_{M} : \go(X, Y) = 0,\quad \forall Y \in TM\}.$$
In the quotient $V = TM/TM^{\go}$, there is a complex structure induced by $\go^{-1}F$ and $(F + i \go)$ is a $(2,0)$-form whose top power is a volume element in $\wedge^{k,0}V$. Finally, if $F=0$, then $M$ is just a Lagrangian submanifold of $X$.

To prove that $M$ is coisotropic, we take $X\in TM^{\go}$, so that $\go(X,\cdot) \in T^0M$, and indeed any element in $T^0M$ is of that form. Due to the invariance of $\tau_F$, for $\go(X,\cdot) \in T^0 M < \tau_F$, we have 
$$\J\go(X,\cdot) = -\go^{-1} \go(X, \cdot)  = - X \in \tau_F.$$
This furnishes both that $X \in TM$, and hence $M$ is coisotropic, and that $F(X,\cdot)=0$, and hence $F$ is also annihilated by the distribution $TM^{\go}$.

To find the complex structure on the quotient $TM/TM^{\go}$ we take $X \in TM$ and apply \J\ to $X + \xi \in \tau_F$. Invariance implies that
$$-\go^{-1}\xi + \go(X,\cdot) \in \tau_F,$$
which, is the same as $-F\circ\go^{-1}\circ F (X)|_{M} = \go(X,\cdot)|_{M}$ where we  are seeing $F$ as a map from $TM$ to $T^*M$. In $TM/TM^{\go}$, there is an inverse $\go^{-1}$ and hence, in the quotient, we have the identity
$$- X = (\go^{-1}F)^2(X),$$
showing that $\go^{-1}F$ induces a complex structure on $TM/TM^{\go}$.

For an $X \in \wedge^{0,1}(TM/TM^{\go})$ we have $\go^{-1}F(X) = -iX$. Applying \go\ we get $F(X,\cdot)+i\go(X,\cdot) =0$ and hence $F + i\go$ is annihilated by $\wedge^{0,1}(TM/TM^{\go})$ and thus is a $(2,0)$-form. Finally, for $X = X_1 + i X_2 \in \wedge^{1,0}(TM/TM^{\go})$, as before we obtain
$(F - i\go)(X,\cdot) =0$, which spells out as
$$F(X_1,\cdot) = -\go(X_2,\cdot) \qquad \mbox{ and } \qquad F(X_2,\cdot) = \go(X_1,\cdot),$$
and therefore $(F+i\go)(X,\cdot)  = -2\go(X_2,\cdot) + 2 \go(X_1,\cdot) \neq 0$, as \go\ is nondegenerate in $TM/TM^{\go}$. Thus $F+i \go$ is a nondegenerate 2,0-form.

If $F$ vanishes, $0$ is a complex structure in $TM/TM^{\go}$ which must therefore be the trivial vector space and hence $M$ is Lagrangian.
\end{ex}

\section{Generalized K\"ahler Manifolds}

The key to generalize the concept of K\"ahler manifold is the observation that for such a manifold there are two \gcss, $\J_1, \J_2$ --- one coming from the complex structure and one from the symplectic --- which commute, as the symplectic form $\go$ is of type $(1,1)$, and such that
$$ \langle \J_1\J_2 v ,v \rangle > 0, \quad \mbox{ for } v\neq 0.$$
Further, a K\"ahler structure is a Calabi-Yau structure if
$$(e^{i\go},e^{-i\go})= \frac{(-2i\go)^n}{n!} = k\gO \wedge \overline{\gO} = (-1)^{|\gO|} k (\gO,\overline{\gO}),$$
for some constant $k \in \C$.

\begin{definition}
A {\it \gk\ structure} on a manifold $M$ is a pair of commuting \gcss, $\J_1$ and $\J_2$ such that
$$ \langle \J_1 \J_2 \cdot, \cdot \rangle$$
is a positive definite metric on $TM\oplus T^*M$.

A {\it generalized Calabi-Yau metric structure} is a \gk\ structure determined by two closed forms $\rho_1$ and $\rho_2$ whose Mukai pairings are related by
$$(\rho_1, \overline{\rho_1}) = k (\rho_2,\overline{\rho_2}),$$
for some constant $k$.
\end{definition}

Letting $G = \J_1\J_2$, for a \gk\ structure, it is clear that $G$ is orthogonal, as a composition of orthogonal maps. Further, $G^2 = \Id$, showing that $G$ is self adjoint and the \gk\ metric is compatible with the natural pairing. Therefore $G$ determines a decomposition of $TM\oplus T^*M$ into $\pm 1$-eigenspaces $C_{\pm}$. As $\J_1$ and $\J_2$ commute, they can be diagonalized simultaneously. If we denote their $+i$-eigenspaces by $L_{1/2}$, this means that
$$ L_1 = L_1 \cap L_2 \oplus L_1 \cap \overline{L_2}\qquad L_2 = L_1 \cap L_2 \oplus \overline{L_1} \cap L_2.$$
Hence
$$C_+ \tensor \C = L_1\cap \overline{L_2} \oplus \overline{L_1} \cap L_2,$$
$$C_- \tensor \C = L_1\cap L_2 \oplus \overline{L_1} \oplus  \overline{L_2}.$$

This shows that $L_1 \cap L_2$ and $L_1 \cap \overline{L_2}$ are $n$-complex-dimensional spaces.

Moreover, according to work on Section \ref{metric}, $C_{\pm}$ are the graphs of $b\pm g$ over $TM$ and, $\J_1$ can be used as a complex structure on $C_{\pm}$ orthogonal with respect to the natural pairing. Hence, projecting to $TM$, we obtain two almost complex structures on $TM$, $J_{\pm}$, orthogonal with respect to $g$: one coming from $C_+$ the other coming from $C_-$. But the $+i$-eigenspaces for the complex structures in $C_{\pm}$ are  $L_1 \cap \overline{L_2}$ and $L_1 \cap L_2$, which are integrable with respect to the Courant bracket, hence their projections over the tangent bundle are integrable with respect to the Lie bracket and $J_{\pm}$ are {\it integrable}. This shows that a generalized K\"ahler manifold has a bihermitian structure.

Finally, Gualtieri shows that if we let $\go_{\pm}$ be the 2-forms associated with the hermitian structures:
$$\go_{\pm} = g(J_{\pm}\cdot, \cdot),$$
then
$$d^c_{+}\go_+ = - d^c_{-}\go_- = db,$$
Where $d^c_{\pm}$ is the $d^c$ operator for the complex structure $J_{\pm}$. And indeed all this structure is enough to describe a \gk\ manifold.

\begin{theo}{proposition5.17}
{\em (Gualtieri \cite{Gu03}, Proposition 5.17)} A bihermitian manifold $(M,g,J_{\pm})$ admits a \gk\ structure for which $J_{\pm}$ are the complex structures induced by the \gk\ structure if and only if 
$$d^c_{+}\go_+ = - d^c_{-}\go_- = db,$$
for some 2-form $b$.

In this case, the generalized complex structures $\J_1$ and $\J_2$ are given by
\begin{equation}\label{E:gb to J}
\J_{1/2} = \frac{1}{2}
\begin{pmatrix} 1 & 0 \\ b & 1\end{pmatrix}
\begin{pmatrix}
J_+ \pm J_- & -(\go_+^{-1} \mp \go_-^{-1})\\
\go_+ \mp \go_- & -(J_+^* \pm J_-^*)
\end{pmatrix}
\begin{pmatrix}
1&0\\-b & 1
\end{pmatrix}.
\end{equation}
\end{theo}
\begin{remark}
For a generic bihermitian structure, $d^c_{+}\go_+$ may not even be closed, hence the \gk\ condition is that not only $d^c_{+}\go_+ = - d^c_{-}\go_-$, but also that this is an exact form.
\end{remark}

\begin{ex}{HK to GK} (Gualtieri \cite{Gu03}, Example 6.30) Let $(M,I,J,K,g)$ be a hyperk\"ahler manifold. Then $g$ is hermitian with respect to the three complex structures, and the respective K\"ahler forms, $\go_I$, $\go_J$ and $\go_K$ are $d^c$-closed. Therefore, choosing, say, $I$ and $J$ for $J_{\pm}$ we see that the hypotheses of the above theorem are fulfilled with $b=0$, hence these two complex structures furnish the hyperk\"ahler manifold with a generalized K\"ahler structure. The generalized complex structures $\J_{1/2}$ are given by
$$\J_{1/2} = \frac{1}{2}\begin{pmatrix}
I \pm J & -(\go_I^{-1} \mp \go_J^{-1})\\
\go_I \mp \go_J & - (I^* \pm J^*).
\end{pmatrix} $$
Recalling that the complex structures $I, J$ and $K$ satisfy relations like $I = -\go_J^{-1}\go_K$, we can re-write $\J_{1/2}$ as:
$$ \J_1 = \begin{pmatrix}1 & 0\\ \go_K & 1 \end{pmatrix}
\begin{pmatrix}
0 & -\frac{1}{2}(\go_I^{-1} - \go_J^{-1})\\
\frac{1}{2}(\go_I - \go_J) & 0
\end{pmatrix}
\begin{pmatrix}1 & 0\\ -\go_K & 1
\end{pmatrix}
$$
$$ \J_2 = \begin{pmatrix}1 & 0\\ -\go_K & 1 \end{pmatrix}
\begin{pmatrix}
0 & -\frac{1}{2}(\go_I^{-1} + \go_J^{-1})\\
\frac{1}{2}(\go_I + \go_J) & 0
\end{pmatrix}
\begin{pmatrix}1 & 0\\ \go_K & 1
\end{pmatrix}.
$$
These are $B$-field transforms of symplectic structures, and the respective differential forms defining the structures are:
$$\rho_1 = \exp(\go_K + \frac{i}{2}(\go_I - \go_J));$$
$$\rho_2 = \exp(-\go_K + \frac{i}{2}(\go_I + \go_J)).$$
\end{ex}
This example shows that it is possible to have two symplectic structures giving rise to a \gk. The same is not true for complex structures as Gualtieri shows that the sum of the types of $\J_1$ and $\J_2$ can not exceed $n$, for a $2n$-dimensional manifold. Also, this is an example of a \gk\ structure which is not merely a $B$-field transform of an actual K\"ahler structure. However, one can still argue that, although the structure is different, the manifold does admit a K\"ahler structure.

\vskip6pt
{\sc Open problem}: {\it Find examples of \gk\ manifolds that do not admit K\"ahler structures}.
\vskip6pt

Another interesting problem in this area concerns the deformation of \gk\ structures. It is a result of Kodaira and Spencer that the K\"ahler condition is open \cite{KoSp58,KoSp60}, in the sense that if one deforms a complex structure compatible with a K\"ahler metric, it is possible to find another K\"ahler  metric with which the deformed complex structure compatible, as long as the deformation is small enough.

\begin{conj}{gk is open}
The \gk\ condition is open.
\end{conj}

\begin{remark}\label{T:upq}
One last and important remark concerns the decomposition of forms in the \gk\ case. Each \gcs\ gives rise to a decomposition of $\wedge^{\bullet}TM \tensor \C$ into $V_j^k$  the $ik$-eigenspaces of $\J_j$. As $\J_{1}$ and $\J_2$ commute, so do their Lie algebra actions. Hence their action on forms can be simultaneously diagonalized. If we let $U^{p,q} = V_1^p \cap V_2^q$, then
$$ V_1^p = \sum_q U^{p,q},\qquad \mbox{and} \qquad V_2^q = \sum_p U^{p,q}.$$
Further, as $L_1\cap L_2$ is $n$-complex-dimensional, we get that $V_1^{\pm n} \in V_2^0$ and vice-versa. And in general, from there,
$$V_1^{\pm (n-k)} \subset \oplus_{j=0}^{k}V_2^{-k+2j}.$$
\end{remark}

\subsection{Twisted Generalized K\"ahler Manifolds}

By considering twisted structures, we obtain the twisted version of the \gk\ condition.

\begin{definition}
An {\it $H$-twisted \gk\ structure} on a manifold $M$ is a pair of commuting $H$-twisted \gcss, $\J_1$ and $\J_2$ such that
$$ \langle \J_1 \J_2 \cdot, \cdot \rangle$$
is a definite metric on $TM\oplus T^*M$.
\end{definition}

These structures can also be described in bihermitian terms (in the twisted case, $J_{\pm}$ will also be integrable complex structures).
\begin{theo}{proposition5.17twisted}
{\em (Gualtieri \cite{Gu03}, Theorem 5.37)} A bihermitian manifold $(M,g,J_{\pm})$ admits an $H$-twisted \gk\ structure for which $J_{\pm}$ are the complex structures induced by the twisted \gk\ structure if and only if 
$$d^c_{+}\go_+ = - d^c_{-}\go_- = H + db,$$
for some 2-form $b$.

In this case, the twisted generalized complex structures $\J_1$ and $\J_2$ are given by Equation \eqref{E:gb to J}.
\end{theo}

\begin{remark}
If $(M,g,J)$ is a hermitian manifold, we can consider the associated 2-form, $\go = g(J \cdot, \cdot)$. This structure is called {\it strong K\"ahler with torsion (SKT)} if $d^c \go$ is a nonzero closed form. By the above, a twisted \gk\ manifold is SKT with respect to both induced complex structures.
\end{remark}

\begin{ex}{liegroups}
(Gualtieri \cite{Gu03}, Example 6.39)
Let 
$J_L$ and $J_R$ be left and right invariant complex structures as constructed in Theorem \ref{T:invariant structures on lie groups} on an even dimensional, compact Lie group $G$. If $G$ is semi-simple, these can be chosen to be Hermitian with respect to the invariant metric induced by the Killing form $K$. We claim that this bihermitian structure furnishes $G$ with a twisted generalized K\"ahler structure with twist $H$, the bi-invariant Cartan 3-form: $H(X,Y,Z) = K([X,Y], Z)$. To prove this claim we just have to compute $d^c_L \go_L$:
\begin{align*}
A = d^c_{L} \go_L(X,Y,Z) &= d_L \go_{L}(J_LX, J_L Y, J_L Z)\\
 & = - \go_L([J_L X, J_L Y], J_L Z) + cp.\\
& = - K([J_L X, J_L Y], Z)  + cp.\\
& = - K(J_L [J_L X, Y] + J_L[X,J_L Y] + [X,Y],Z) + cp.\\
& = (2 K([J_L X,J_L Y], Z)+ cp.) - 3 H(X,Y,Z)\\
& = -2A - 3 H,
\end{align*}
where $cp.$ stands for cyclic permutations. This proves that $d^c_L \go_L = -H$. Since the right Lie algebra is antiholomorphic to the left, the same calculation yields $d^c_R \go_R = H$ and by Theorem \ref{T:proposition5.17twisted}, this bihermitian structure induces a twisted \gk\ structure.

Observe that the \tgks\ obtained in this example is neither left nor right  invariant under the Lie group action, as each of the induced complex structures entering formula \eqref{E:gb to J} is only left/right invariant. However, as both complex structures are invariant under the maximal torus action, the induced \tgk\ structure is also invariant under this action.
\end{ex}


\chapter{Generalized Complex Structures on Nilmanifolds}\label{nilmanifolds}

When Thurston gave the first example of a symplectic manifold with no K\"ahler structure in \cite{Th76}, he laid the ideas of what became known as symplectic fibrations and we have already seen in Theorem \ref{T:Thurston} that this construction also works in the generalized complex setting. But Thurston's simplest example, a principal torus bundle over a torus, can be generalized in a different direction to what is called a nilmanifold. Later work by Cordero, Fern\'andez and Gray \cite{CFG86,CFG91} brought these manifolds to the attention of differential geometers and the study of invariant geometries on these spaces has been an interesting source of examples \cite{BG88,FIL94,IRTU01,KS03,Sa01}. 

A nilmanifold is a homogeneous space $M=\Gamma\backslash G$, where
$G$ is a simply connected nilpotent real Lie group and $\Gamma$ is
a lattice of maximal rank in $G$.
Such groups $G$ of dimension $\leq 7$ have been classified, and 6
is the highest dimension where there are finitely many. According
to \cite{Mag86,Mo58} there are 34 isomorphism classes of connected,
simply-connected 6-dimensional nilpotent Lie groups. This means
that, with respect to invariant geometry, there are essentially 34
separate cases to investigate.

The question of which 6-dimensional nilmanifolds admit 
a symplectic structure was settled by Goze and
Khakimdjanov \cite{GoKh96}: exactly 26 of the 34 classes admit
 symplectic forms.  Subsequently, the question of
left-invariant complex geometry was solved by Salamon~\cite{Sa01}:
he proved that exactly 18 of the 34 classes admit an invariant
complex structure.  While the torus is the only nilmanifold
admitting a K\"ahler structure, 15 of the 34 nilmanifolds admit both
complex and symplectic structures.  This leaves us with 5 classes
of 6-dimensional nilmanifolds admitting neither complex nor
symplectic left-invariant geometry.  See Figure 1 for
illustration.

It was this result of Salamon which inspired Gualtieri and myself to ask whether the
5 outlying classes might admit a \emph{generalized complex
structure}.  The main result of this Chapter
is to answer this question in the affirmative: all 6-dimensional
nilmanifolds admit generalized complex structures.

\setlength{\unitlength}{1mm}
\begin{picture}(60,51)(-65,-26)
\linethickness{1pt}
\thinlines
\qbezier(10,15)(25,15)(25,0)
\qbezier(25,0)(25,-15)(10,-15)
\qbezier(10,-15)(-5,-15)(-5,0)
\qbezier(-5,0)(-5,15)(10,15)
\qbezier(-10,15)(-25,15)(-25,0)
\qbezier(-25,0)(-25,-15)(-10,-15)
\qbezier(-10,-15)(5,-15)(5,0)
\qbezier(5,0)(5,15)(-10,15)
\put(13,10){\makebox(0,0){\footnotesize Symplectic}}
\put(15,0){\makebox(0,0){$11$}}
\put(-15,0){\makebox(0,0){$3$}}
\put(0,0){\makebox(0,0){$15$}}
\put(0,-18){\makebox(0,0){$5$}}
\put(-11,10){\makebox(0,0){\footnotesize Complex}}
\put(-30,-25){\framebox(60,45)[bc]{{\footnotesize Generalized complex}}}
\end{picture}
\begin{center}
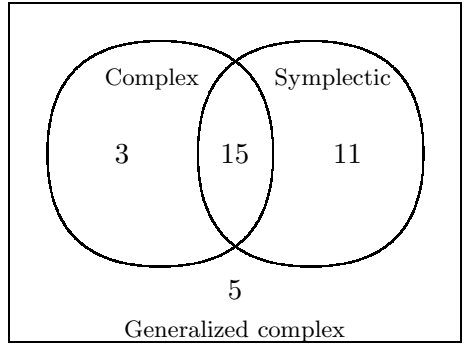
{\small {\it Figure 1: Left-invariant structures on the 34
six-dimensional nilpotent Lie groups.}}
\end{center}
\vskip12pt

We begin in Section~\ref{nil} with a brief introduction to nilmanifolds. Some results about generalized complex
structures on nilmanifolds in arbitrary dimension appear in
Section~\ref{gen}. Section~\ref{six} contains our main result: the
classification of left-invariant generalized complex structures on
6-dimensional nilmanifolds. In Section~\ref{iwa} we show that
while the moduli space of left-invariant complex structures on the
Iwasawa nilmanifold is disconnected (as shown in~\cite{KS03}), its
components can be joined using generalized complex structures. In
the final section, we provide an 8-dimensional nilmanifold which
does not admit a left-invariant generalized complex structure,
thus precluding the possibility that all nilmanifolds admit
left-invariant generalized complex geometry.

Except for Example \ref{T:type 2 to symplectic}, all the results on this chapter are present in previously advertised work \cite{CG04a}, in collaboration with Gualtieri.


\section{Nilmanifolds}\label{nil}

A nilmanifold is the quotient $M=\Gamma\backslash G$ of a
connected, simply-connected nilpotent real Lie group $G$ by the
left action of a maximal lattice $\Gamma$, i.e. a discrete
cocompact subgroup. By results of Malcev \cite{Mal51}, a nilpotent
Lie group admits such a lattice if and only if its Lie algebra has
rational structure constants in some basis.  Moreover, any two
nilmanifolds of $G$ rationally compatible with a lattice can be expressed as finite covers of a third
one.

A connected, simply-connected nilpotent Lie group is diffeomorphic
to its Lie algebra via the exponential map and so is contractible.
For this reason, the homotopy groups $\pi_k$ of nilmanifolds
vanish for $k>1$, i.e. nilmanifolds are Eilenberg-MacLane spaces
$K(\Gamma,1)$.  In fact, their diffeomorphism type is determined
by their fundamental group.  Malcev showed that this fundamental
group is a finitely generated nilpotent group with no element of
finite order.  Such groups can be expressed as central $\Z$
extensions of groups of the same type, which implies that any
nilmanifold can be expressed as a circle bundle over a nilmanifold
of lower dimension.  Because of this, one may easily use Gysin
sequences to compute the cohomology ring of any nilmanifold.
Nomizu used this fact to show that the rational cohomology of a
nilmanifold is captured by the subcomplex of the de Rham complex
$\Omega^\bullet(M)$ consisting of forms descending from
left-invariant forms on $G$:

\begin{theo}{Nomizu}
{\em (Nomizu \cite{No54})} The de Rham complex $\Omega^\bullet(M)$
of a nilmanifold $M=\Gamma\backslash G$ is quasi-isomorphic to the
complex $\wedge^{\bullet}\Gg^*$ of left-invariant forms on $G$,
and hence the de Rham cohomology of $M$ is isomorphic to the Lie
algebra cohomology of~$\Gg$.
\end{theo}

In this chapter we will search for generalized complex structures on
$\Gamma\backslash G$ which descend from left-invariant ones on
$G$, which we will call \emph{left-invariant} generalized complex
structures. This will require detailed knowledge of the structure
of the Lie algebra $\Gg$, and so we outline its main properties in
the remainder of this section.

Nilpotency implies that the central descending series of ideals
defined by $\Gg^0=\Gg,\ \Gg^k = [\Gg^{k-1},\Gg]$ reaches $\Gg^s=0$
in a finite number $s$ of steps, called the {\it nilpotency index}, $\nil(\Gg)$  (also called the nilpotency index of any nilmanifold associated to $\Gg$). Dualizing, we obtain a filtration of $\Gg^*$ by the
annihilators $V_i$ of $\Gg^i$, which can also be expressed as
\begin{equation*}
V_i = \left\{v\in\Gg^*\ :\ dv\in \wedge^2 V_{i-1}\right\},
\end{equation*}
where $V_0=\{0\}$.  Choosing a basis for $V_1$ and extending
successively to a basis for each $V_k$, we obtain a {\it Malcev
basis} $\{e_1,\ldots,e_{n}\}$ for $\Gg^*$. This basis satisfies
the property
\begin{equation}\label{E:basis}
de_i \in {\wedge}^2\left<e_1, \ldots, e_{i-1}\right>\ \ \forall i.
\end{equation}
The filtration of $\Gg^*$ induces a filtration of its exterior
algebra, and leads to the following useful definition:
\begin{definition}
With $V_i$ as above, the {\it nilpotent degree} of a $p$-form \ga,
which we denote by  $\idx(\ga)$, is the smallest $i$ such that
$\ga \in {\wedge}^p V_i$.
\end{definition}

\begin{remark}
If \ga\ is a 1-form of nilpotent degree $i$ then $\idx(d\ga) =
i-1$.
\end{remark}

In this thesis, we specify the structure of a particular nilpotent
Lie algebra by listing the exterior derivatives of the elements of
a Malcev basis as an $n$-tuple of 2-forms $(de_k = \sum c^{ij}_k
e_i\wedge e_j)_{k=1}^m$. In low dimensions we use the shortened
notation $ij$ for the 2-form $e_i\wedge e_j$, as in the following
6-dimensional example: the 6-tuple $(0,0,0,12,13,14+35)$ describes
a nilpotent Lie algebra with dual $\Gg^*$ generated by 1-forms
$e_1,\ldots, e_6$ and such that $de_1=de_2=de_3=0$, while
$de_4=e_1\wedge e_2$, $de_5=e_1\wedge e_3$, and $de_6=e_1\wedge
e_4+e_3\wedge e_5$.  We see clearly that
$V_1=\left<e_1,e_2,e_3\right>$,
$V_2=\left<e_1,e_2,e_3,e_4,e_5\right>$, and $V_3=\Gg^*$, showing
that the nilpotency index of $\Gg$ is 3.

\section{Generalized Complex Structures on Nilmanifolds}\label{gen}

In this section, we present two results concerning generalized
complex structures on nilmanifolds of arbitrary dimension.  In
Theorem \ref{pure spinor}, we prove that any left-invariant
generalized complex structure on a nilmanifold must be generalized
Calabi-Yau, i.e. the canonical bundle $U_L$ has a closed
trivialization.  In Theorem \ref{T:maximal} we prove an upper bound
for the type of a left-invariant generalized complex structure,
depending only on crude information concerning the nilpotent
structure.

We begin by observing that a left-invariant generalized complex
structure must have constant type $k$ throughout the nilmanifold
$M^{2n}$, and its canonical bundle $U_L$ must be trivial.  Hence,
 we may choose a global trivialization of the form
\begin{equation*}
\rho = e^{B+i\om}\Omega,
\end{equation*}
where $B,\omega$ are real left-invariant 2-forms and $\Omega$ is a
globally decomposable complex $k$-form, i.e.
\begin{equation*}
\Omega=\theta_1\wedge\cdots\wedge\theta_k,
\end{equation*}
with each $\theta_i$ left-invariant. These data satisfy the
nondegeneracy condition
$\omega^{2n-2k}\wedge\Omega\wedge\overline{\Omega}\neq 0$ as well
as the integrability condition $d\rho = (X + \xi)\cdot\rho$ for
some section $X + \xi\in C^\infty(T\oplus T^*)$.  Since $\rho$ and
$d\rho$ are left-invariant, we can choose $X+\xi$ to be
left-invariant as well.

It is useful to order $\{\theta_1,\ldots,\theta_k\}$ according to
nilpotent degree, and also to choose them in such a way that
$\{\theta_j\ : \ \idx(\theta_j)>i\}$ is linearly independent
modulo $V_i$; this is possible according to the following lemma.
\begin{lem}\label{T:trivial}
It is possible to choose a left-invariant decomposition
$\Omega=\theta_1\wedge\cdots\wedge\theta_k$ such that
$\idx(\theta_i)\leq\idx(\theta_j)$ for $i<j$, and such that
$\{\theta_j\ :\ \idx(\theta_j) > i\}$ is linearly independent
modulo $V_i$ for each $i$.
\end{lem}
\begin{proof}
Choose any left-invariant decomposition
$\Omega=\theta_1\wedge\cdots\wedge\theta_k$ ordered according to
nilpotent degree, i.e. $\idx(\theta_i)\leq\idx(\theta_j)$ for
$i<j$.  Certainly $\{\theta_1,\ldots,\theta_k\}$ is linearly
independent modulo $V_0=\{0\}$.  Now let $\pi_i: \Gg^* \into
\Gg^*/V_i$ be the natural projection, and suppose
$\{\pi_{i}(\theta_j)\ :\ \idx(\theta_j)>i\}$ is linearly
independent for all $i<m$.  Consider $X=\{\pi_m(\theta_j)\ :\
\idx(\theta_j)>m\}$. If there is a linear relation
$\pi_m(\theta_{p}) = \sum_{l\neq p} \ga_l\pi_m(\theta_l)$ among
these elements, then we may replace $\theta_{p}$ with
$\tilde{\theta}_p=\theta_p - \sum_{l\neq p}\ga_l \theta_l$, which
does not change $\Omega$ or affect linear independence modulo
$V_i,\ i<m$.   However, note that $\idx(\tilde\theta_p)\leq m$,
i.e. we have removed an element from $X$.  Reordering by degree
and repeating the argument, we may remove any linear relation
modulo $V_m$ in this way, proving the induction step.
\end{proof}

We require a simple linear algebra fact before moving on
to the theorem.
\begin{lem}\label{T:injection}
Let $V$ be a subspace of a vector space $W$, let
$\alpha\in\wedge^p V$, and suppose
$\{\theta_1,\ldots,\theta_m\}\subset W$ is linearly independent
modulo $V$. Then
$\alpha\wedge\theta_{1}\wedge\cdots\wedge\theta_{m}=0$ if and only
if $\alpha=0$.
\end{lem}
\begin{proof}
Let $\pi:W\rightarrow W/V$ be the projection and choose a
splitting $W\cong V\oplus W/V$;
$\alpha\wedge\theta_{1}\wedge\cdots\wedge\theta_{m}$ has a
component in $\wedge^p V\otimes\wedge^m(W/V)$ equal to
$\alpha\otimes \pi(\theta_1)\wedge\cdots\wedge\pi(\theta_m)$,
which vanishes if and only if $\alpha=0$.
\end{proof}

\begin{theorem}\label{pure spinor}  Any left-invariant generalized
complex structure on a nilmanifold must be generalized Calabi-Yau.
That is, any left-invariant global trivialization $\rho$ of the
canonical bundle must be a closed differential form.  In
particular, any left-invariant complex structure has
holomorphically trivial canonical bundle.
\end{theorem}
\begin{proof}
Let $\rho=e^{B+i\omega}\Omega$ be a left-invariant trivialization
of the canonical bundle such that
$\Omega=\theta_1\wedge\cdots\wedge\theta_k$, with
$\{\theta_1,\ldots,\theta_k\}$ ordered according to
Lemma~\ref{T:trivial}. Then the integrability condition
$d\rho=(X+\xi)\cdot\rho$ is equivalent to the condition
\begin{equation}\label{E:equi}
d(B+i\omega)\wedge\Omega + d\Omega =
(i_X(B+i\omega))\wedge\Omega+i_X\Omega+\xi\wedge\Omega.
\end{equation}
The degree $k+1$ part of~(\ref{E:equi}) states that
\begin{equation}\label{E:eq1}
d\Omega = i_X(B+i\om)\wedge\Omega + \xi\wedge\Omega.
\end{equation}
Taking wedge of \eqref{E:eq1} with $\theta_i$ we get
\begin{equation*}
d\theta_i \wedge \theta_1\wedge\cdots\wedge\theta_k=0,\ \ \forall
i.
\end{equation*}
Now, let $\{\theta_1,\ldots,\theta_j\}$ be the subset with
nilpotent degree $\leq \idx(d\theta_i)$.  Note that $j<i$ since
$\idx(\theta_i)=\idx(d\theta_i)+1$.  Then since
\begin{equation*}
(d\theta_i\wedge\theta_1\wedge\cdots\wedge\theta_j)\wedge\theta_{j+1}\wedge\cdots\wedge\theta_k
= 0,
\end{equation*}
we conclude from Lemma~\ref{T:injection} that
\begin{equation}\label{E:important}
d\theta_i\wedge\theta_1\wedge\cdots\wedge\theta_j=0, \mbox{ with } j<i.
\end{equation}
Since this argument holds for all $i$, we conclude that
$d\Omega=0$. The degree $k+3$ part of~(\ref{E:equi}) states that
$d(B+i\omega)\wedge\Omega=0$, and so we obtain that $d\rho =
e^{B+i\omega}d\Omega=0$, as required.
\end{proof}

Equation \eqref{E:important} shows that the integrability
condition satisfied by $\rho$ leads to constraints on
$\{\theta_1,\ldots,\theta_k\}$. Since these will be used
frequently in the search for \gcss, we single them out as follows.

\begin{cor}\label{T:ideal}
Assume $\{\theta_1,\ldots,\theta_k\}$ are chosen according to
Lemma \ref{T:trivial}. Then
\begin{equation}\label{idealtheta}
d\theta_i \in \mc{I}(\{\theta_j : \idx(\theta_j) <
\idx(\theta_i)\}),
\end{equation}
where $\mc{I}(~)$ denotes the
ideal generated by its arguments. Since $\nil(\theta_i)$ is weakly
increasing, we have, in particular,
\begin{equation*}
d\theta_i\in\mathcal{I}(\theta_1,\ldots,\theta_{i-1}).
\end{equation*}
\end{cor}

\begin{example}\label{T:simple case1}
Since $d\theta_1 \in \mc{I}(0)$, we see that $\theta_1$ is always
closed, and therefore it lies on $V_1$ or, equivalently, $\idx(\theta_1)=1$.
\end{example}

So far, we have described constraints deriving from the
integrability condition on $\rho$.  However, nondegeneracy (in
particular, $\Omega\wedge\overline{\Omega}\neq 0$) also places
constraints on the $\theta_i$ appearing in the decomposition of
$\gO$.  The following example illustrates this.

\begin{example}\label{T:simple case2}
If $\theta_1,\ldots,\theta_j\in V_i$, then nondegeneracy implies
that $\dim V_i \geq 2j$.  For a fixed nilpotent algebra, this
places an upper bound on the number of $\theta_j$ which can be
chosen from each $V_i$.
\end{example}

In the next lemma we prove a similar, but more subtle constraint
on the 1-forms~$\theta_i$.

\begin{lem}\label{T:highest index}
Assume that $\{\theta_1,\ldots,\theta_k\}$ are chosen according to
Lemma \ref{T:trivial}. Suppose that no $\theta_i$ satisfies
$\idx(\theta_i)=j$, but that there exists $\theta_l$ with
$\idx(\theta_l) =j+1$. Then $\theta_l\wedge\overline\theta_l\neq
0$ modulo $V_j$ (i.e. in $\wedge^2(V_{j+1}/V_j)$), and in
particular $V_{j+1}/V_j$ must have dimension 2 or greater.
\end{lem}

\begin{proof}
Assume that the hypotheses hold but
$\theta_l\wedge\overline\theta_l=0$ modulo $V_j$.  Because of
this, it is possible to decompose $\theta_l=v+\alpha$, where
$\idx(\alpha)< j+1$, and such that $v\wedge\overline v = 0$.
Therefore, up to multiplication of $\theta_l$ (and therefore
$\Omega$) by a constant, $v$ is real.

By Corollary \ref{T:ideal}, $d\theta_l\in
\mathcal{I}(\{\theta_i : \nil(\theta_i)< j+1\})$. By hypothesis, there is no $\theta_i$ with nilpotent degree $j$, therefore we obtain
\begin{equation*}
dv +d\ga = \sum_{\idx(\theta_i) < j} \xi_i \wedge \theta_i,
\end{equation*}
where $\xi_i\in\Gg^*$.   Since $\idx(dv)=j$, there is an element
$w \in \Gg^{j-1}$ such that $i_w dv \neq 0$.  On the other hand,
the nilpotent degrees of $d\ga$ and the $\{\theta_i\}$ in the sum above are
less than $j$, hence interior product with $w$ annihilates all
these forms.  In particular,
\begin{equation*} 0 \neq i_w dv = \sum_{\idx(\theta_i) < j}(i_w\xi_i)\theta_i.
\end{equation*}
Therefore, $i_w\xi_i$ is nonzero for some $i$. But, the left hand side
is real, thus
\begin{equation*}0 = i_w dv \wedge \overline{i_w dv}
= \left(\sum_{\idx(\theta_i) < j}(i_w\xi_i)\theta_i\right) \wedge
\left(\sum_{\idx(\theta_i) <
j}(i_w\overline{\xi_i})\overline{\theta_i} \right).
\end{equation*}
By the nondegeneracy condition, the right hand side is nonzero,
which is a contradiction.  Hence
$\theta_l\wedge\overline\theta_l\neq 0$ modulo $V_j$.
\end{proof}

From this lemma, we see that if $\dim V_{j+1}/V_j = 1$ occurs in a
nilpotent Lie algebra, then it must be the case, either that some
$\theta_i$ has nilpotent degree $j$, or that no $\theta_i$ has
nilpotent degree $j+1$.  In this way, we see that the size of the
nilpotent steps $\dim V_{j+1}/V_j$ may constrain the possible
types of left-invariant generalized complex structures, as we now
make precise.

\begin{theorem}\label{T:maximal}
Let $M^{2n}$ be a nilmanifold with associated Lie algebra
$\frak{g}$. Suppose there exists a $j>0$ such that, for all $i\geq
j$,
\begin{equation}\label{jump}
\dim \left({V_{i+1}}/{V_{i}}\right) = 1.
\end{equation}
Then $M$ cannot admit left-invariant generalized complex
structures of type $k$ for $k \geq 2n - \nil(\Gg) +j $.

In particular, if $M$ has maximal nilpotency index (i.e. $\dim
V_1=2,\ \dim V_i/V_{i-1}=1\ \forall i>1$), then it does not admit
\gcss\ of type $k$ for $k\geq 2$.
\end{theorem}
\begin{proof}

First observe that for any nilpotent Lie algebra $\nil(\Gg)\leq
2n-1$, so the theorem only restricts the existence of structures
of type $k>1$.

According to Lemma~\ref{T:highest index}, if none of the
$\theta_i$ have nilpotent degree $j$, then there can be none with
nilpotent degree $j+1, j+2, \ldots$ by the condition~(\ref{jump}).
Hence we obtain upper bounds for the nilpotent degrees of
$\{\theta_1,\ldots,\theta_k\}$, as follows.  First, $\theta_1$ has
nilpotent degree 1 (by Example~\ref{T:simple case1}).  Then, if
$\idx(\theta_2)\geq j+2$, this would imply that no $\theta_i$ had
nilpotent degree $j+1$, which is a contradiction by the previous
paragraph.  Hence $\idx(\theta_2) < j+2$.  In general,
$\idx(\theta_i) < j+i$.  Suppose that $M$ admits a generalized
complex structure of type $k>1$.  Then we see that
$\idx(\theta_k)< j+k$. By Example~\ref{T:simple case2}, we see
this means that $\dim V_{j+k-1}\geq 2k$.

On the other hand, $\dim V_{j+k-1} = 2n - \dim \Gg^*/V_{j+k-1}$,
and since $\Gg^*=V_{\nil(\Gg)}$, we have
\begin{equation*}
\Gg^*/V_{j+k-1} = V_{\nil(\Gg)}/V_{j+k-1} \cong
V_{\nil(\Gg)}/V_{\nil(\Gg)-1}\oplus\cdots\oplus V_{j+k}/V_{j+k-1},
\end{equation*}
and the $\nil(\Gg)-j-k+1$ summands on the right have dimension 1,
by hypothesis.  Hence $\dim V_{j+k-1} = 2n-\nil(\Gg)+j+k-1$, and
combining with the previous inequality, we obtain
\begin{equation*}
k < 2n-\nil(\Gg)+j,
\end{equation*}
as required. For the last claim, observe that $M^{2n}$ has maximal
nilpotency index when $\nil(M) = 2n-1$, in which case
$\eqref{jump}$ holds for $j=1$.
\end{proof}

Constraints beyond those mentioned here may be obtained if one
considers the fact that $\gO\wedge \overline{\gO}$ defines a
foliation and that $\go$ restricts to a symplectic form on each
leaf.  Both the leafwise nondegeneracy of $\omega$ and the
requirement of being closed on the leaves lead to useful
constraints on what types of generalized complex structure may
exist, as we shall see in the following sections.

\section{Generalized Complex Structures on
6-nilmanifolds}\label{six}

In this section, we turn to the particular case of 6-dimensional
nilmanifolds. The problem of classifying those which admit
left-invariant complex (type 3) and symplectic (type 0) structures
has already been solved \cite{Sa01,GoKh96}, so we are left with the task
of determining which 6-nilmanifolds admit left-invariant
generalized complex structures of types $1$ and $2$. The result of
our classification is presented in Table~1, where explicit
examples of all types of left-invariant generalized complex
structures are given, whenever they exist.  The main results
establishing this classification are Theorems \ref{T:2,1} and
\ref{T:1,2}.  Throughout this section we often require the use of
a Malcev basis $\{e_1,\ldots,e_6\}$ as well as its dual basis
$\{\del_1,\ldots,\del_6\}$.  We use $e_{i_1\cdots i_p}$ as an
abbreviation for $e_{i_1}\wedge\cdots\wedge e_{i_p}$.

\subsection{Generalized Complex Structures of Type 2}

By the results of the last section, a left-invariant structure of
type $2$ is given by a closed form $\rho = \exp(B+i\om) \theta_1
\theta_2$, where
$\om\wedge\theta_1\theta_2\overline{\theta}_1\overline{\theta}_2\neq
0$. As a consequence of Theorem \ref{T:maximal}, any 6-nilmanifold
with maximal nilpotence step cannot admit this kind of structure.
We now rule out some additional nilmanifolds, using Lemma
\ref{T:highest index}.

\begin{lem}\label{lemmafirst}
If a 6-nilmanifold $M$ has nilpotent Lie algebra given by
$(0,0,0,12,14,-)$, and has nilpotency index $4$, then $M$ does not
admit left-invariant generalized complex structures of type $2$.
\end{lem}
\begin{proof}
Suppose that $M$ admits a structure $\rho$ of type $2$.  Since $M$
has nilpotency index $4$, $V_{i+1}/V_i$ is 1-dimensional for $i
\geq 1$. From $d\theta_1=0$ and Lemma \ref{T:highest index}, we
see that $\theta_1= z_1e_1+z_2e_2+z_3e_3$ and $\nil(\theta_2) \leq
2$, thus $\theta_2= w_1e_1+w_2e_2+w_3e_3+w_4e_4$. The conditions
$d(\theta_1 \theta_2)=0$ and
$\theta_1\theta_2\overline{\theta_1}\overline{\theta_2} \neq 0$
together imply $z_3=0$.  Further, the annihilator of
$\theta_1\theta_2\overline{\theta_1}\overline{\theta_2}$ is
generated by $\{\del_5,\del_6\}$.  Hence, the nondegeneracy
condition $\omega^2\wedge\Omega\wedge\overline{\Omega}\neq 0$
implies that
$$B+i\om = (k_1e_1 +
\dots k_4 e_4 + k_5e_5)e_6 + \alpha,$$ where $k_5 \neq 0$ and $\ga
\in \wedge^2\left<e_1,\cdots,e_5\right>$. But, using the structure
constants, we see that $d\rho$ must contain a nonzero multiple of
$e_6$, and so cannot be closed.
\end{proof}

\begin{lem}\label{lemmasecond}
Nilmanifolds associated to the algebras defined by
\begin{align*}
&(0,0,0,12,14,13-24),\\
&(0,0,0,12,14,23+24)
\end{align*}
do not admit left-invariant generalized complex structures of type
$2$.
\end{lem}
\begin{proof}
Each of these nilmanifolds has $\nil(\Gg)=3$, with $\dim V_1=3$
and $\dim V_2/V_1=1$.  Suppose either nilmanifold admitted a
structure $\rho$ of type $2$.  If $\nil(\theta_2) = 2$, we could
use the argument of the previous lemma to obtain a contradiction.
Hence, suppose $\idx(\theta_2)=3$.  Lemma \ref{T:highest index}
then implies that $\theta_2\wedge\overline{\theta_2}\neq 0\pmod
{V_2}$, which means it must have nonzero $e_5$ and $e_6$
components.

Now, if $\theta_1$ had a nonzero $e_3$ component, then
$d\theta_2\wedge\theta_1$ would have nonzero $e_{234}$ component,
contradicting~(\ref{idealtheta}).  Hence
\begin{equation}
\theta_1= z_1 e_1 + z_2 e_2 \qquad \theta_2= \sum_{i=1}^6 w_i e_i.
\end{equation}
But for these, the coefficient of $e_{123}$ in
$d\theta_2\wedge\theta_1$ would be $- z_2 w_6$ (for the first
nilmanifold) or $z_1 w_6$ (for the second), in each case implying
$\theta_1\wedge\overline{\theta_1}=0$, a contradiction.
\end{proof}

\begin{lem}\label{T:lemmathird}
Nilmanifolds associated to the algebras defined by
\begin{align*}
&(0,0,12,13,23,14),\\
&(0,0,12,13,23,14-25)
\end{align*}
do not admit left-invariant generalized complex structures of type
$2$.
\end{lem}
\begin{proof}
Each of these nilmanifolds have $\nil(\Gg)=4$, with $\dim V_1=2$,
$\dim V_2=3$, and $\dim V_3=5$.  Suppose either nilmanifold
admitted a structure $\rho$ of type $2$. $V_4/V_3$ is
1-dimensional and so Lemma \ref{T:highest index} implies that
$\idx(\theta_2) \neq 4$. Since $\theta_1$ is closed and satisfies
$\theta_1\overline{\theta}_1\neq 0$, we may rescale it to obtain
$\theta_1=e_1+z_2e_2$. The condition $d\theta_2 \in
\mathcal{I}(\theta_1)$ implies we can write $\theta_2 = w_2e_2 +
w_3e_3 + w_4(e_4 + z_2e_5)$. Now let
\begin{equation*}
B+i\om=\sum_{i<j}k_{ij}e_{ij}
\end{equation*}
and differentiate $\rho$ using the structure constants.  In both
cases, $d\rho=0$ implies $B+i\om = \xi_1\theta_1 + \xi_2\theta_2$
for 1-forms $\xi_i$. Therefore \go\ is degenerate on the leaves
defined by
$\Ann(\theta_1\theta_2\overline{\theta_1}\overline{\theta_2})$,
contradicting the requirement
$\omega^2\wedge\Omega\wedge\overline\Omega\neq 0$.
\end{proof}

\begin{theo}{2,1}
The only 6-dimensional nilmanifolds not admitting left-invariant
generalized complex structures of type 2 are those with maximal nilpotency index and those excluded by
Lemmas~\ref{lemmafirst}, \ref{lemmasecond}, and \ref{T:lemmathird}.
\end{theo}
\begin{proof}
In Table 1, in the the end of this chapter, we provide explicit left-invariant generalized complex
structures of type 2 for all those not excluded by the preceding
lemmas and Theorem \ref{T:maximal}.
\end{proof}

\subsection{Generalized Complex Structures of Type 1}

A left-invariant generalized complex structure of type $1$ is
given by $\rho =\exp(B+i\om)\theta_1$, where $\om^2\wedge
\theta_1\overline{\theta_1} \neq 0$.  Note that this implies that
$\omega$ is a symplectic form on the $4$-dimensional leaves of the
foliation determined by $\theta_1\wedge\overline{\theta_1}$.

\begin{theo}{1,2}
The only 6-nilmanifolds which do not admit left-invariant
generalized complex structures of type $1$ are those associated to
the algebras defined by $(0,0,12,13,23,14)$ and
$(0,0,12,13,23,14-25)$.
\end{theo}
\begin{proof}
In Table 1, at the end of this chapter, we provide explicit forms defining type 1 structures
for all 6-nilmanifolds except the two mentioned above.

Suppose the nilmanifold is associated to the Lie algebra
$(0,0,12,13,23,14)$. Then up to a choice of Malcev basis, a
generalized complex structure of type 1 can be written
\begin{equation*}
\rho = \exp(B+i\om)(e_1+z_2e_2), \qquad B+i\om =
\sum_{i<j}k_{ij}e_{ij}.
\end{equation*}
The condition $d\rho=0$ implies that
\begin{equation*}
(-k_{36}e_{314} +k_{45}e_{135} - k_{45}e_{423} + k_{46}e_{136} +
k_{56}e_{236}-k_{56}e_{514})(e_1+z_2e_2)=0.
\end{equation*}
The vanishing of the $e_{1245}$, $e_{1236}$, $e_{1235}$, and
$e_{1234}$ components imply successively that $k_{56}$, $k_{46}$,
$k_{45}$, and $k_{36}$ all vanish.

The leaves of the distribution $\Ann(\theta_1\overline{\theta_1})$
are the tori generated by $\del_3,\del_4,\del_5,\del_6$, and the
previous conditions on $B+i\omega$ imply that on these leaves,
\go\ restricts to $\im(k_{34})e_{34} +
\mbox{Im}(k_{35})e_{35}$ which is degenerate, contradicting
$\om^2\wedge\theta_1\overline{\theta}_1\neq 0$. An identical
argument holds for the nilpotent algebra $(0,0,12,13,23,14-25)$.
\end{proof}

\section{$\beta$-transforms of Generalized Complex
Structures}\label{iwa}

In this section, we will use Theorem~\ref{T:deformation} to
show that any left-invariant complex structure on a nilmanifold
may be deformed into a left-invariant generalized complex
structure of type 1.  By connecting the type 3 and type 1
structures in this way, we go on to show that the two disconnected
components of the left-invariant complex moduli space on the
Iwasawa manifold may be joined by paths of generalized complex
structures.

\begin{theo}{complex}
Every left-invariant complex structure on a $2n$-dimensional
nilmanifold can be deformed, via a \gb-field, into a
left-invariant generalized complex structure of type $n-2$.
\end{theo}
\begin{proof}
According to Theorem \ref{T:deformation}, such a deformation
can be obtained if we find a holomorphic Poisson structure.  Let
us construct such a bivector $\beta$.

By Theorem~\ref{pure spinor}, a left-invariant complex structure
on a nilmanifold has a holomorphically trivial canonical bundle.
Let $\Omega=\theta_1\wedge\cdots\wedge\theta_n$ be a holomorphic
volume form, and choose the $\theta_i$ according to
Lemma~\ref{T:trivial}, so that, by Corollary~\ref{T:ideal}, the
differential forms
$\theta_1,\theta_1\wedge\theta_2,\ldots,\theta_1\wedge\cdots\wedge\theta_n$
are all holomorphic.  Now let $\{x_1,\ldots,x_n\}$ be a dual basis
for the holomorphic tangent bundle.  By interior product with
$\Omega$, we see that the multivectors $x_n,x_{n-1}\wedge
x_n,\ldots,x_1\wedge\cdots\wedge x_n$ are all holomorphic as well.
In particular, we have a holomorphic bivector $\beta=x_{n-1}\wedge
x_n$.  Calculating the Schouten bracket of this bivector with
itself, we obtain
\begin{equation*}
[\beta,\beta]=[x_{n-1}\wedge x_n,x_{n-1}\wedge
x_n]=2[x_{n-1},x_n]\wedge x_{n-1}\wedge x_n=0,
\end{equation*}
where the last equality follows from the fact that
$[x_{n-1},x_n]\in\left<x_{n-1},x_n\right>$, since
$\theta_i([x_{n-1},x_n])=-d\theta_i(x_{n-1},x_n)=0$ for $i<n-1$,
by Corollary~\ref{T:ideal}.

Hence $\beta$ gives rise to a deformation of the generalized
complex structure. The resulting structure $\tilde\rho$ is given
by the following differential form:
\begin{equation*}
\tilde\rho=e^{i_\beta}\rho = \rho + i_\beta\rho =
e^{\theta_{n-1}\wedge\theta_n}\theta_1 \wedge \dots \wedge
\theta_{n-2},
\end{equation*}
and we see immediately that it is a left-invariant generalized
complex structure of type $n-2$.
\end{proof}

\begin{ex}{iwasawa}
In \cite{KS03}, Ketsetzis and Salamon study the space of
left-invariant complex structures on the Iwasawa nilmanifold. This
manifold is the quotient of the complex 3-dimensional Heisenberg
group of unipotent matrices by the lattice of unipotent matrices
with Gaussian integer entries. As a nilmanifold, it has structure
$(0,0,0,0,13-24,14+23)$. Ketsetzis and Salamon observe that the
space of left-invariant complex structures with fixed orientation
has two connected components which are distinguished by the
orientation they induce on the complex subspace $\langle \del_1,
\del_2,\del_3, \del_4 \rangle$.

\vskip6pt \noindent {\bf Connecting the two components.} Consider
the left-invariant complex structures defined by the closed
differential forms $\rho_1 = (e_1+ie_2)(e_3+ie_4)(e_5+ie_6)$ and
$\rho_2 = (e_1+ie_2)(e_3-ie_4)(e_5-ie_6)$.  These complex
structures clearly induce opposite orientations on the space
$\langle \del_1, \dots, \del_4 \rangle$, so lie in different
connected components of the moduli space of left-invariant complex
structures.

By Theorem \ref{T:complex}, the first complex structure can be
deformed, by the $\beta$-field $\beta_1 = \frac{-1}{4}(x_3 -
ix_4)(x_5-ix_6)$ into the generalized complex structure
$$e^\beta \rho_1 = e^{-(e_{35} - e_{46}) - i(e_{45} + e_{36})}(e_1+ie_2),$$
and then, by the action of the closed $B$-field $B_1= e_{35} - e_{46}$, into
$$\rho = e^{ - i(e_{45} + e_{36})}(e_1+ie_2).$$

On the other hand, the second complex structure can be deformed,
via the $\beta$-field $\beta_2 = \frac{1}{4}(x_3 +
ix_4)(x_5+ix_6)$, into the type 1 generalized complex structure
$$e^\beta \rho_2 = e^{(e_{35} - e_{46}) - i(e_{45} + e_{36})}(e_1+ie_2),$$
and then by the $B$-field $B_2 = -(e_{35} - e_{46})$ into
$$ \rho =e^{ - i(e_{45} + e_{36})}(e_1+ie_2),$$
which is the same generalized complex structure obtained from
$\rho_1$.

Therefore, by deforming by $\beta$- and $B$-fields, the two
disconnected components of the moduli space of left-invariant
complex structures can be connected, through generalized complex
structures.
\end{ex}

\begin{ex}{type 2 to symplectic}
The standard complex structure on the 4-torus or in the nilmanifold $(0,0,0,12)$ can be deformed into the $B$-field transform of a symplectic structure, by Theorem \ref{T:complex}. If a nilmanifold $N$ is a 2-torus bundle over either of these and there is a cohomology class $a \in H^2(N)$ which restricts to the symplectic form on the 2-torus, then, by Theorem \ref{T:Thurston}, $N$ will have a type 2 structure which can deformed into a symplectic structure, by Example \ref{T:fibrations}. Moreover, in the case of invariant structures, it is easy to see that if the symplectic form of a type 2 structure is closed, then the structure can be deformed into a symplectic one. Using the forms from Table 1, this explains why the following nilmanifolds have both type 2 and symplectic structures:
\begin{align*}
&(0,0,0,12,13,14+23) &\quad& (0,0,0,12,13,24)\\
&(0,0,0,12,13,14) && (0,0,0,12,13,23)\\
&(0,0,0,0,12,13) && (0,0,0,0,0,12).
\end{align*}
On the other hand, the nilmanifold $(0,0,0,0,0,12+34)$ admits a type 2 structure determined by the form $(1+i2)(3+i4)\exp i(56)$ and admits no symplectic structure. This shows that in Example \ref{T:fibrations} we actually need the 2-form $\tau$ to be closed for deformations on the base to give rise to deformations of the total space.
\end{ex}

\section{An 8-dimensional Example}\label{8d}

We have established that all 6-dimensional nilmanifolds admit
generalized complex structures.  One might ask whether every
even-dimensional nilmanifold admits a left-invariant generalized
complex structure.  In this section, we answer this question in the
negative by presenting an 8-dimensional nilmanifold which does
not admit any type of left-invariant generalized complex
structure.

\begin{ex}{8d}
Consider a nilmanifold $M$ associated to the 8-dimensional
nilpotent Lie algebra defined by
\begin{equation*}
(0,0,12,13,14,15,16,36-45-27).
\end{equation*}

Since it has maximal nilpotency index, Theorem \ref{T:maximal}
implies that it may only admit left-invariant generalized complex
structures of types 1 and 0.  We exclude each case in turn:

\begin{itemize}
\item \emph{Type 1}: Suppose there is a type 1 structure, defined
by the left-invariant form $\rho = e^{B+i\om}\theta_1$.  Then
$d\theta_1$ =0 and $\theta_1\overline{\theta}_1\neq 0$ imply that
$\theta_1\overline{\theta}_1$ is a multiple of $e_{12}$ and
therefore $\om$ must be symplectic along the leaves defined by
$\langle\del_3,\ldots,\del_8\rangle$.  These leaves are actually
nilmanifolds associated to the nilpotent algebra defined by
$(0,0,0,0,0,12+34)$, which admits no symplectic structure, and so
we obtain a contradiction.
\item \emph{Type 0}:
The real second cohomology of $M$ is given by
\begin{equation*}
H^2(M,\RR) = \langle e_{23},e_{34}-e_{25},e_{17}\rangle,
\end{equation*}
and since $e_8$ does not appear in any of its generators, it is
clear that any element in $H^2(M,\RR)$ has vanishing fourth power,
hence excluding the existence of a symplectic structure.
\end{itemize}
In this way, we see that the 8-dimensional nilmanifold $M$ given
above admits no left-invariant generalized complex structures at
all.
\end{ex}
\setlength{\oddsidemargin}{0 cm}
\setlength{\evensidemargin}{0 cm}
\setlength{\topmargin}{-1.0 cm}
\setlength{\footskip}{0.0cm}
\setlength{\parindent}{2em}
\setlength{\textwidth}{19 cm}
\setlength{\textheight}{23.5 cm}
\renewcommand{\baselinestretch}{1}
\pagestyle{empty}
\begin{landscape}

\section{Table of Generalized Complex Structures on 6-nilmanifolds}
In the table below we provide explicit examples of pure forms defining generalized complex structures for all nilmanifolds of dimension 6 admiting such structures. In the first column are the structure constants of the Lie algebra associated to the nilmanifold, $b_1$ and $b_2$ are the first and second Betti numbers and in the last four columns we provide differential forms defining generalized complex structures of type 3,2,1 and 0. Observe that for the symplectic case we just write the symplectic form $\go$ instead of $e^{i \go}$.
\vskip18pt
\begin{center}
{\tiny
\begin{tabular}{|l c c l l l l|}
\hline
Nilmanifold class & $b_1$ & $b_2$ & Complex (type 3) &Type 2 & Type 1 & Symplectic (type 0)\\
\hline $(0,0,12,13,14,15)$ & 2 & 3 &  -- & -- &
$(1+i2)\exp{i(36-45)}$ & $16 + 34
- 25$ \\
$(0,0,12,13,14,34+52)$ & 2 & 2 & -- & -- & $(1+i2)\exp(- 45 +  36 + i(36+ 45))$ & --\\
$(0,0,12,13,14,23+15)$ & 2 & 3 & -- & -- & $(1+i2)\exp i(36-45)$ &
$16 +
24 + 34 -25$\\
$(0,0,12,13,23,14)$ & 2 & 4 & -- &  -- & -- & $15 + 24 + 34 - 26$\\
$(0,0,12,13,23,14-25)$ & 2 & 4 & -- & -- & -- & $15 + 24 - 35 + 16$\\
$(0,0,12,13,23,14+25)$ & 2 & 4 & $(1+i2)(4+i5)(3+i6)$ & $(1 +
i2)(4+i5)\exp i(36)$ &
$(1 + i2)\exp(43-56+i(46 - 35))$ & $15 + 24 + 35 + 16$\\
$(0,0,12,13,14+23,34+52)$ & 2 & 2 & -- & -- & $(1+i2)\exp(45-35+36+i(-36+45-16))$  & --\\
$(0,0,12,13,14+23,24+15)$ & 2 & 3 & -- & -- & $(1+i2)\exp(2 \times 35 +
i(36- 45))$ & $16 + 2 \times 34 - 25$\\
\hline $(0,0,0,12,13,14+35)$ & 3 & 5 & -- & $(1+2+i3)(4+i5)\exp
i(26) $ &$(1+i2)\exp i(
 36 + 45) $ & --\\

$(0,0,0,12,13,14+23)$ & 3 & 6 & $(1+i2)(3 - 2 \times i4)(5+
2\times i6)$ & $(2+i3)(4+i5)\exp i(16-34) $ & $(1+i2)\exp
i(36+45)$ & $16 -2\times
34 - 25$ \\

$(0,0,0,12,13,24)$ & 3 & 6 & $(1+i2)(3+4+i4)(5+6-i6)$ &
$(1+2+i3)(4+i5)\exp i(26) $ &
 $(1+i2)\exp i(35+46)$ & $26+14+35 $ \\

$(0,0,0,12,13,14)$ & 3 & 6 & $(1+i2)(3+i4)(5+i6)$ & $(2+i3) (4+i5)
 \exp i(16)$ & $(1+i2)\exp(35-46+i(36+45))$ & $16+24+35$  \\

$(0,0,0,12,13,23)$ & 3 & 8 & $(1+i2)(3+i4)(5+i6)$ &
$(1+i2)(5+i6)\exp i(16-34)$ & $(1+i2)\exp(35-46 + i(36 + 45))$ &
$15 + 24
+36$ \\

$(0,0,0,12,14,15+23)$ & 3 & 5 & -- & -- & $(1+i3)\exp i(26-45)$ & $13 +26
- 45$ \\

$(0,0,0,12,14,15 + 23 + 24)$ & 3 & 5 & -- & -- & $(1+i3) \exp i(26-45)$ &
$13+26-45$\\

$(0,0,0,12,14,15+24)$ & 3 & 5 & -- & -- & $(1+i3) \exp i(26-45)$ &
$13+26-45$\\

$(0,0,0,12,14,15)$ & 3 & 5 & -- & -- & $(1+i3) \exp i(26-45)$ &
$13+26-45$\\

$(0,0,0,12,14,24)$ & 3 & 5 & $(1+i2)(3+i4)(5+i6)$ &
$(1+i2)(5+i6)\exp i(34)$ & $(1+i2)\exp(35-46 + i(45+36))$ & --\\

$(0,0,0,12,14,13+42)$ & 3 & 5 & $(1+i2)(3+i4)(2 \times 5 - i6)$ &
-- &  $(1+i2)\exp i(35 + 46)$ & $15 +
26 + 34$\\

$(0,0,0,12,14,23+24)$ & 3 & 5 & $(1+i2)(3+4+i3)(5+6+i6)$ & -- &
$(1+i2)\exp(35+46+i(35-46))$ & $16-34+25$\\

$(0,0,0,12,23,14+35)$ & 3 & 5 & -- & $(1+2+i3)(5+i4)\exp(3+i1)6$ & $(1+i2)\exp(36 + 45 +i(36 - 45))$ & --\\

$(0,0,0,12,23,14-35)$ & 3 & 5 & $(1+i3)(4-i5)(2+i6)$ &
$(1+i3)(4-i5)\exp i(26)$ & $(1+i3)\exp(24 + 56 + i(25+46))$ & --\\

$(0,0,0,12,14-23,15+34)$ & 3 & 4 & -- &  $(1+i2)(3+i4)\exp i(56)$
&
 $(2+i3) \exp i(16 + 35 + 45-26)$  & $16+35+24$ \\

$(0,0,0,12,14+23,13+42)$ & 3 & 5 & $(1+i2)(3-i4)(5+i6)$ &
$(1+i2)(3-i4)\exp i(56)$ & $(1+i2)\exp(35+46 + i(36-45))$ &
 $15+2\times 26+34$\\

\hline

$(0,0,0,0,12,15+34)$ & 4 & 6 & -- & $(1+i2)(3+i4)\exp i(56)$ & $(3+i4)\exp i(25+16)$ & --\\
$(0,0,0,0,12,15)$ & 4 & 7 & -- & $(1+i2)(3+i4)\exp i(56)$ &
$(1+i2)\exp i(34 + 56)$ & $16+25+34$\\
$(0,0,0,0,12,14+25)$ & 4 & 7 & $(1+i2)(4+i5)(3+i6)$ &
$(1+i2)(4+i5)\exp i(36)$ & $ (1+i2) \exp(34 + 56 + i(35 - 46))$& $13+26+45$\\
$(0,0,0,0,12,14+23)$ & 4 & 8 & $ (1+i2)(3-i4)(5+i6)$ &
$(1+i2)(3-i4)\exp i(56)$ & $(1+i2)\exp(35 + 46 +i(36 - 45))$ & $13+26+45$\\
$(0,0,0,0,12,34)$ & 4 & 8 & $(1+i2)(3+i4)(5+i6)$&
$(1+i2)(3+i4)\exp i(56)$ & $(1+i2)\exp(35-46 + i(45 + 36))$ & $15+36+24$\\
$(0,0,0,0,12,13)$ & 4 & 9 & $(2+i3)(1+i4)(5+i6)$ &
$(2+i3)(5+i6)\exp i(14)$ & $(2+i3)\exp(15-46 +i(16+45))$ & $16+25+34$ \\
($0,0,0,0,13+42,14+23)$ & 4 & 8 & $(1+i2)(3-i4)(5+i6)$ &
$(1+i2)(3+i4)\exp i(56)$ & $(1+i2)\exp(35+46+i(36-45))$ & $16+25+34$ \\
\hline $(0,0,0,0,0,12+34)$ & 5 & 9 & $(1+i2)(3+i4)(5+i6)$ &
$(1+i2)(3+i4)\exp i(56)$ &
$(1+i2)\exp i(36+45)$ & --\\
$(0,0,0,0,0,12)$ & 5 & 11 & $(1+i2)(3+i4)(5+i6)$ &
$(1+i2)(5+i6)\exp i(34)$ & $(1+i2)\exp i(36+45)$& $16+23+45$\\
 \hline
$(0,0,0,0,0,0)$ & 6 & 15 & $(1+i2)(3+i4)(5+i6)$ &
$(1+i2)(3+i4)\exp i(56)$ &
$(1+i2)\exp i(36+45)$ & $12+34+56$\\
\hline
\end{tabular}\label{table}}
\vskip18pt {\small {\it Table 1: Differential forms defining
left-invariant Generalized Calabi-Yau structures. The symbol `-'
denotes nonexistence.}}
\end{center}
\end{landscape}

\pagestyle{headings}

\setlength{\oddsidemargin}{0.5 cm}
\setlength{\evensidemargin}{0 cm}
\setlength{\topmargin}{0.0 cm}
\setlength{\footskip}{1.0cm}
\setlength{\parindent}{2em}
\setlength{\textwidth}{16 cm}
\setlength{\textheight}{21.5 cm}
\renewcommand{\baselinestretch}{1.5}

\chapter{The $dd^{\mathcal{J}}$-lemma}\label{ddj-lemma}

From now we move away from finding examples of \gcss\ and focus on a generalization of the $dd^c$-lemma. The key to generalize this property is Theorem \ref{T:U_k to U_k+1}, which allowed us to define generalized versions of \del, \delbar\ and $d^c$ --- denoted by \del, \delbar\ and $d^{\J}$ --- which agree with those operators for a \gcs\ induced by a complex structure on a manifold. Recall that if the complex $dd^c$-lemma holds on a complex manifold, then  not only will its cohomology split into $H^{p,q}$ according to the splitting of the exterior algebra $\wedge^{\bullet}TM\tensor \C$ into $\wedge^{p,q}TM$, but also the manifold is {\it formal}. Ultimately, we are interested in determining which of these properties will also be a consequence of the generalized $dd^c$-lemma.

In this chapter we establish the first (and easier one): if the generalized $dd^c$-lemma holds,  the cohomology of the manifold splits according to the splitting of $\wedge^{\bullet}TM\tensor \C$ into $U^k$ (as defined on page \pageref{decomposition of forms}). We also introduce an analogue of the Fr\"olicher spectral sequence, which we call the {\it canonical spectral sequence}, to give a converse to the above result: if the canonical spectral sequence degenerates at $E_1$ and there is a decomposition of cohomology, the generalized $dd^c$-lemma holds.

With this decomposition at hand, we also prove that the Poincar\'e dual of a generalized complex submanifold, suitably changed by its associated 2-form, is always represented by an element in a specific component of the splitting of the cohomology. This result is analogous to the fact that a complex submanifold of a complex manifold represents a $(p,p)$-class and a Lagrangian submanifold of a symplectic manifold represents a primitive class in middle dimension.

The relation between the generalized $dd^c$-lemma and formality will have to wait for the further two chapters, and will come from the study of this property in the symplectic setting.

This chapter is organized as follows. In the first section we introduce the $dd^{\J}$-lemma and prove that it induces a splitting in cohomology. We also give one condition that implies that this lemma holds which will be useful in the symplectic case. In  Section \ref{canonical spectral sequence}  we introduce the canonical spectral sequence and use a result by Deligne {\it et al} \cite{DGMS75} on double complexes to prove that the splitting of cohomology together with the degeneracy of the canonical spectral sequence at $E_1$ implies the $dd^{\J}$-lemma. Finally, in the last section, we prove that generalized complex submanifolds are represented by classes in $\mc{U}^0$.

\section{The $dd^{\J}$-lemma}\label{ddj lemma}

In Theorem \ref{T:U_k to U_k+1}, we saw that, in the presence of a \gcs\ \J, the exterior derivative decomposes as $d = \del + \delbar$ and there we defined an operator $d^{\J} = - i(\del-\delbar)$. In the complex case, these operators are just the standard \del, \delbar\ and $d^c$. Based on this, we find the following definition relevant:

\begin{definition}
A \gc\ manifold {\it satisfies the $dd^{\J}$-lemma} if
$$ \im d \cap \ker d^{\J} =\im d^{\J} \cap \ker d = \im dd^{\J}.$$
\end{definition}

\begin{remarks}
\item We could equivalently have said that $M$ satisfies the $\del\delbar$-lemma if
$$ \im \del \cap \ker \delbar =\im \delbar \cap \ker \del = \im \del\delbar.$$
It is easy to see that these properties are equivalent.
\item As $d^{\J} = \J^{-1}d \J$, one can easily check that if $\im d \cap \ker d^{\J} = \im dd^{\J}$, then the second equality follows. Indeed, if $\ga \in \im d^{\J} \cap \ker d$, then $\J \ga \in \im d \cap \ker d^{\J}$ and hence $\J \ga = dd^{\J} \beta$, for some $\beta$. Therefore
$$\ga = -\J dd^{\J}\beta = \J d \J^{-1} d \J \beta = d^{\J} d (\J \beta),$$
and hence $\ga \in \im dd^{\J}$.
\end{remarks}

As we will see in Chapter \ref{symplectic_ddbar}, the $dd^{\J}$-lemma in the symplectic case is equivalent to the strong Lefschetz property, or just Lefschetz property, for short. Also, both the complex $dd^c$-lemma and the Lefschetz property hold for compact K\"ahler manifolds, therefore it is reasonable to conjecture that this property holds for any twisted \gk\ manifold. This is actually true and has been proved recently by Gualtieri.

\begin{theo}{gk&ddbar}{\em (Gualtieri \cite{Gu04})}
Any compact twisted generalized K\"ahler manifold satisfies the $dd^{\J}$-lemma with respect to both complex structures involved.
\end{theo}

On the other hand, when it comes to examples of nonk\"ahlerian symplectic manifolds, the most standard ones are nontoroidal nilmanifolds. There are essentially two arguments showing that these are not K\"ahler:
\begin{renumerate}
\item Any nontoroidal nilmanifold is not formal, hence can not have a complex structure satisfying the $dd^c$-lemma \cite{CFC87,Has89};
\item Any nontoroidal nilmanifold admiting a symplectic structure does not satisfy the Lefschetz property \cite{BG88}.
\end{renumerate}

Therefore it is reasonable to conjecture:
\begin{conj}{nilmanifolds&ddbar}
If a nilmanifold admits an (invariant) \gcs\ for which the $dd^{\J}$-lemma holds, this nilmanifold is a torus.
\end{conj}

One further advantage of having the expression $d^{\J} =  \J^{-1} d \J$ is that this allows us to define the operator $d^{\J}$ without having to determine the decomposition of $\wedge^{\bullet}T^*M$ into $U_k$ beforehand, i.e, we can define $d^{\J}$ without prior knowledge about $\del$ and $\delbar$. But still the $dd^{\J}$-lemma has implications about a decomposition of the cohomology of the manifold.

\begin{theo}{decomposition2}
The following properties are equivalent for a \gc\ manifold $(M,\J)$:
\begin{enumerate}
\item $M$ satisfies the $dd^{\J}$-lemma;
\item The inclusion of the complex of $d^{\J}$-closed forms $\gO_{d^{\J}}$ into the complex of differential forms $\gO$ induces an isomorphism in cohomology:
$$(\gO_{d^{\J}},d) \stackrel{i}{\hookrightarrow} (\gO,d),\qquad H^{\bullet}(\gO_{d^\J}) \stackrel{i^*}{\cong}H^{\bullet}(\gO).$$
\end{enumerate}
And both imply that every cohomology class $a \in H^{\bullet}(M)$ can be represented by a form $\ga = \sum \ga_k$, with $\ga_k \in \mc{U}^k$ such that $d\ga_k=0$. If $a=0$ is the trivial cohomology class, then for any such \ga\ each $\ga_k$ is exact. In particular, each cohomology class can be represented by a $d$ and $d^{\J}$-closed form.
\end{theo}
\begin{proof}
We start with the implication (1) $\Rightarrow$ (2), which uses part of the argument from Theorem \ref{T:formality1}.
\begin{renumerate}
\item $i^*:H^{\bullet}(\gO_{d^\J}) \into H^{\bullet}(\gO)$ is injective:

If $i^* \ga$ is exact, then \ga\ is $d^{\J}$-closed and exact, hence by
the $dd^{\J}$-lemma $\ga = dd^{\J}\beta$, so \ga\ is the derivative of a
$d^{\J}$-closed form and hence its cohomology class in $\gO_{d^{\J}}$ is
also zero.

\item $i^*:H^{\bullet}(\gO_{d^\J}) \into H^{\bullet}(\gO)$ is surjective:

Let \ga\ be a closed form and set $\beta = d^{\J}\ga$. Then $d\beta = d
d^{\J} \ga = - d^{\J} d \ga = 0$, so $\beta$ satisfies the conditions of the
$dd^{\J}$-lemma, hence $\beta = d^{\J} d \gamma$. Let $\tilde{\ga} = \ga -
d\gamma$, then $d^{\J}\tilde{\ga} = d^{\J}\ga - d^{\J} d \gamma = \beta - \beta
=0$, so $[\ga] \in \im(i^*)$.
\end{renumerate}

For the converse, assume that $dd^{\J} \ga =0$, i.e., $d^{\J}\ga \in \im d^{\J} \cap \ker d$. We want to prove that $d^{\J}\ga = dd^{\J} \beta$, for some $\beta$.
Since $d^{\J} d\ga= 0$, $d\ga \in \mc{E}_{d^{\J}}$ is a closed form in $\gO_{d^{\J}}$ and represents the trivial cohomology class in $\gO$, hence it also represents the trivial cohomology class in $\gO_{d^{\J}}$, i.e., there is $\gamma_1 \in \gO_{d^{\J}}$ such that $d\ga = d \gamma_1$.  Therefore $d(\ga - \gamma_1) =0$ 
and the cohomology class $[\ga - \gamma_1]$ has a $d^{\J}$-closed representative $\gamma_2$:
$$ \ga - \gamma_1 = \gamma_2 + d\beta.$$
Applying $d^{\J}$, we get, $d^{\J} \ga = d^{\J} d \beta,$
as we wanted.

To prove the last claim, let $a$ be a cohomology class. Then there is a representative $\ga$ for it which is $d$- and $d^c$-closed. Since $d = \del + \delbar$ and $d^{\J}= -i(\del - \delbar)$, we conclude that \ga\ is both \del\ and \delbar\ closed, and so must be each of its components $\ga_k$ relative to the splitting. Hence we obtain a splitting for the cohomology class $a = \sum [\ga_k]$. If $a$ was the trivial cohomology class, any such $\ga$ would be $dd^{\J}\beta$, for some $\beta$, from ({\it ii\,}) above. So the decomposition of \ga\ would be $\ga_k = dd^{\J}\beta_k$, showing that each of the summands is exact.

\end{proof}

\begin{definition}
If $(M,\mc{J})$ is a \gc\ manifold satisfying the $dd^{\J}$-lemma, we define the {\it generalized cohomology of $M$}, $HH^{k}(M)$, as the cohomology classes in $H^{\bullet}(M)$ that can be represented by sections of $\mc{U}^k$.
\end{definition}

\begin{remark}
In the case of a Calabi-Yau manifold, the generalized cohomology agrees with the Hochschild homology \cite{Cal03}. This result motivates our notation.
\end{remark}

Now, the Mukai pairing on forms is essentially the wedge product and therefore induces a pairing in cohomology, i.e., if one represents two cohomology classes $a$ and $b$ by forms $\ga$ and $\gb$, then the product
$$(a,b) : = [(\ga,\gb)],$$
is independent of the representatives chosen, where $[\cdot]$ indicates taking the cohomology class of the argument. With that said, Proposition \ref{T:mukai pairing on U_k} gives us the following:

\begin{prop}{<HHk,HH-k>}
On a compact \gc\ manifold $M$ satisfying the $dd^{\J}$-lemma, the Mukai pairing vanishes in $HH^k(M)\times HH^l(M)$ unless $k = -l$, in which case it is nondegenerate.
\end{prop}
\begin{proof}
Just represent classes $a \in HH^k(M)$ and $b \in HH^{l}(M)$ by sections \ga\ and \gb\ of $U^k$ and $U^l$ respectively. Then, from Proposition \ref{T:mukai pairing on U_k},
$$(a,b) = [(\ga,\gb)] = [0],$$
if $k \neq -l$.
Now, still with $a$ and \ga\ as above, by Poincar\'e duality, there is a closed form \gb\ such that $[(\ga,\gb)]\neq 0$. Decomposing $[\beta]$ into its generalized cohomology components, $b_l \in HH^l(M)$ we get that
$$(a,b_{-k}) = (a, \sum b_l) = (a,[\gb])= [(\ga,\gb)] \neq 0,$$
proving nondegeneracy.
\end{proof}

Next we give another condition equivalent to the $dd^{\J}$-lemma, this time involving the operators \del\ and \delbar. This is essentially the argument used by Merkulov in the symplectic case to prove the equivalence of the Lefschetz property and the $dd^{\J}$-lemma (cf. Theorem \ref{T:ddelta lemma}):

\begin{theo}{merkulov}
If in a \gc\ manifold $M^{2n}$ every $\del$-cohomology class has a $\delbar$-closed representative and
$$ \im \del \cap \ker \delbar = \im \delbar \cap \ker \del,$$
then $M$ satisfies the $\del\delbar$-lemma.
\end{theo}
\begin{proof}
This is proved by induction. I'll use superscripts to keep track of the spaces that the forms belong to, e.g. $\ga^k \in \mc{U}^k$. Assume that $\ga^{n} \in \im \delbar \cap \ker \del  \cap \mc{U}^{n}$. Then $\ga^{n} = \delbar\beta^{n+1}=0$, as $\mc{U}^{n+1} = \{0\}$. Therefore $\ga^{n} =\del\delbar 0$.

Going one step down, let  $\ga^{n-1} \in \im \delbar \cap \ker \del \cap \mc{U}^{n-1}$, say
$$\ga^{n-1} = \delbar\beta^{n}.$$
Since $\del\beta^{n} \in \mc{U}^{n+1} =\{0\}$, by hypothesis, there is a \delbar-closed form $\tilde{\beta}^n$ in the same $\del$-cohomology class as $\beta^{n}$, i.e. $\beta^{n} = \tilde{\beta}^{n} + \del\xi^{n-1}$ and hence
$$\ga^{n-1} = \delbar \beta^{n} = \delbar(\tilde{\beta}^{n} + \del\xi^{n-1}) = \delbar\del\xi^{n-1}.$$

For the general case, let $\ga^{k} \in \im \delbar \cap \ker \del \cap \mc{U}^{k}$, say
$$\ga^{k} = \delbar \beta^{k+1}.$$
As $\del\delbar \beta^{k+1}= \del \ga =0$, we get that
$$\del\beta \in \im \del \cap \ker \delbar \cap \mc{U}^{k+2} = \im \delbar \cap \ker \del \cap \mc{U}^{k+2},$$
by the hypothesis of the theorem. And then, by the induction hypothesis, $\del \beta^{k+1} = \del\delbar \xi^{k+2}$, hence $\del(\beta^{k+1} - \delbar \xi_1^{k+2})=0$. As any \del-cohomology class has a \delbar-closed representative, we can find $\eta^{k}$ and \delbar-closed $\xi_2^{k+1}$ such that
$$\beta^{k+1} - \delbar \xi_1^{k+2} = \xi_2^{k+1} + \del \eta^{k},$$
and, by applying $\delbar$,
$$ \ga^k = \delbar \beta^{k+1} = \delbar\del \eta^{k},$$
which proves the induction step.
\end{proof}

\section{The Canonical Spectral Sequence}\label{canonical spectral sequence}

In \cite{DGMS75}, Deligne {\it et al} give a series of properties equivalent to the $dd^c$-lemma for a complex manifold. These properties are stated in the context of double complexes and the result we are concerned about here is a converse to Theorem \ref{T:decomposition2} stating when a decomposition of cohomology into $HH^k$ implies the $dd^{\J}$-lemma.

The only obstacle is that while in the complex case there is a natural bigrading for the complex of differential forms, we have found only one grading for forms on a \gc\ manifold, namely, the one given by the $U^k$. A way to remedy this is to mimic Brylinski \cite{Br88} and Goodwillie \cite{Goo85}: introduce a formal element $\beta$ of degree 2 and consider the complex:
$$\mc{A} = \gO^{\bullet}(M) \tensor \odot \mbox{span}\{\beta,\beta^{-1}\},$$
and to change the differential to:
$$d^{\beta}(a \beta^k) = \del a \beta^k + \delbar a \beta^{k+1}.$$
The complex $\mc{A}$, which we call the {\it canonical complex}, has a bigrading given by $\mc{A}^{p,q} = \mc{U}^{p-q}\beta^q$, and the differential $d^{\beta}$ decomposes as $\del^{\beta}$ and $\delbar^{\beta}$, where $\del^{\beta}:\mc{A}^{p,q} \into \mc{A}^{p+1,q}$ and $\delbar^{\beta}:\mc{A}^{p,q}\into \mc{A}^{p,q+1}$. The complex of differential forms sits inside $\mc{A}$ as the $\beta$-periodic elements:
$$\tau:\gO \into \mc{A};\qquad \tau(\ga) = \sum_{k \in\Z} \ga \beta^{k}.$$
And this is a map of differential algebras which preserves the decompositions of $d$ and $d^{\beta}$:
$$ \tau(\del \ga) = \del^{\beta} \tau(\ga) \qquad \mbox{ and } \qquad \tau(\delbar \ga) = \delbar^{\beta} \tau(\ga).$$
One can easily check that the $\del\delbar$-lemma holds for $\gO^{\bullet}$ if and only if the corresponding lemma holds for $\mc{A}$.

Also, the bigrading gives a filtration of $\mc{A}$:
$$F^p\mc{A} = \sum_{p'\geq p}\mc{A}^{p',q};$$
$$F^{p}\mc{A}^m = \sum_{p' \geq p} \mc{A}^{p',m-p'}.$$
Which is preserved by $d^{\beta}$, i.e., $d^{\beta}:F^p\mc{A} \into F^p\mc{A}$. For each $m$, $F^{p}\mc{A}^{m} = \{0\}$ for $2p \geq n+m$ and $F^{p}\mc{A}^{m} = \mc{A}^{m}$ for $2p \leq m-n$, $2n$ is the dimension of the manifold.
This means that the filtration is {\it bounded} and hence the induced spectral sequence, which we call the {\it canonical spectral sequence}, converges to the cohomology of the operator $d^{\beta}$. This spectral sequence is periodic in the sense that $E^{p,q}_r \cong E^{p-q,0}_r$. If we define:

\begin{definition}
The {\it $\del$-cohomology} of a \gcm\ is given by:
$$H_{\del}^k = \frac{\ker \del:U^k \into U^{k+1}}{\im \del:U^k \into U^{k+1}}.$$
\end{definition}

Then the first term $E^{p,q}_1 \cong H_{\del}^{p-q}$ is just the $\del$-cohomology of the manifold, which is finite dimensional, since, by a result of Gualtieri \cite{Gu03}, \del\ is an elliptic operator. The second term is the cohomology induced by $\delbar$ in $H_\del$, and the sequence goes on. However, if the $\del\delbar$-lemma holds this sequence degenerates at $E_1$. Conversely, Deligne's `theorem' (\cite{DGMS75}, Proposition 5.17 and remark 5.21) tells us that the degeneracy at $E_1$ together with the decomposition of cohomology imply the $dd^{\J}$-lemma:

\begin{theo}{canonical spectral sequence} If the canonical spectral sequence degenerates at $E_1$ and the decomposition of forms into subbundles $U^k$ induces a decomposition in cohomology, then the $dd^{\J}$-lemma holds.
\end{theo}

\begin{remark}
For a \gcs\ induced by a complex structure the canonical spectral sequence is just the Fr\"olicher spectral sequence repeated over and over. In particular, the degeneracy of the canonical spectral sequence at $E_r$ is equivalent to the degeneracy of the Fr\"olicher spectral sequence at the same stage. In Chapter \ref{symplectic_ddbar} we will see that, in the symplectic case, this spectral sequence is isomorphic to the spectral sequence coming from what Brylinski calls the `canonical complex' \cite{Br88}. This is what motivates our terminology.
\end{remark}

As an application of this theorem, Parshin \cite{Pa66} shows that if two compact complex manifolds are birationally equivalent via a holomorphic map, then one of them satisfies the $dd^{\J}$-lemma if and only if the other one does. Recalling that in the symplectic setting there is a blow-up operation, which is topologically similar to the complex blow-up, one can ask if the same is true in the symplectic setting. As we will see in Chapter \ref{blowup} this is not the case.

It is possible that the canonical spectral sequence degenerates at $E_1$ even though the $dd^{\J}$-lemma does not hold. One example of such phenomenon is given by a result by Kodaira \cite{Kod64} stating that the Fr\"olicher --- and hence the canonical --- spectral sequence always degenerates at $E_1$ for complex surfaces, although not all of those satisfy the $dd^c$-lemma.  We will encounter another example where this happens in Chapter \ref{symplectic_ddbar} where we will prove that it also degenerates at $E_1$ for  any symplectic structure --- the same result was established before by Brylinski \cite{Br88}.

Finally, following Fr\"olicher \cite{Fro55}, we have the following result.
\begin{prop}{Frolicher}
If $M^{2n}$ admits a \gcs, then the Euler characteristic of $M$ is given by
$$\chi(M) = \pm \sum (-1)^{k} \dim H_{\del}^{k}(M).$$
Where the sign is $+$ if the elements in $U^0$ are even forms and $-$ otherwise. 
\end{prop}
\begin{proof}
Given the periodic condition, $E^{p,q}_r \cong E^{p-q,0}_r$, this spectral sequence is equivalent to long exact sequences:
$$\cdots \into \mc{U}^{k-1} \stackrel{d_1= \del}{\into} \mc{U}^{k} \stackrel{d_1 = \del}{\into} \mc{U}^{k+1} \into \cdots;$$
$$\cdots \into H_{\del}^{k-3} \stackrel{d_2}{\into} H_{\del}^{k} \stackrel{d_2}{\into} H_{\del}^{k+3} \into \cdots;$$
$$ \cdots \into H_{d_{r-1}}^{k-2r+1} \stackrel{d_r}{\into} H_{d_{r-1}}^{k} \stackrel{d_r}{\into} H_{d_{r-1}}^{k+2r-1} \into \cdots.$$
As $d_r$ maps $ev/od$ to $od/ev$, the Euler characteristic is preserved and hence can be computed from the first sequence where the spaces involved are finite dimensional.
\end{proof}

\subsection{Twisted Cohomology and Lie Groups}\label{twisted cohomology}

As Theorem \ref{T:U_k to U_k+1} also holds in the twisted case, we still have a decomposition, $d_H= \del_H + \delbar_H$ and hence we have operators $\del_H$, $\delbar_H$ and $d_H^{\J} = -i(\del_H - \delbar_H) = \J^{-1}d_H\J$ defined for \tgcm. By an abuse of language, we will still denote $\del_H$, $\delbar_H$ and $d^{\J}_H$ by \del, \delbar\ and $d^{\J}$, with the understanding that they also depend on the twisting 3-form, if the \gcs\ is twisted. Given that, it also makes sense to state the `$dd^{\J}$-lemma' in this setting.

\begin{definition}
A \tgcm\ $(M,\J,H)$ satisfies the $d_Hd^{\J}$-lemma if
$$ \im d_H \cap \ker d^{\J} = \im d^{\J} \cap \ker d_H = \im d_Hd^{\J}.$$
\end{definition}

Theorems \ref{T:decomposition2} and \ref{T:merkulov} also hold in this setting, only with the obvious change of replacing $d$ by $d_H$ and the ordinary cohomology by the $d_H$-cohomology in Theorem \ref{T:decomposition2}. It is interesting to notice that the twisted cohomology normally has only a $\Z_{2}$ grading, but the presence of the $d_{H}d^{\J}$-lemma  gives a whole $\Z$ grading.

\begin{ex}{liegroups2} {\sc --- Lie groups and the $d_Hd^{\J}$-lemma.} Given Gualtieri's theorem  (Theorem \ref{T:gk&ddbar}), twisted generalized K\"ahler manifolds will have their $d_H$-cohomology splitting according to Theorem \ref{T:decomposition2}. In particular, this holds for semi-simple Lie groups (cf. Example \ref{T:liegroups}), but, from Example \ref{T:liegroups dh-cohomology}, the $d_H$-cohomology of those groups is trivial. Note that this provides no contradiction to the existence of such a splitting, on the contrary, as the twisted cohomology is trivial, every cohomology class has a $d^{\J}$-closed representative. 
\end{ex}

\section{Submanifolds}\label{submanifolds2}

Recall from Section \ref{submanifolds} that a \gc\ submanifold $(M,0)$ of a complex manifold $N$ is just a standard complex submanifold. Assuming that the ambient manifold satisfies the $dd^{\J}$-lemma, the Poincar\'e dual to any such submanifold, $PD(M)$ will be an integral $(p,p)$-cohomology class, hence any submanifold (of arbitrary dimension) is represented by an integral class in $\oplus H^{p,p}(N) = HH^0(N)$.

In this section we show that this is not mere coincidence and indeed, whenever the cohomology of a \gc\ manifold splits into $HH^k$, \gc\ submanifolds will be represented by elements in $HH^0$. The need to consider the $2$-form is also evident.

\begin{theo}{submanifolds in HH0}
Let $(N^{2n},H,\J)$ be a twisted generalized complex manifold, let $(M^m,F)$ be a generalized complex submanifold and $\ga \in \mc{U}^k(N)$. Then, for each $p \in M$, the degree $m$ part of $e^F \wedge \ga|_{T_pM}$ vanishes if $k \neq 0$. If the twisted cohomology of $N$ splits into generalized cohomology, then $[e^{-F} PD(M)]$ is a class in $HH^0(N)$.
\end{theo}

\begin{proof}
For $p \in M$, let $\mu \in \wedge^{2n-m}T_p^*N$ be a volume form of a vector space complementary to $T_p M$ so that, for $X \in T_pM$, $X\cdot \mu =0$. The restriction $\xi$ of the degree $m$ part of $e^F \wedge \ga$ to $T_p M$ can be computed using the Mukai pairing: $\xi = (1,e^F \wedge \ga|_{T_pM})_{M}$.

This form is nonzero if and only if it wedges nonzero with $\mu$, i.e., if and only if
 $$ 0 \neq \mu \wedge (1,e^F \wedge \ga|_{T_pM})_{M}=  \pm (\mu,e^F \wedge \ga)_{N} = \pm (e^{-F}\wedge \mu,\ga)_{N}.$$

Since $e^{-F}\mu$ annihilates the generalized tangent space $\tau_F$ of $M$ at $p$, Lemma  \ref{T:maximal real in U_0} implies that $e^{-F}\mu \in U^0_p$ and hence its Mukai pairing with any form $\ga \in U^k_p$ vanishes, if $k \neq 0$, according to Proposition \ref{T:mukai pairing on U_k}.

If the twisted cohomology of $N$ splits into generalized cohomology, then the first part of the theorem shows that a \gc\ submanifold induces a functional on $HH^{\bullet}(N)$ which vanishes in $HH^{k}(N)$, for $k \neq 0$. From the definition of Poincar\'e dual, it is clear that this functional coincides with taking Mukai pairing with the $d_H$-cohomology class given by $[e^{-F}PD(M)]$. Therefore $[e^{-F}PD(M)]$ is an element of $HH^0(N)$.
\end{proof}

One question can be raised from this theorem:
\vskip6pt\noindent
{\sc Open question}: {\it Which elements in $HH^0(N)$ can be represented by submanifolds?}
\vskip6pt
Observe that the problem of finding which cohomology classes can be represented by submanifolds comes with the requirement that the cohomology class has to be integral. This is {\it not} the case when one considers a submanifold with $2$-form $(M,F)$ as there is no integrality restriction on $F$ and therefore neither on $[e^{-F}PD(M)]$.

This problem might be compared to the Hodge conjecture --- the search for algebraic cycles representing integral cohomology classes in algebraic varieties --- or, in symplectic geometry, to the problem of characterizing which primitive integral cohomology classes can be represented by Lagrangian submanifolds.
 

\chapter{Symplectic $d\delta$-lemma and the Lefschetz Property}\label{symplectic_ddbar}

As we have seen in Chapter \ref{gcss}, Section \ref{linear algebra of gcs}, a generalized almost complex structure induces a splitting of the differential forms into subbundles $U^k$ and, in the presence of a $dd^{\J}$-lemma, this splitting induces a splitting in cohomology (Theorem \ref{T:decomposition2}). In the complex case, the $dd^{\J}$-lemma is just the standard $dd^c$-lemma and the decomposition of forms is just given by $U^k = \sum_{p-q=k} \Omega^{p,q}$.

In the symplectic case, the corresponding decomposition has not yet appeared in the literature, although the corresponding $dd^{\J}$-lemma has. In \cite{Kos85} Koszul introduced an operator for Poisson manifolds which agrees with $d^{\J}$ in the symplectic case (in which case we denote it by $\delta$). In \cite{Br88}, Brylinski studied further this operator in conjunction with Libermann's $*$ operator \cite{Lib54,LiMa87}, which resembles the Hodge $*$. Brylinski called forms that are both $d$ and $\delta$ closed {\it symplectic harmonic}. Later, Yan \cite{Ya96} and Mathieu \cite{Ma90} proved independently that the existence of a harmonic representative in each cohomology class is equivalent to the strong Lefschetz property (or just Lefschetz property, for short).

\vskip6pt
\noindent
{\bf Lefschetz property.} A symplectic manifold $(M^{2n},\go)$ satisfies the {\it Lefschetz property at level $k$} if the map
$$[\go^{n-k}] : H^{k} \into H^{2n-k}$$
is surjective. It satisfies the {\it Lefschetz property} if these maps are surjective for $0 \leq k \leq n$.
\vskip6pt

 Finally, Merkulov \cite{Me98} proved that for a compact manifold the existence of symplectic harmonic forms in each cohomology class (and therefore the Lefschetz property) is indeed equivalent to the $d\delta$-lemma:

\setlength{\unitlength}{1mm}
\begin{picture}(60,10)(-65,-5)
\linethickness{1pt}
\thinlines
\put(-50,0){\makebox(0,0){Lefschetz property}}
\put(5,0){\makebox(0,0){Harmonic representatives}}
\put(55,0){\makebox(0,0){$d\delta$-lemma}}
\put(35,0){\makebox(0,0){$\Leftrightarrow$}}
\put(-25,0){\makebox(0,0){$\Leftrightarrow$}}
\end{picture}

One of the points of this chapter is to introduce the theory of symplectic harmonic forms, as studied by Brylinski, Yan and Merkulov and give complete proofs to theorems quoted above. Although most of this is already standard and used in many other places, Merkulov's theorem is still contested with remarkable frequency. This seems to happen mostly for two reasons:  his proof seems to be missing the proof for a claim (here Lemmata \ref{T:lemma2c3} and \ref{T:lemma3} and Proposition \ref{T:weak ddelta lemma}) and the final conclusion about formality does not follow, at least if one is considering formality in the sense of Sullivan  as defined on page \pageref{formalitypage}. So I hope these issues will be settled by the presentation here.

The second purpose of this chapter is to use this theory to describe symplectic geometry from a generalized complex point of view. We determine how the exterior algebra of the cotangent bundle decomposes into the subbundles $U^k$ and show that $U^k$ is canonically isomorphic to $\wedge^{n-k}T^*M$ via a map $\gf$. Further we show that, besides the usual relations between $d$, $\delta$, \del\ and \delbar, there are more striking ones:
$$ \gf(d \ga) = \del \gf(\ga) \qquad \mbox{ and } \qquad \gf(\delta \ga) = -2i \delbar \gf(\ga).$$
These relations show that in some sense it is equivalent to work with $d$, $\delta$ and the degree decomposition of forms and \del, \delbar\ and the $U^k$ decomposition. Also, these relations can be used to prove that the canonical spectral sequence always degenerates at $E_1$ in the symplectic case, as $\gf$ provides an isomorphism of $\del$ and $d$ cohomologies --- this result had been established before by Brylinski by different means.

This chapter is organized as follows: In Section \ref{delta*} we review the theory of symplectic harmonic forms, as introduced by Brylinski and in Section \ref{sl2} we present Yan's theory (the $\frak{sl}(2,\C)$ representation on $\wedge^{\bullet}T^*M$). Then we are able to characterize the bundles $U^k$ and establish the forementioned  relations between $d$, $\delta$, \del\ and \delbar\ in Section \ref{decomposition}. After that we finally prove the main theorems of this chapter stating that the Lefschetz property, existence of harmonic representatives on every cohomology class and the $d\delta$-lemma are equivalent.

\section{The $\delta$ and $*$ Operators}\label{delta*}

In this section we follow Brylinski's approach \cite{Br88} and revisit the symplectic $d^{\J}$, introduce a new operator, the symplectic $*$, and show some relations between them.

On a symplectic manifold, one can consider two different actions of $\omega^{-1}$ on forms. One of them is the one that appears in the definition of the generalized complex structure $\omega^{-1}:T^*M \into TM$ (which can be extended to the whole exterior algebra of $T^*M$), the other is the corresponding Lie algebra action in the spin representation $\omega^{-1}:\wedge^k T^*M \into \wedge^{k-2}T^*M$. To avoid confusion, we will denote the second action by $\Lambda$. Moreover, for a $k$-form \ga, $\go^{-1}\ga$ is a $k$-vector, which naturally evaluates on $k$-forms, so we denote $\go^{-1}(\ga,\beta) =(\go^{-1}\ga)\lfloor\beta$. According to the definition of $d^{\J}$, which we denote by $\delta$ in the symplectic case, we have (see Example \ref{T:symplectic lie algebra action})
$$\delta = [d,\J] = \Lambda d - d \Lambda.$$

Another important operator can be obtained by considering the action of $\omega^{-1}$. We define the symplectic $*$ on acting on $k$-forms by the following identity
$$ \alpha \wedge *\beta = \omega^{-1}(\alpha, \beta) v_M,$$
for any $k$-form $\beta$, where $v_M$ is the volume form induced by the symplectic form, $v_M = \omega^n/n!$. Clearly $*$ is linear and $*\ga(p)$ depends only on the value of \ga\
at $p$. It is also worth mentioning that for a function $f$,
$*f=f v_M$. 

Assume that $(M_i^{2n_i},\go_i)$, $i=1,2$, are symplectic manifolds of dimension
$2n_i$ and let $(M,\go)$ be the product with the symplectic form $\go =
\go_1+\go_2$. Further, denote by $*_i$ the operator $*$ in
each $M_i$ and by $*$ the same operator in $M$. Choosing Darboux coodinates
$(x^i,y^i)$ for $M_i$, we also have Darboux coordinates for $M$
given by $(x^1,x^2,y^1,y^2)$ and these coordinates determine a
decomposition of any $k$-form \ga\ in $M$ as $\ga = \sum \ga^1_j \wedge
\ga^2_{k-j}$ with $\deg(\ga^i_j)= j$, $\ga^i_{*} \in
\mbox{span}\{dx^i_1, \cdots dy^i_{n_i}\}$. So, in order to know the
effect of $*$ on a $k$-form it is enough to study its effect on forms
of the type $\ga^1 \wedge \ga^2$, with $\ga^i \in
\mbox{span}\{dx^i_1, \cdots dy^i_{n_i}\}$. Let $k_i = \deg(\ga^i)$ and
$k= k_1+k_2$. Then:
\begin{align*}
(\beta^1 \wedge \beta^2) \wedge * (\ga^1 \wedge \ga^2) &= \omega^{-1}(\beta^1
\wedge \beta^2, \ga^1 \wedge \ga^2) v_M\\
&= \omega_1^{-1}(\beta^1, \ga^1) \omega_2^{-1}(\beta^2, \ga^2) v_{M_1} \wedge
v_{M_2}\\
&= \omega_1^{-1}(\beta^1, \ga^1)v_{M_1}\wedge (\omega_2^{-1}(\beta^2, \ga^2)) v_{M_2}\\
&= \beta^1 \wedge *_1 \ga^1 \wedge \beta^2 \wedge *_2 \ga^2 =
(-1)^{k_1 k_2} \beta^1 \wedge \beta^2 \wedge *_1 \ga^1 \wedge *_2 \ga^2.
\end{align*}
Hence, 
\begin{equation}\label{E:product}
* (\ga^1 \wedge \ga^2)  =  *_2 \ga^2  \wedge*_1 \ga^1.
\end{equation}

\begin{lem}\label{T:symplectic*}(Brylinski \cite{Br88}):
The symplectic $*$ satisfies $** = id$.
\end{lem}

\begin{proof}
If M is two dimensional, then choosing symplectic
coordinates, we have:

\begin{renumerate}
\item For a function, $*f=f dx \wedge dy$; 
\item For the 1-form $dx$,
\begin{align*}
 dx \wedge *dx &= \omega^{-1}(dx, dx) dx \wedge dy  =
0 \mbox{ and }\\
dy \wedge *dx &= \go^{-1}(dy, dx) dx \wedge dy =  dx\wedge dy
\end{align*}  
so $*dx = -dx$ and, analogously, $*dy = - dy$. Hence, for an arbitrary
1-form, $*(fdx+gdy)= f*dx+ g*dy = -(fdx+gdy)$;
\item For the volume form $v_M= dx \wedge dy$,
$$ v_M \wedge *v_M = \omega^{-1}(v_M,v_M)v_M = \det
\begin{pmatrix}
\overset{0}{\overbrace{\go^{-1}(dx, dx)}} &
\overset{-1}{\overbrace{\go^{-1}(dx, dy)}}\\
\overset{1}{\overbrace{\go^{-1}(dy, dx)}} &
\overset{0}{\overbrace{\go^{-1}(dy, dy)}}\\
\end{pmatrix} v_M = v_M. $$ 
\end{renumerate}

Thus, $*f v_m = f$, and for the case of a two dimensional manifold we
are done.

For the general case, by induction over the dimension of $M$, we use
Darboux coordinates to decompose a neighbourhood of a point as a
product of lower dimensional symplectic open sets and use \eqref{E:product}
$$
 ** (\ga^1 \wedge \ga^2) = *(*_2 \ga^2 \wedge *_1 \ga^1) = *_1*_1 \ga^1
\wedge *_2*_2 \ga^2 = \ga^1 \wedge \ga^2,
$$
as required.
\end{proof}

\begin{prop}{adjoint}
{\em (Brylinski \cite{Br88}):}
When acting on $k$ forms, $\delta = (-1)^{k+1} * d *$.
\end{prop}
\begin{proof}
We use induction again and start with a 2-dimensional manifold in which
case, after choosing Darboux coordinates, we have

\begin{renumerate}
\item For functions, $\delta f= 0$ and $* d * f = * d (f dx
\wedge dy) = 0$;
\item For 1-forms of the type $fdx$,
\begin{align*}
\delta(f dx) &= (\Lambda d - d \Lambda)(f dx) = \Lambda(f_y dy \wedge dx) = f_y\\
&= * f_y dx \wedge dy = * d(- f dx) = * d * f dx
\end{align*}
The case of forms of type $f dy$ is analogous; 
\item For 2 forms, letting again $v_M = dx \wedge dy$,
$$ \delta(f v_M) = (\Lambda d - d\Lambda)(f v_M) = - d \Lambda f dx \wedge dy  =
df = d*(f dx \wedge dy) = - * d * (f v_M).$$ 
\end{renumerate}

In higher dimensions, we use more or less the same argument used
before to split a neighbourhood of a point as a product of lower
dimensional symplectic manifolds $V = V_1 \times V_2$ and then proceed
by induction, but now we have to be aware that $\delta$ and $d$ are
not tensors. So instead of considering arbitrary $\ga^1$ and $\ga^2$
we only look at the ones that are pull-backs of forms in $V_1$ and
$V_2$. For such forms,
\begin{align*}
*d*(\ga^1 \wedge \ga^2) &= *d (*_2 \ga^2 *_1 \ga^1) =
 *(d(*_2\ga^2)\wedge *_1\ga^1 + (-1)^{k_2} *_2 \ga^2 d(*_1 \ga^1))\\
&= \ga^1 \wedge *_2 d *_2 \ga^2 + (-1)^{k_2} *_1 d *_1 \ga^1 \wedge
\ga^2\\
&= (-1)^{k_2+1} \ga^1 \wedge \delta \ga^2 + (-1)^{k_2+k_1+1}
\delta \ga^1 \wedge \ga^2 \\
& = (-1)^{k_1+k_2+1}\delta(\ga^1 \wedge \ga^2) = (-1)^{k+1}
\delta(\ga^1 \wedge \ga^2),
\end{align*}
where we have used the fact that in these circumstances (i.e., $\ga^i
= \pi_i^*\beta^i$ for some $\beta^i \in \Omega^*V_i$) one can prove
$$ \delta(\pi_1^*\beta^1\wedge\pi_2^*\beta^2) = \delta
(\pi_1^*\beta^1)\wedge\pi_2^*\beta^2 +(-1)^{k_1}
\pi_1^*\beta^1\wedge\delta(\pi_2^*\beta^2)$$
from the definition of $\delta$.

So, in the space generated by the $k$-forms of the form
$\pi_1^*\beta^1\wedge\pi_2^*\beta^2$ we have the relation $\delta =
(-1)^k * d *$. Therefore the same relation
is valid in the $C^{\infty}$ closure of this space. Since $V$ can be taken
to be small,
this closure is the whole space of differential forms in $V$.
\end{proof}

\begin{remark}
It is important to notice that $\delta$ does not satisfy the Leibniz rule. Indeed, in \cite{LinRe03}, Lin and Sjamaar prove that if $f$ is a function and $\ga$ a form then
$$ \delta(f\ga) = f\delta \ga - X_f \lfloor \ga,$$
where $X_f$ is the Hamiltonian vector field of $f$.
\end{remark}

 Based on Proposition \ref{T:adjoint}, Brylinski drew an analogy with Riemannian geometry and called $\delta$ the {\it symplectic adjoint} of $d$. We will not use this name in this thesis, but this explains the next definition.

\begin{definition}
A form is called {\em symplectic harmonic} if it is $d$ and $\delta$-closed.
\end{definition}

Observe that any closed $1$-form is also $\delta$ closed and therefore symplectic harmonic, hence there may be more than one symplectic harmonic form in a given cohomology class. Also the `laplacian', $d\delta+ \delta d$, vanishes. This shows that the parallel with harmonic forms is not a good one and we should think of $\delta$ more appropriately as an analogue of $d^c$ than one of $d^*$. Still, historically, it was this analogy with Riemannian geometry that justified the question: (Brylinski \cite{Br88}) {\it Is there a symplectic harmonic form representing each cohomology class?}

\section{The $\frak{sl}(2,\C)$ Representation}\label{sl2}

In this section we present Yan's theory developed to answer the question above. We shall postpone the proof of Yan's theorem (Theorem \ref{T:yan2}, stating that there is a symplectic harmonic form in each cohomology class if and only if the manifold satisfies the Lefschetz property) to Section \ref{ddelta}, as the theory developed turns out to be very interesting and useful also in the generalized complex point of view and for clarity we shall go first through the generalized complex framework.

We start by letting $L$ be the wedge product with $\go$ and $H=[L,\Lambda]$. Then these operators satisfy the following commutation relations:
\begin{renumerate}
\item $[d,L] = 0$, since \go\ is closed;

\item If $\Pi_k:  \gO^{\bullet}(M) \into  \gO^{\bullet}(M)$ is the
projection on the summand  $\gO^k(M)$, then
$$H = \sum (n-k) \Pi_k.$$

We will show this by induction on the dimension of $M$. If $M$ is two dimensional, then, for
functions, we have
$$[L, \Lambda]f = L\Lambda f - \Lambda L f = -f \Lambda dx\wedge dy =
f.$$
For a 1-form \ga\, $[L, \Lambda] \ga = 0$ since $L\ga$ is a
3-form in a 2-dimensional manifold and $\Lambda \ga$ a $-1$-form. And
finally, for the volume form
$$ [L, \Lambda] dx\wedge dy = L \Lambda dx\wedge dy = L(-1) =
-dx\wedge dy$$
and the two dimensional case is done.

In higher dimensions, in appropriate coordinates, we have the splitting
$L=L_1+L_2$, $\Lambda= \Lambda_1+ \Lambda_2$ and if $\ga= \ga^1\wedge
\ga^2$ is a $k$-form ($k= k_1+k_2$) as before, we have
\begin{align*}
[L,\Lambda] \ga &= ([L_1, \Lambda_1] + [L_2, \Lambda_1] + [L_1,
\Lambda_2] + [L_2, \Lambda_2]) \ga^1 \wedge \ga^2\\
&= (n_1- k_1) \ga^1 \wedge \ga^2 + (n_2- k_2) \ga^1 \wedge \ga^2 = (n-k) \ga,
\end{align*}
where we have used that $L_1$ and $\Lambda_2$ commute, as well as
$L_2$ and $\Lambda_1$;

\item $[L, \delta] = - d$. Indeed, for a $k$-form \ga\ we have
\begin{align*}
[L,\delta]\ga &= L\delta \ga - \delta L \ga =  L \Lambda d \ga - d L
\Lambda \ga - \Lambda L d \ga + d\Lambda L \ga\\
&= [L, \Lambda] d \ga - d [L,\Lambda] \ga = \mbox{ (by item 2)}
= (m-k-1) d\ga - (m-k) d\ga = -d\ga;
\end{align*}

\item $\Lambda = -*L*$, which is proved again by induction. In
dimension two the unique nontrivial case is for $f v_M$, since both
$-*L*$ and $\Lambda$ have degree $-2$. Then, again letting $v_M = dx\wedge dy$
$$ - *L * v_M = - * L 1 = -* v_M = -1 = \Lambda(v_M).$$
The argument for the general case is the same as that already used;

\item $[\Lambda,\delta] =0$ and $[\Lambda,d] =
\delta$. In fact, the second claim is the very definition of
$\delta$, while the first follows from item {\it iv$\,$}) and Proposition
\ref{T:adjoint}, since when acting on $k$-forms
$$ [\Lambda,\delta] = - (-1)^{k+1}* L * * d * + (-1)^{k+1} * d * * L * =
(-1)^{k+1} *[d,L]* = 0;$$ 

\item $[L,H] = 2 L$ and $[\Lambda, H]= - 2 \Lambda$. This follows
from the expression for $H$ obtained in item 2.
\end{renumerate}

We summarize all of this in the next

\begin{prop}{commuting} {\em (Yan \cite{Ya96}):}
The operators $L$, $\Lambda$ and $H =
[L,\Lambda]$ satisfy the following commutation relations
\begin{alignat*}{3}
[L,d] &= 0 &\qquad  [L,\delta] &= -d &\qquad  L &= - * \Lambda *\\
[\Lambda,\delta] &= 0 &\qquad [\Lambda,d] &= \delta &\qquad \Lambda&= - * L *\\
H= \sum & (n-k) \Pi_k&\qquad [L,H] &= 2L &\qquad [\Lambda,H] &= -2 \Lambda
\end{alignat*}

Moreover, if we let
$$ X= \begin{pmatrix}
0& 1\\ 0&0
\end{pmatrix}, ~~
Y= \begin{pmatrix}
0& 0\\ 1&0
\end{pmatrix} \mbox{~~and~~}
Z= \begin{pmatrix}
-1& 0\\ 0&1
\end{pmatrix}$$
be a basis for $\frak{sl}(2,\C)$ and $\gf$ be defined by
$$\gf(X) = \Lambda,~~\gf(Y) = L \mbox{~~and~~} \gf(Z)= H$$
then we get a representation of $\frak{sl}(2,\C)$ in the space (of linear
operators in) $\gO^{\bullet}(M)$.
\end{prop}
With this, we have available all the representation theory of
$\frak{sl}(2,\C)$. We recall that
\begin{definition}
In a \frak{sl}(2,\C) representation \gf\ over a vector space $V$, a
vector $v \in V$ is said to be a {\it root} if $v$ is an eigenvector for
$\gf(Z)$. If besides that $\gf(X) v = 0$, $v$ is said to be a {\it
primitive root} or simply {\it primitive}. A representation of
$\frak{sl}(2,\C)$ is said to be of {\it finite spectrum} if
\begin{enumerate}
\item $\gf(Z)$ has only finitely many eigenvalues,
\item The space $V$ can be decomposed as a direct sum of $\gf(Z)$ eigenspaces
\end{enumerate}
\end{definition}

The expression for $H$ in Proposition \ref{T:commuting} shows that the
representation we are dealing with is of finite spectrum. And for
these representations we have the following result:

\begin{prop}{decomposition}
Any representation of $\frak{sl}(2,\C)$ of finite spectrum can be
decomposed as a direct product of irreducible finite dimensional
subspaces. Any irreducible subspace $W$ contains one and only one, up
to multiplicative constant, primitive root $v$ and $v, Yv, \cdots, Y^n
v$ form a basis for $W$, where, $n = \dim W-1$. All eigenvalues of
$\gf(Z)$ are integers.

Letting $V_k$ be the $k$-eigenspace of $V$, $\gf(Y^k):V_k \into
V_{-k}$ and $\gf(X^k):V_{-k} \into V_k$ are isomorphisms for every $k$
and the space $V_k$ is decomposed as
$$V_k = P_k \oplus \gf(Y)P_{k+2} \oplus \gf(Y) P_{k+4} \oplus\cdots,$$
where $P_k$ stands for the space of primitive roots in $V_k$. 
\end{prop}

\begin{remark}
One should notice that this is precisely the same theory used in
complex geometry, with $\Lambda$ replaced by $L^*$
to prove, for example, that K\"ahler manifolds satisfy the
Lefschetz property. See \cite{We73}.
\end{remark}

\begin{cor} (Yan \cite{Ya96}):
In a symplectic manifold $(M^{2n},\go)$, the maps
$$L^k:\gO^{n-k}(M) \into \gO^{n+k}(M)$$
$$\Lambda^k: \gO^{n+k}(M) \into \gO^{n-k}(M)$$
are isomorphisms.
\end{cor}

This $\frak{sl}(2,\C)$ representation can be restricted to the space of
symplectic harmonic forms. Indeed, if \ga\ is harmonic, by Proposition
\ref{T:commuting}, $dL\ga= Ld\ga =0$ and $\delta L \ga = L \delta \ga
+d \ga = 0$ and analogously for $\Lambda$. So we still have
representation theory available here.

\begin{cor}\label{T:harmonic representation} (Yan \cite{Ya96}):
In a symplectic manifold $(M^{2n},\go)$, the maps
$$L^k:H_{hr}^{n-k} \into H_{hr}^{n+k} ~~\mbox{ and }~~ \Lambda^k:
H_{hr}^{n+k} \into H_{hr}^{n-k}$$
are isomorphisms, where $H_{hr}^k$ is the set of symplectic harmonic
$k$-forms.
\end{cor}

\section{The Decomposition of Forms}\label{decomposition}

In this section we use Yan's operators $L$ and $\Lambda$ to give a concrete description of the subbundles $U^k$. With that, we can explicitly determine the splitting $d = \del + \delbar$, which shows a new subtle relation between these operators and suggests that, in some sense, when working with forms, degrees and the operators $d$ and $\delta$, one is actually using $\del$ and $\delbar$.

We start by describing the spaces $U^k$ for a symplectic vector space $(V,\go)$. Recall that $U^n$ is the line bundle defining the generalized complex structure, so, in this case, $U^n = \mbox{span}\{e^{i \omega}\}$. The next bundles are defined by
$$U^{n-k} = \wedge^k\overline{L}\cdot e^{i\go}= \overline{L} \cdot U^{n-k+1}.$$

For a symplectic structure $e^{i \go}$, the bundle $\overline{L}$ is given by
$$\overline{L} = \{X + i X\lfloor \go \mid X \in V\tensor \C\}.$$
Before we can give an explicit description of $U^k$ we need the following lemma:

\begin{lem}\label{T:symplectic decomposition}
For any vector $X \in V\tensor \C$ and complex $k$-form \ga\ the following identities hold:
$$\Lambda((X\lfloor \omega)\wedge \ga) = X\lfloor \ga + (X \lfloor \omega)\wedge \Lambda \alpha;$$
$$e^{\frac{\Lambda}{2i}} ((X\lfloor\go) \wedge \ga) = \frac{1}{2i}e^{\frac{\Lambda}{2i}} X\lfloor \ga + (X\lfloor \go) e^{\frac{\Lambda}{2i}}\ga.$$
\end{lem}
\begin{proof}
We start with the first identity. It is enough to take Darboux coordinates so that \go\ is standard and check for $X = \del_{x_i}$ and $\del_{y_i}$. As both cases are similar we will do only the first.

Write $\ga = \ga_0 + dx_i\ga_x + dy_i \ga_y + dx_1\wedge dy_i \ga_{xy}$. The left hand side is
\begin{align*}
\Lambda((X\lfloor \omega)\wedge \ga) &= \Lambda(dy_i \wedge \alpha) = \Lambda(dy_i\ga_0+ dy_i \wedge dx_i \ga_x)\\
&=dy_i \Lambda \ga_0 + \ga_x +dy_i\wedge dx_i \Lambda \ga_{x}
\end{align*}
And the right hand side is
\begin{align*}
X\lfloor \ga + (X\lfloor \omega)\wedge \Lambda \alpha &=\ga_x + dy_i\ga_{xy} + dy_i(\Lambda \ga_0 + dx_i \Lambda \alpha_x + dy_i \Lambda \ga_y - \ga_{xy} + dx_i \wedge dy_i \Lambda \ga_{xy})\\
&= \ga_x + dy_i \ga_{xy} + dy_i \Lambda \ga_0 + dy_i \wedge dx_i \Lambda \ga_x - dy_i \ga_{xy}.
\end{align*}
So, the first identity follows.

By induction, from the first identity, we get that
$$\Lambda^k(X\lfloor \go \wedge \ga) = k \Lambda^{k-1} X\lfloor\ga + (X\lfloor \go) \wedge\Lambda^k \ga .$$
Therefore, by expanding the exponential in Taylor series, we obtain the second identity.
\end{proof}

\begin{theo}{symplectic U_k}
The decomposition of $\wedge^{\bullet}V\tensor \C$ for a symplectic vector space $V$ is given by
$$U^{n-k} = \{e^{i\go}(e^{\frac{\Lambda}{2i}}\alpha) \mid \alpha \in \wedge^k V\}.$$
Hence, the natural isomorphism
$$\gf:\wedge^{\bullet}V \tensor \C \into \wedge^{\bullet}V \tensor \C \qquad \gf(\ga) = e^{i\go} e^{\frac{\Lambda}{2i}}\ga,$$
is such that $\gf: \wedge^k V \cong U^{n-k}$.
\end{theo}
\begin{proof}
This is done by induction. For $\alpha$ a 0-form the expression above agrees with $U^n$. If $U^{n-k}$ is as described above, then $U^{n-k-1}= \overline{L}\cdot U^{n-k}$, so its elements are linear combinations of terms of the form, with $\ga \in \wedge^k V \tensor \C$,
\begin{align*}
(X + i X \lfloor \go)e^{i\go}(e^{\frac{\Lambda}{2i}}\ga)&=(2i (X\lfloor \omega) \wedge  e^{\frac{\Lambda}{2i}}\ga + X\lfloor e^{{\Lambda}{2i}}\ga)e^{i\omega}\\
&=(-e^{\frac{\Lambda}{2i}}X\lfloor \ga + 2i e^{\frac{\Lambda}{2i}}(X\lfloor\go \wedge\ga) + e^{\frac{\Lambda}{2i}}X\lfloor\ga)e^{i\go}\\
&=2i e^{i\go} e^{\frac{\Lambda}{2i}}((X\lfloor \omega)\wedge \ga),
\end{align*}
where the second equality follows from the previous lemma. Now, since $\ga$ can be chosen to be any $k$-form and $\go$ is nondegenerate, the space generated by the forms above is $\{e^{i\go}e^{\frac{\Lambda}{2i}}\beta\mid \beta \in \wedge^{k+1}V\tensor \C\}$, and the theorem is proved.
\end{proof}

In general we have the identity
$$ \ga = e^{i\go}e^{\frac{\Lambda}{2i}}\left(e^{\frac{-\Lambda}{2i}}e^{-i\go} \ga \right),$$
and hence, letting $\ga_k$ be the component of $e^{\frac{-\Lambda}{2i}}e^{-i\go}\ga$ of degree $k$ we obtain the decomposition of the form \ga\ into its $\mc{U}^k$ components:
$$\ga = \sum_k e^{i\go}e^{\frac{\Lambda}{2i}} \ga_k.$$
There will be no need to find the explicit expression for the $\ga_k$.

Now we move on to the operators \del\ and \delbar.

\begin{theo}{del and delbar}
For any form $\ga$,
$$d(e^{i\go}e^{\frac{\Lambda}{2i}} \ga) = e^{i\go}e^{\frac{\Lambda}{2i}} (d \ga -\frac{1}{2i} \delta \ga).$$ 
Therefore
$$\del(e^{i\go}e^{\frac{\Lambda}{2i}} \ga) = e^{i\go}e^{\frac{\Lambda}{2i}} d \ga$$
$$\delbar(e^{i\go}e^{\frac{\Lambda}{2i}} \ga) = - e^{i\go}e^{\frac{\Lambda}{2i}} \frac{1}{2i}\delta \ga.$$
Hence, the natural isomorphism $\gf$ of Theorem \ref{T:symplectic U_k} is such that
$$\gf(d \ga) = \del \gf(\ga) \qquad \mbox{and}\qquad \gf(\delta \ga) = -2i \delbar \gf(\ga).$$
\end{theo}
\begin{proof}
By definition of $\delta$, $d\Lambda = \Lambda d - \delta$. Then, by induction, and using that $\delta$ and $\Lambda$ commute, $d\Lambda^{k} = \Lambda^k d - k \Lambda^{k-1} \delta$.
Therefore,
\begin{align*}
d(e^{i\go}e^{\frac{\Lambda}{2i}} \ga)&= e^{i\go} (e^{\frac{\Lambda}{2i}} \ga)
=e^{i\go} \sum d\left(\frac{\Lambda^k}{(2i)^k k!} \ga\right)\\
&=e^{i\go} \sum \left(\frac{\Lambda^k}{(2i)^k k!} d\ga - \frac{1}{2i} \frac{\Lambda^{k-1}}{(2i)^{k-1} (k-1)!} \delta\ga\right)\\
&=e^{i\go} e^{\frac{\Lambda}{2i}}(d\ga -\frac{1}{2i}\delta \ga).
\end{align*}
The rest of the theorem follows from the fact that \del\ and \delbar\ are the projections of $d$ onto $\mc{U}^{k+1}$ and $\mc{U}^{k-1}$ respectively.
\end{proof}


As a first application of this theorem we remark that the canonical spectral sequence (see Chapter \ref{ddj-lemma}, Section \ref{canonical spectral sequence}) always degenerates at $E_1$ for a symplectic structure. Indeed, the $E_1$ term is just the $\del$-cohomology and the canonical spectral sequence converges to the $d$-cohomology. But the expression above shows that these two cohomologies are isomorphic.

\begin{ex}{Upq to omegapq}
As a second application we relate the decomposition of forms into $\mc{U}^{p',q'} = \mc{U}_{\go}^{p'} \cap \mc{U}_J^{q'}$ in a K\"ahler manifold (see page \pageref{T:upq}) with the standard decomposition into $\gO^{p,q}$. Observe that although the map \gf\ introduced in Theorem \ref{T:symplectic U_k} does not preserve degrees, it does preserve the difference $p-q$, since in this case $\go$ is of type $(1,1)$. Therefore, this map preserves the $\mc{U}_J^{q'}$ (see Example \ref{T:complex decomposition}). As the elements in $\gO^{p,q}$ have degree $p+q$, $\gf(\gO^{p,q}) \in \mc{U}_{\go}^{n-p-q}$. Thus,
\begin{equation}\label{E:Upq to omegapq}
\gf(\gO^{p,q}) = \mc{U}^{n-p-q,p-q}
\end{equation}
gives the generalized decomposition from the standard $\gO^{p,q}$ decomposition.
\end{ex}

\section{Lefschetz Property, Harmonic Representatives and the $d\delta$-lemma}\label{ddelta}

In this section we present both Yan's and Merkulov's theorem, using the $\frak{sl}(2,\C)$ representation from Section \ref{sl2}. We complete Merkulov's asserted proof of his theorem with Lemmata \ref{T:lemma2c3} and \ref{T:lemma3} and Proposition \ref{T:weak ddelta lemma} (Merkulov's argument has already been used to prove Theorem \ref{T:merkulov}, hence it is not present in this section). These theorems state that the existence of symplectic harmonic representatives in each cohomology class, the $d\delta$-lemma and the Lefschetz property are equivalent.

To prove Yan's theorem, we let $\tilde{H}^k_{hr}(M)=H^k_{hr}/\im d \cap H^k_{hr}$. Then
$\tilde{H}^k_{hr}(M)$ is a subspace of $H^k(M)$.

\begin{lem}\label{T:yan} (Yan \cite{Ya96})
If \ga\ is a closed $n-k$-form for which $L^{k+1}\ga$ is exact, then
there is a symplectic harmonic form in the same cohomology class as \ga.
\end{lem}
\begin{proof}
By hypothesis, there exists $\tilde{\beta}$ such that $L^{k+1}\ga
= d\tilde{\beta}$. Then, since $L^{k+1}:\gO^{n-k-1}(M) \into
\gO^{n+k+1}(M)$ is an isomorphism, we can find $\beta$
satisfying $L^{k+1} \beta = \tilde{\beta}$. Letting $\xi = \ga - d
\beta$, we see that $\xi$ is closed and in the same cohomology class as
\ga. Moreover $L^{k+1}\xi = L^{k+1}\ga - L^{k+1}d\beta = 0$, thus
$\xi$ is (closed and) primitive and hence symplectic harmonic.
\end{proof}

\begin{theo}{yan2}{\em (Yan \cite{Ya96}):}
The following are equivalent for a symplectic manifold $(M^{2n},\go)$
\begin{enumerate}
\item $M$ satisfies the Lefschetz property;

\item There is a symplectic harmonic representative in each cohomology class.
\end{enumerate}
\end{theo}
\begin{proof} Assume that (2) holds and
consider the following commutative diagram
$$
\begin{array}{ccc}
 H^{n-k}_{hr}& \stackrel{L^k}{\longrightarrow} & H^{n+k}_{hr} \\
\Big\downarrow & & \Big\downarrow \\
 H^{n-k}& \stackrel{L^k}{\longrightarrow}& H^{n+k} \\
\end{array}
$$
Then the two vertical arrows are surjective by hypothesis, and the
top horizontal arrow is an isomorphism by the previous corollary, therefore the bottom arrow is surjective and (1) holds.

Now we assume that $M$ satisfies the Lefschetz property. Then initially we observe that there is a splitting 
\begin{equation}\label{E:primitive splitting}
H^{n-k}(M) = \im L+ P_{n-k},
\end{equation}
where $P_{n-k} = \{[\ga]| L^{k+1} [\ga] = 0\}$ and $\im L$ is the image of $L:H^{n-k-2}(M) \into H^{n-k}(M)$. 
Indeed, given a closed $n-k$-form $\ga$, letting
$\beta= L^{k+1} \ga$, $\beta$ is a closed $n+k+2$-form. Thus, by
the Lefschetz property, there is a closed $n-k-2$-form $\gamma$ such
that
$$0 = [L^{k+2} \gamma - \beta] = [L^{k+2}\gamma - L^{k+1}\ga] =
L^{k+1}[\go \gamma - \ga].$$
Therefore, $\ga = L \gamma + (\ga - \go \gamma).$ The first term is in
$\im L$ and the second is in $P_{n-k}$.

Now we proceed by induction. Any closed 0- or 1-form is symplectic
harmonic so there is nothing to be done in these cases. Assuming
that every cohomology class in $H^j(M)$ has a symplectic harmonic
representative for $j < n-k$, we are going to prove the same for
$n-k$. Let $\ga$ be an 
$n-k$-closed form, then we can decompose $\ga = L \gamma +
\tilde{\ga}$, where $[L^{k+1}\tilde{\ga}] =0$. By hypothesis, we can
find $\tilde{\gamma}$ symplectic harmonic in the same cohomology class
as $\gamma$ and, by the previous lemma, we can find $\xi$ symplectic
harmonic in the same cohomology class as $\tilde{\ga}$.
\end{proof}

Now we want to use Theorem \ref{T:merkulov} to conclude that the Lefschetz property implies the $d\delta$-lemma. We provide next in Lemmata \ref{T:lemma2c3} and \ref{T:lemma3} and Proposition \ref{T:weak ddelta lemma} the parts missing to make Merkulov's proof complete.

\begin{lem}\label{T:lemma2c3}
There is a
nonvanishing constant $C_{j,n-k}$ such that if \ga\ is a primitive $(n-k)$-form then for all $j \leq k$, $\Lambda^j L^j \ga = C_{j,n-k} \ga$.
\end{lem}

\begin{proof}
This is a consequence of the representation of $\frak{sl}(2)$, but can also be obtained by using the relation $[L,\Lambda] = \sum (n-k) \Pi_k$ and induction. 
\end{proof}

\begin{lem}\label{T:lemma3}
If $\ga = \sum L^j \ga_j$  then
$$d\delta\ga = \sum L^j d\delta\ga_j.$$
Moreover, if $\ga_j$ is primitive so is $d\delta\ga_j$.
\end{lem}
\begin{proof}
Since we know that $d\delta$ is linear (over \RR), we
only have to prove that
$$d\delta L^r \ga = L^r d\delta \ga$$
But in this case, using Proposition \ref{T:commuting}, we have
$$d\delta L^r \ga = d(L\delta +d)L^{r-1} \ga = dL\delta
L^{r-1}\ga = \cdots = L^rd\delta\ga.$$

Finally, if \ga\ is a primitive $k$-form then $L^{n-k+1}\ga =0$, hence
$$ L^{n-k+1} d \delta \ga = d \delta L^{n-k+1} \ga = 0$$
and consequently $d \delta \ga$ is primitive.
\end{proof}

\begin{prop}{weak ddelta lemma}
If $(M,\go)$ is compact and the Lefschetz property holds then
$$\im \delta \cap \ker d  =\im d \cap \im \delta,$$
$$\im d \cap \ker \delta  =\im d \cap \im \delta.$$
\end{prop}
\begin{proof}
We will prove only the first identity as the second is analogous and also follows as an easy consequence of the first. The result can be restated in the following way: if $d\delta\ga=0$, then $\delta\ga$ is exact. We will prove it by induction on the degree of \ga.

Assume \ga\ is a 0-form. Then $d\delta \ga=0$
and $\delta \ga =0 = d0$, so $\delta \ga$ is exact for every 
0-form.

Now, let \ga\ be a 1-form with $d\delta \ga =0$. Then $\delta\ga$ is
closed 0-form, and consequently constant: $\delta\ga = c$. This means
$$ c = \delta \ga = -*d*\ga$$
Applying $*$ to both sides, and using that $** = \Id$ (Lemma \ref{T:symplectic*}) we get
$$c v_M = -d*\ga,$$
which can not happen in a compact manifold.

For a general $k$-form $\ga$ with $d\delta\ga =0$, decompose $\ga =
\sum L^r \ga_r$ with $\ga_r$ primitive, then, by Lemma \ref{T:lemma3}
$$0 = d\delta \ga = \sum L^r d\delta \ga_r,$$
with $d\delta \ga_r$ primitive, hence this is a direct sum and each
term must vanish, hence we must have $d \delta \ga_r =0$ for every
$r$.

By the induction hypothesis, $\delta \ga_r = d \gf_r$ for $r>0$ and so
$$ \delta L^r \ga_r = (L\delta + d) L^{r-1}\ga_r = \cdots = r d
L^{r-1} \ga_r + L^r\delta \ga_r = d(r L^{r-1}\ga_r + L^r\gf_r).$$

Now we are left with the case of \ga\ a primitive $k$-form. For this,
define $\beta$ by
\begin{equation}\label{E:definitionofbeta}
L^{n-k+1}\beta = d L^{n-k} \ga.
\end{equation}
Then $L^{n-k+2} \beta = d L^{n-k+1}\ga =0$, consequently $\beta$ is
primitive, hence, applying $\Lambda^{n-k+1}$ and using Lemma
\ref{T:lemma2c3}, we get
\begin{align*}
C_{n-k+1,k-1} \beta &= \Lambda^{n-k+1}dL^{n-k} \ga =
\Lambda^{n-k}(d\Lambda -\delta)L^{n-k} \ga = \cdots\\
&=\big(\overset{0,~\ga~\text{primitive}}{\overbrace{d\Lambda^{n-k+1}}} -(n-k+1) \delta \Lambda^{n-k}\big)L^{n-k}
\ga\\
&=\delta C_{n-k,k}\ga.
\end{align*}
Applying $L^{n-k+1}$ to both sides and using \eqref{E:definitionofbeta} we obtain that, for a 
suitable constant $C$, 
$$C L^{n-k+1} \delta  \ga = L^{n-k+1}\beta =  d(L^{n-k}\ga).$$
Since by hypothesis $\delta \ga$ is closed, Lefschetz
implies that $\delta \ga$ is exact.
\end{proof}

\begin{theo}{ddelta lemma} {\em (Merkulov \cite{Me98}):} A compact symplectic manifold satisfies the Lefschetz property if and only if it satisfies the $d\delta$-lemma.
\end{theo}
\begin{proof}
If the Lefschetz property holds, then, by Theorem \ref{T:yan2}, each cohomology class has a $\delta$-closed representative. Further, according to Proposition \ref{T:weak ddelta lemma}, $\im d \cap \ker \delta = \im \delta \cap \ker d$. Applying the map \gf\ from Theorem \ref{T:del and delbar}, $d$ and $\delta$ change to \del\ and \delbar\ and, by Theorem \ref{T:merkulov}, the manifold  satisfies the $d\delta$-lemma.

Conversely, if the manifold satisfies the $d\delta$-lemma, by Theorem \ref{T:decomposition2}, each cohomology class can be represented by a $\delta$-closed form, i.e., a symplectic harmonic form and thus Theorem \ref{T:yan2} implies that the manifold satisfies the Lefschetz property.
\end{proof}

\chapter{The Symplectic Blow-up and the Lefschetz Property}\label{blowup}

It is a result of Parshin \cite{Pa66} that the $dd^c$-lemma is preserved by rational equivalence, in particular, by the blow-up along a complex submanifold. Although the blow-up construction does not seem to have an analogue in generalized complex geometry, there is a {\it symplectic} blow-up, as introduced by McDuff \cite{McD84}. This operation is similar to the complex blow-up from a topological point of view, but different in many senses from the differential geometric point of view. For example, there is a canonical complex structure in the complex blow-up, while there is no canonical choice of symplectic structure in the symplectic counterpart.

One further difference is the behaviour of the symplectic $d\delta$-lemma which, as we have seen, is equivalent to the Lefschetz property (cf. Theorem \ref{T:ddelta lemma}): In her paper, McDuff introduces the symplectic blow-up and uses it to give the first simply connected example of a non-K\"ahler symplectic manifold, which happens to be the blow-up of $\C P^5$ along a symplectically embedded Thurston manifold (the nilmanifold with structure $(0,0,0,12)$). This example of McDuff fails to be K\"ahler, amongst other reasons, because it does not satisfy the Lefschetz property. This means that even if the ambient space satisfies the $d\delta$-lemma, the blown-up manifold may not do the same. This is also related to the fact that symplectic submanifolds are more flexible than complex submanifolds. For example, any complex submanifold of $\C P^n$ is K\"ahler (algebraic), whereas any symplectic manifold can have its structure changed so it embeds in $\C P^n$ \cite{Gr71,Ti77} which shows that there are no topological obstructions to the symplectic embedding.

One of the purposes of this chapter is to study systematically how the Lefschetz property behaves under the blow-up, in particular we seek conditions under which we can assure that the blown-up manifold will satisfy the Lefschetz property. We prove that this is the case if both submanifold and ambient manifold satisfy the Property (cf. Theorem \ref{T:lefschetz}). Moreover we study the blow-down map and show that even if the blown-up manifold satisfies the Lefschetz property, the original ambient manifold will not necessarily do so (cf. Theorem \ref{T:lefschetz in m2} and Proposition \ref{T:i=2d}), which, together with McDuff's example, shows that we can not decide whether the blow-up will or will not satisfy the Lefschetz property based solely on the ambient manifold. 

The second purpose of this chapter is to answer the question of the correlation between the $d\delta$-lemma and formality. One of the implications has been known to be false for some time, as Gompf produced a simply-connected symplectic  6-manifold which does not satisfy the Lefschetz property \cite{Go95}, but which is formal by Miller's result \cite{Mi79}. The converse implication had been conjectured by Babenko and Taimanov \cite{BT00a} and was the object of study of other papers \cite{IRTU03,LO94}.

Our starting point in this case is that Babenko and Taimanov studied thoroughly the behaviour of Massey products under blow-up \cite{BT00a,BT00b} and these tend to `survive' in the blow-up, which is markedly different from the behaviour of the Lefschetz property. Therefore, we produce an example of a nonformal symplectic manifold satisfying the Lefschetz property by blowing-up a nilmanifold along a suitable torus. We also produce a 4-dimensional example using Donaldson submanifolds \cite{Do96} and results of Fern\'andez and Mu\~noz \cite{FM02}.

This chapter is organized as follows. In the first section we explain briefly how the blow-up is done in symplectic geometry and derive the cohomology ring (also at the form level) of the blow-up based on the ambient manifold, submanifold and Chern and Thom classes of the normal bundle of the submanifold. In Section \ref{S:lefschetz}, we study how the Lefschetz property behaves under blow-up, initially in the case of the blow-up along an embedded surface and later the general case. In Section \ref{S:massey} we recall the definition of Massey products, their relation with formality and their behaviour under blow-up and finish in Section \ref{S:nonformal examples} with examples of nonformal symplectic manifolds satisfying the Lefschetz property.

The results of this chapter were presented first in \cite{Ca04}.

\section{The Symplectic Blow-up}\label{S:blowup}

We begin by giving a description of the cohomology ring of the blown-up
manifold in terms of the cohomology rings of the ambient manifold and the
embedded submanifold and the Chern and Thom classes of the normal
bundle of the embedding. We shall outline the blow-up construction in
order to fix some notation. For a detailed presentation we refer to
\cite{McD84}.

Assume that $i:(M^{2d},\sigma) \hookrightarrow
(X^{2n},\go)$ is a symplectic
embedding, with $M$ compact. Let $k=n-d$. In these circumstances we can 
choose a
complex structure in $TX$ that restricts to one in $TM$, and hence
also to the normal bundle $E
\stackrel{\pi}{\rightarrow} M$. Therefore $E$ is a complex
bundle over $M$ and one can form its projectivization
$$ \mathbb{C}P^{k-1} \longrightarrow \tilde{M} \longrightarrow M$$
and also form the ``tautological'' line bundle $\tilde{E}$ over $\tilde{M}$:
 the subbundle of $\tilde{M}\times E$ whose fibers are the
elements $\{([v],\lambda v), \lambda \in \mathbb{C}\}$. We have
the following commutative diagram
\begin{equation}\label{E:diagram}
\begin{array}{ccccc}
\tilde{E_0} & \longrightarrow &\tilde{E} &
\stackrel{q}{\longrightarrow}& \tilde{M}\\
\Big\downarrow \vcenter{\rlap{$\gf$}}& & \Big\downarrow \vcenter{\rlap{$\gf$}}& &
\Big\downarrow \vcenter{\rlap{$p$}}\\
E_0 & \longrightarrow &(E,\go) &
\stackrel{\pi}{\longrightarrow}& (M,\sigma)
\end{array}
\end{equation}
where $q$ and \gf\ are the projections over $\tilde{M}$ and $E$
respectively, $E_0$ is the complement of the zero section in $E$ and
$\tilde{E}_0$ the complement of the zero section in $\tilde{E}$.

It is easily seen that $E_0$
and $\tilde{E}_0$ are diffeomorphic via \gf. Furthermore, if we let $V$
be a sufficiently small disc subbundle in $E$ with its canonical
symplectic structure \go, then it is symplectomorphic to a neighbourhood
of $M \subset X$ and we identify the two from now on. Letting
$\tilde{V} = \gf^{-1}(V)$,
we can form the manifold
$$ \tilde{X} = \overline{X-V} \cup_{\del V}\tilde{V}.$$
Then, the map \gf\ can be extended to a map  $f:\tilde{X}
\into X$, being the identity in the complement of $\tilde{V}$. The manifold $\tilde{X}$ is the blow-up of $X$ along $M$ and $f:\tilde{X}
\into X$ is the projection of the blow-up.

\begin{lem} (McDuff \cite{McD84})
There is a unique class $a \in H^2(\tilde{M})$ which restricts to the standard K\"ahler class on each fiber of $\tilde{M} \into M$ and pulls back to the trivial class in $\tilde{E}_0$. Moreover, $H^{\bullet}(\tilde{E}) \cong H^{\bullet}(\tilde{M})$ is a free module over $H^{\bullet}(M)$ with generators $1,a,\cdots,a^{k-1}$. 
\end{lem}

\begin{theo}{cohomology}
{\em (McDuff \cite{McD84})} If the codimension of $M$ is at least 4, the
fundamental groups of $X$ 
and the blown up manifold $\tilde{X}$ are isomorphic. Further, there
is a short exact sequence
$$0 \into H^*(X) \rightarrow H^*(\tilde{X}) \into A^* \into 0,$$
where $A^*$ is the free module over $H^*(M)$ with generators $a,
\cdots, a^{k-1}$, with $a$ as in the previous lemma.
Moreover, there is a representative $\ga$ of $a$ with support in the 
tubular neighbourhood $V$ such that, for \e\ small enough, the form 
$\tilde{\go} = f^*(\go) + \e \ga$ is a symplectic form in $\tilde{X}$.
\end{theo}

\begin{remark}
As is observed by McDuff \cite{McD84}, the Leray-Hirsch theorem implies that $a^k$ is related to $a, \cdots, a^{k-1}$ in $\tilde{E}$ by
$$
a^k = -c_k -c_{k-1} a + \cdots - c_1 a^{k-1},
$$
where the $c_j$'s are the Chern classes of the normal bundle $E$.

In \cite{RT00} it is shown that in $\tilde{X}$ this relation becomes
$$
a^k = -f^*(t) -c_{k-1} a + \cdots - c_1 a^{k-1},
$$
where $t$ is the Thom class of the embedding $M \hookrightarrow X$ and $f:\tilde{X} \into X$ the projection of the blow-up.
\end{remark}

With this, we have a complete description of the cohomology ring of
$\tilde{X}$. For $v_1, v_2 \in H^*(X)$ and $u_1, u_2 \in H^*(M)$,
\begin{equation}\label{E:product rules}
\begin{cases}
f^*(v) \wedge f^*(w) &= f^*(v\wedge w);\\
f^*(v) a &= i^*(v) a;\\
a u_1 \wedge a u_2 &= a^2 u_1 \wedge u_2\\
a^k &= -f^*(t) -c_{k-1} a + \cdots - c_1 a^{k-1};\\
f^*(t) \wedge u_1 &=  f^*(t \wedge u_1), \mbox{ the Thom map extended 
to } X.
\end{cases}
\end{equation}

\section{The Lefschetz Property and Blowing-up}\label{S:lefschetz}

Now we move on to study how the Lefschetz property behaves under blow-up. The first case to look at would be the blow-up of a point, but, as
we will see, this does not change the kernel of the Lefschetz map at
any level (cf. Theorem \ref{T:i>2d}). The next case would be a
surface. Here, on the one hand, the situation is simple enough for us
to be able to give a fairly complete account of what happens, and, on the other, we can already see that in this case it is possible to decrease the
dimension of the kernel of the Lefschetz map.

\subsection{Blowing up along a Surface}
Assume that $i:(M^2,\sigma) 
\hookrightarrow
(X^{2n},\go)$ is 
a connected surface symplectically embedded in $X$, $M$ and $X$ are compact, and
let $\tilde{X}$ be the 
blow-up of $X$ along $M$. In $H^1(\tilde{X})$ things go as follows
\begin{equation}\label{E:h1}
\begin{aligned}
(f^*(\go) + \e a)^{n-1} f^*(v) &= f^*(\go^{n-1}v) + \e^{n-1}
a^{n-1}i^*v \\
&= f^*(\go^{n-1}v - \e^{n-1} t v)
\end{aligned}
\end{equation}
and if Lefschetz holds for $X$ and $\e$ is small enough Lefschetz will
also hold for $\tilde{X}$, or, more generally,
$\dim(\ker(\tilde{\go}^{n-1})) \leq  
\dim(\ker(\go^{n-1}))$. Now we proceed to show that in certain 
conditions the inequality holds.

\begin{lem}\label{T:lemma1}
Let $i: (M^2,\sigma) \hookrightarrow (X^{2n},\go)$ be a symplectic
embedding, $M$ and $X$ be compact and $t$ be the Thom class of this embedding. The following
are equivalent: 
\begin{enumerate}
\item There are $v_1, v_2 \in H^1(X)$ in $\ker(\go^{n-1})$
such that $i^*(v_1\wedge v_2) \neq 0$;
\item There exists $v_1 \in \ker(\go^{n-1})$ such that $t\wedge v_1 \not
\in \im(\go^{n-1})$.
\end{enumerate}
\end{lem}
\begin{proof}
Assuming (1), by the defining property of the Thom class,
$$ \int_X t \wedge v_1 \wedge v_2 = \int_M i^*(v_1 \wedge v_2) \neq 0,$$
and, since both $v_1$ and $v_2$ pair trivially with $\im(\go^{n-1})$,
but pair nontrivially with $t \wedge v_i$, we
see that $t \wedge v_i \not \in \im(\go^{n-1})$. So (1) implies (2).

On the other hand, assume that there is a $v_1$ satisfying (2). Let $\{a_i\}$ be a
basis for $\ker (\go^{n-1})$ and $\{\tilde{a}_i\}$ be a basis for a
complement. Since $(H^1)^* \cong H^{2n-1}$, we can view the dual basis
$\{a_i^*,\tilde{a}_i^*\}$ as a basis for $H^{2n-1}$. Then we note that
$\im(\go^{n-1}) \subset \mbox{span}\{\tilde{a}_i^*\}$, and since these spaces
have the same dimension they are the same. Therefore, the condition
$t\wedge v_1 
\not \in \im(\go^{n-1})$ implies that it pairs nontrivially with some
of the $a_i$. Let  $v_2$ be such an $a_i$. Then again by the defining
property of the Thom class we have
$$\int_M i^*(v_1 \wedge v_2) = \int_X t \wedge v_1 \wedge v_2 \neq 0,$$
and $i^*(v_1\wedge v_2) \neq 0$.
\end{proof}

\begin{lem}\label{T:h1}
If the equivalent conditions (1) and (2) of the previous lemma are
satisfied and \e\ is 
small enough, then
$$ \dim(\ker(\tilde{\go}^{n-1}:H^1(\tilde{X}) \into
H^{2n-1}(\tilde{X}))) \leq \dim(\ker(\go^{n-1}:H^1(X) \into
H^{2n-1}(X))) - 2.$$
\end{lem}
\begin{proof}
Let $V$ be a complement of $\ker(\go^{n-1})$ in $H^1(X)$ and
$v_1$ and 
$v_2$ the cohomology classes satisfying condition (1) of 
lemma \ref{T:lemma1}. Then, since neither $t \wedge v_1$ or $t \wedge
v_2$ is in 
$\im(\go^{n-1})$, for \e\ small enough, equation \eqref{E:h1}
shows that
$f^*(t \wedge v_i) \not \in \im(\tilde{\go}|_V)$, 
since
$\tilde{\go}^{n-1}|_V$ is simply a perturbation of the injection 
$\go^{n-1}|_V$. On the other hand, 
$\tilde{\go}^{n-1} v_i = -\e^{n-1} f^*(t\wedge v_i)$, and therefore $f^*(t
\wedge v_i)$ is in the image of $\tilde{\go}^{n-1}$, so
$$\dim(\im(\tilde{\go}^{n-1})) \geq \dim(\im(\go^{n-1})) + 2$$ 
and the result follows.
\end{proof}

Now we move on to $H^2(\tilde{X})$, where we have
\begin{equation}\label{E:h2}
\begin{aligned}
(f^*(\go) + \e a)^{n-2} (f^*(v_2) + a v_0) =&\\
= f^*(\go^{n-2}v_2  -  \e^{n-2} t v_0)& + \e^{n-3}a^{n-2}((n-2) \sigma v_0
 + \e( i^*v_2 - c_1 v_0))
\end{aligned}
\end{equation}
and then we observe that the map above is a perturbation of
$$f^*(v_2) + a v_0 \mapsto f^*(\go^{n-2}v_2) + \e^{n-3}a^{n-2}(n-2) \sigma v_0.$$
Therefore for \e\ small enough, Lefschetz will hold for $\tilde{X}$
if it holds for $X$, or more generally $\dim(\ker(\tilde{\go}^{n-2}))
 \leq \dim(\ker(\go^{n-2}))$.

Again, we may have the inequality.
\begin{lem}\label{T:lemma2}
Let $i: (M^2,\sigma) \hookrightarrow (X^{2n}\go)$ be a symplectic
embedding, $M$ and $X$ be compact and $t$ be the Thom class of this embedding. The following
are equivalent: 
\begin{enumerate}
\item There exists $v \in \ker(\go^{n-2}:H^2(X) \into H^{2n-2}(X))$ such that $i^*v \neq 0$;
\item The Thom class $t$ is not in the image of $\go^{n-2}$.
\end{enumerate}
\end{lem}
\begin{proof}
The proof is the same as the one for lemma
\ref{T:lemma1}. Assuming (1), by definition of the Thom class,
$$ \int_X t v = \int_M i^*v \neq 0.$$
On the other hand any vector in the kernel of $\go^{n-2}$ pairs
trivially with $\im(\go^{n-2})$, so $t\not \in \im(\go^{n-2})$.

Conversely, we let again $\{a_i\}$ be a basis for $\ker(\go^{n-2})$,
$\{\tilde{a}_i\}$ a basis for a complement and
$\{a^*_i,\tilde{a}^*_i\}$ the dual basis and again identify the
dual space with $H^{2n-2}$. Then we see that $\im(\go^{n-2}) =
\mbox{span}\{\tilde{a}^*_i\}$ and, since $t \not \in \im(\go^{n-2})$,
$t$ must pair nontrivially with at least one of the $a_i$'s. Call it $v$.
\end{proof}

\begin{lem}\label{T:h2}
If the equivalent conditions (1) and (2) of the previous lemma are
satisfied and \e\ is small enough then 
$$ \dim(\ker(\tilde{\go}^{n-2}:H^2(\tilde{X}) \into
H^{2n-2}(\tilde{X}))) = \dim(\ker(\go^{n-2}:H^2(X) \into
H^{2n-2}(X))) - 1$$
\end{lem}
\begin{proof}
By conveniently choosing $v_2$ and $v_0$ in \eqref{E:h2},
\begin{equation}\tag{\ref{E:h2}}
\begin{aligned}
(f^*(\go) + \e a)^{n-2} (f^*(v_2) + a v_0) =&\\
= \underset{H^{2n-2}(X)}{\underbrace{f^*(\go^{n-2}v_2 - \e^{n-2} t v_0)}}& + \e^{n-3}a^{n-2}((n-2)
\sigma v_0 + \e( \underset{\neq 0}{\underbrace{i^*v_2}} - c_1 v_0))
\end{aligned}
\end{equation}
the term in $H^{2n-2}(X)$ can be made equal to any
pre-chosen element in $\im(\go^{n-2}) \oplus \mbox{span}\{t\}$. Once
$v_2$ and $v_0$ are chosen, changing $v_2$ by an element in
$\ker(\go^{n-2})$ does not affect the result. On the other hand, by
varying $v_2$ by an element in $\ker(\go^{n-2})$ the coefficient of $a$ can be made
equal to anything in $H^{2}(M)$. Therefore $\dim(\im(\tilde{\go}^{n-2})) =
\dim(\im(\go^{n-2})) + 2$ and $\dim(H^2(\tilde{X})) = \dim(H^2(X)) +
1$, hence the result follows.
\end{proof}

Finally, we finish the study of the blow-up along surfaces claiming
that, for $i>2$, 
$$\dim(\ker(\tilde{\go}^{n-i})) = \dim(\ker({\go}^{n-i})).$$
Indeed, if $v_i \in \ker(\go^{n-i})$ then, $i^*(v_i) =0$, since it has
degree greater than 2, and therefore $a f^*(v_i) =a i^*(v_i) =0$ and 
$$(f^*(\go)+\e a)^{n-i}f^*(v_i)  = f^*(\go^{n-i} v_i) =0,$$
so $f^*(\ker(\go^{n-i})) \subset \ker(\tilde{\go}^{n-i})$.

Conversely, assuming $i$ even (the odd case is analogous),

\begin{align*}
(f^*(\go)+\e a)^{n-i}&(f^*(v_{i}) + a^{\frac{i}{2}} v_0 +
 a^{\frac{i-2}{2}}v_2) =\\
 &f^*(\go^{n-i}v_i) +
\e^{n-i-1} a^{n-\frac{i}{2}-1} ((n-i)\sigma v_0 + \e v_2) + \e^{n-i}
a^{n-\frac{i}{2}}v_0,
\end{align*}
and therefore $f^*(v_i) + a^{\frac{i}{2}} v_0 + a^{\frac{i-2}{2}}v_2$ will be in
$\ker(\tilde{\go}^{n-i})$ if, and only if, $v_0=0$ (by the coefficient
of $a^{n-\frac{i}{2}}$), $v_1 =0$ (by the coefficient of $a^{n-\frac{i}{2}-1}$) and $v_i
\in \ker(\go^{n-i})$, establishing the reverse inclusion.

So we have proved:
\begin{theo}{lefschetz in m2}
Let $i: M^2 \hookrightarrow X^{2n}$ be a symplectic embedding, $M$ and $X$ 
be compact 
and $\tilde{X}$ the blow up of $X$ along $M$. Then, for \e\ small enough,
\begin{itemize}
\item for $i >2$,
$$ \dim(\ker(\tilde{\go}^{n-i}))
= \dim(\ker(\go^{n-i}))
,$$
in particular, Lefschetz holds at level $i$ in $\tilde{X}$ if, and
only if, it does so in $X$;

\item if there is an element in $\ker(\go^{n-2})$ that restricts to a
nonzero element in $H^{2}(M)$ then
$$ \dim(\ker(\tilde{\go}^{n-2}))
= \dim(\ker(\go^{n-2}))
-1,$$
otherwise these kernels have the same dimension;

\item if there are elements $v_1,v_2 \in \ker(\go^{n-1})$ such that
$i^*(v_1 \wedge v_2) \neq 0$, then
$$ \dim(\ker(\tilde{\go}^{n-2}))
\leq \dim(\ker(\go^{n-2}))
-2,$$
otherwise
$$ \dim(\ker(\tilde{\go}^{n-2}))
 \leq \dim(\ker(\go^{n-2}))
.$$

\end{itemize}
\end{theo}

\subsection{The General case}
Now we treat the general case of the blow-up. Our main objective is to prove that if both $M^{2d}$ and $X^{2n}$
satisfy the Lefschetz property so does the blow-up of $X$ along $M$,
although in the course of this proof we obtain slightly  more, including a
generalization of lemma \ref{T:h1}. The first part of the proof was already
encountered at the end of the 2-dimensional case.

\begin{prop}{i>2d}
Assume that $(M^{2d},\sigma) \hookrightarrow (X^{2n},\go)$ is a
symplectic embedding with $M$ and $X$ compact and $2d < n$. Let $\tilde{X}$ be the blown-up
manifold. Then, for $ i> 2d$
$$\dim(\ker(\tilde{\go}^{n-i})) = \dim(\ker(\go^{n-i})).$$
In particular, $\tilde{X}$ will have the Lefschetz property at level
$i > 2d$ if, and only if, $X$ does so.
\end{prop}

\begin{remark}
The condition $2d < n$ is there only so that we can talk about a Lefschetz map at
level $i> 2d$, and this Proposition says that {\em we can not change the
dimension of the kernel of the Lefschetz map beyond the dimension of the
submanifold along which we are blowing-up.}
\end{remark}

\begin{proof}
First we return to our usual notation and let $k =
n-d$. Let $v_i \in \ker(\go^{n-i}) \subset H^{i}(X)$ and consider the
cohomology class $f^*(v_i) \in H^i(\tilde{X})$. The restriction of $v_i$ to $M$ is zero,
since the degree of $v_i$ is greater than the dimension of
$M$. Therefore $a v_i= 0$ and
$$(f^*(\go) + \e a)^{n-i} f^*(v_i) = f^*(\go^{n-i} v_i) = 0.$$
Thus, $f^*(\ker(\go)^{n-i}) \subset \ker(\tilde{\go}^{n-i})$.

On the other hand assume that $v= f^*(v_i) + a v_{i-2} + \cdots + a^{l} v_{i-2l}$ 
is an element of the kernel of $\tilde{\go}^{n-i}$. We may
further assume that the last term above, $v_{i-2l}$, is not zero or else
$v$ is of the form $f^*(v_i)$. From $v \in \ker(\tilde{\go})$ we have
\begin{align*}
0 &= (f^*(\go) + \e a)^{n-i}(f^*(v_i) + a v_{i-2} + \cdots + a^{l} v_{i-2l})\\
&= f^*(\go^{n-i} v_i) + \sum_{j=0,m=1}^{n-i,l} \e^{j}
\begin{pmatrix}
n-i\\j
\end{pmatrix}
a^{j+m}
\sigma^{n-i-j} v_{i-2m}.
\end{align*} 
Since $i> 2(n-k)$, the degree of the element above is $ 2n -i < 2k$
and therefore the highest power of $a$ in the expression above is
still smaller than $k$. Hence the coefficient of $a^{l+n-i}$, which is 
$v_{i-2l}$, must vanish. Thus we had from the beginning $v=f^*(v_i)$
and the expression above reduces to
$$0 = (f^*(\go) + \e a)^{n-i}(f^*(v_i)) = f^*(\go^{n-i}v_i)$$
and $v \in f^*(\ker(\go^{n-i}))$, which shows the reverse inclusion
and proves the proposition.
\end{proof}

\begin{prop}{i=2d}
Assume that $i:(M^{2d},\sigma) \hookrightarrow (X^{2n},\go)$ is a
symplectic embedding with $M$ and $X$ compact and $2d < n$. Let $\tilde{X}$ be the blown-up
manifold. If there is a $v \in \ker(\go^{n-2d})$ such that $i^*v \neq
0$, then
$$\dim(\ker(\tilde{\go}^{n-2d})) =  \dim(\ker(\go^{n-2d})) -1,$$
otherwise these kernels have the same dimension, as long as \e\ is
small enough. In particular, if $X$ 
has the Lefschetz property at level $2d$, so does $\tilde{X}$.
\end{prop}
\begin{proof}
Initially we observe that the same argument used in lemma
\ref{T:lemma2} shows that the existence of $v \in \ker(\go^{n-2d})$
such that $i^*v \neq 0$ is equivalent to the fact that the Thom class, $t$, 
of the embedding is not in the image of $\go^{n-2d}$. Now we let $k= n-d$
and  write down the Lefschetz map at level $2d$
\begin{equation}\label{E:i=2d}
\begin{aligned}
(f^*(&\go) + \e a)^{n-2d}(f^*(v_{2d})+a v_{2d-2} + \cdots + a^d v_0) = f^*(\go^{n-2d}v_{2d} - \e^{n-2d} v_0 t) +\\
&+ \sum_{i=k-d}^{k-1}
a^i\left(\left( \sum_{l \geq i - n+2d}^{d} \begin{pmatrix} n -2d\\i-l
\end{pmatrix} \e^{i-l}\sigma^{n-2d-i+l} v_{2(d-l)}\right)-  \e^{n-2d} v_0
c_{k-i}\right), 
\end{aligned}
\end{equation}
where the $c_i$'s are the Chern classes of the normal bundle of $M$. Then we
claim that we can make it equal to any element in 
\begin{equation}\tag{$*$}
f^*(\im(\go^{n-2d})\oplus\mbox{span}\{t\})\oplus a^{k-d} H^{2d}(M) \oplus \cdots
a^{k-1} H^{2}(M).
\end{equation}

The idea is the following: the system above is triangular
and therefore easy to solve.
Indeed, let $f^*(w_{2(n-d)})+a^{k-d} w_{2d} +\cdots + a^{k-1} w_{2}$ be an
element of the space ($*$). We start by choosing $v_{2d}$ and $v_0$ so that $\go^{n-2d} v_{2d}
- 
\e^{n-2d} t v_0$ equals $w_{2(n-d)}$. Observe that we can still change
$v_{2d}$ by any element in the kernel of $\go^{n-2d}$. Now look at the
coefficient of 
$a^{k-1}$ in \eqref{E:i=2d}:
$$ \e^{n-2d} v_2 + (n-2d) \e^{n-2d-1} \sigma v_0 - \e^{n-2d} v_0 c_1:=\e^{n-2d}v_2 + 
F(v_0).$$
Since we have already chosen $v_0$, we can now choose $v_2$ so that
the expression above equals $w_2$.

Assuming by induction that $v_{2j}$ have already been chosen for $j < j_0 < k-d$ so
that the coefficient of $a^{k-j}$ is $w_{2j}$ we see that the coefficient
of $a^{k-j_0}$ in \eqref{E:i=2d} is of the form
$$ \e^{n-2d} v_{2j_0} + F(v_0, \cdots, v_{2j_0-2}),$$
where $F$ is a function. Then again we can choose $v_{2j_0}$ so as to
have the desired equality.

Finally the coefficient of $a^{k-d}$ is of the form
$$ \e^{n-2d} i^*v_{2d} + F(v_0, \cdots, v_{2d-2}) \in H^{2d}(M).$$
And then, changing $v_{2d}$ by a multiple of the element in
$\ker(\go^{n-2d})$ whose restriction to $M$ is nonvanishing, we can make this
coefficient equal $w_{2d}$.

Now a simple counting of the dimensions involved shows that
$$\dim(\ker(\tilde{\go}^{n-2d})) = \dim(\ker(\go^{n-2d}))-1.$$

In order to prove the ``otherwise'' case, we start observing that if
for every $v\in \ker(\go^{n-2d})$, $i^*v =0$ then $f^*(\ker(\go)) \subset
\ker(\tilde{\go}^{n-2d})$. Therefore we immediately 
have $\dim(\ker(\go^{n-2d})) \leq \dim(\ker(\tilde{\go}^{n-2d}))$.

The reverse inequality is similar to what we have done so far and also to the subject of Proposition \ref{T:i<2d}, so we shall omit its proof.
\end{proof}

Before we can tackle the case $i < 2d$ we have to recall that, from Theorem \ref{T:yan2}, if $M$ satisfies the Lefschetz property there is a splitting of a cohomology class into primitive elements:
\begin{equation}\tag{\ref{E:primitive splitting}}
H^{i} = P_i \oplus \im(\go),
\end{equation}
where $P_i$ is defined by
$$P_i = \{v \in H^{i} \mid \go^{d-i+1}v =0 \},$$
if $i\leq d$ and $P_i = \{0\}$ otherwise. The elements in $P_i$
are called {\em primitive $i$--cohomology classes}.

Hence, we can write every $v \in H^{i}$, $i\leq j$ in a unique way as $v = v^0 + v^1
\sigma + \cdots + v^{[i/2]}\sigma^{[i/2]}$, with $v^j$
primitive. Observe that if $i>d$, then the first few terms in this
decomposition will vanish simply because $P_j = \{0\}$ for
$j>d$. Again, the notation for the splitting above will be used
consistently in the sequence.

\begin{prop}{i<2d}
Let $i: (M^{2d},\sigma) \hookrightarrow (X^{2n},\go)$ be a symplectic
embedding with $M$ and $X$ compact and $2d < n$. Assume further that $M$ satisfies the
Lefschetz property. Then, for \e\ small enough and $i \leq 2d$,
$$\dim(\ker(\tilde{\go}^{n-i})) \leq \dim(\ker(\go^{n-i})).$$ 
In particular, if $X$ satisfies the Lefschetz property at level $i$ so does
$\tilde{X}$.
\end{prop}
\begin{proof}
Firstly we observe that the cases of $i$ odd and $i$ even
can be treated similarly, but for simplicity we shall work out only the 
even
case: $2i$. 

We want to take the limit $\e \into 0$ in the map
$\tilde{\go}^{n-2i}$, but, as it stands, the
resulting map will clearly have a big kernel. So, what we shall do is
to find linear maps $A_{\e}$ and $B_{\e}$ such that $\lim_{\e \into 0}
B_{\e}\tilde{\go}^{n-2i}A_{\e}$ has kernel $f^*(\ker(\go^{n-2i}))$. From this we shall
conclude that the dimension of the kernel of $\tilde{\go}^{n-2i}$ is
at most the dimension of the kernel of $\go^{n-2i}$ as long as \e\ is
small enough.

We define $A_{\e}:H^{2i}(\tilde{X}) \into H^{2i}(\tilde{X})$ by
$$A_{\e}\left(f^*(v_{2i}) + \sum_{j=0}^{i-1} a^{i-j} v_{2j}\right) =
f^*(v_{2i})  + 
\sum_{j=0}^{i-1} \frac{1}{\e^j}a^{i-j} v_{2j}.$$
And $B_{\e}: H^{2n-2i}(\tilde{X}) \into H^{2n-2i}(\tilde{X})$ by
$$B_{\e} \left(f^*(v_{2n-2i}) + \sum_{j=0}^{i-1} a^{n-d-i+j}
v_{2d-2j}\right) =f^*(v_{2n-2i}) + \sum_{j=0}^{i-1}
\frac{1}{\e^{n-d-2i+j}}a^{n-d-i+j} v_{2d-2j} $$

Now we move on to write the map $\lim_{\e\into 0}B_{\e} \tilde{\go}^{n-2i}
A_{\e}$:
$$\lim_{\e\into 0}B_{\e} \tilde{\go}^{n-2i}
A_{\e} = f^*(\go^{n-2i}v_{2i}) + \sum_{j=0}^{i-1}
a^{n-d-i+j}\sum_{l=0}^{i-1}b_{d-j-l}\sigma^{d-j-l}v_{2l},$$
where $b_j= \begin{pmatrix} n-2i \\ j \end{pmatrix}$ are the binomial
coefficients.

We can further split the cohomology classes $v_{2l}$ into their primitive
parts, according to lemma \ref{T:yan}, $v_{2l} = v_{2l}^0 + \sigma
v_{2l}^1 + \cdots +\sigma^l v_{2l}^l$. With that, elements of
$H^{2i}(\tilde{X})$ will be in the kernel of the map above only if the
coefficients of $a^j\sigma^l$ vanish. The only terms that will give us information about primitives of degree $2l$ are the coefficients of $a^{k-i+l} \sigma^{d-2l}, a^{k-i+l+1} \sigma^{d-2l-1}, \dots, a^{k-1}\sigma^{d-l-i+1}$, and the vanishing of these is equivalent to the following:
$$\begin{pmatrix}
b_{d-2l} & b_{d-2l-1} & \cdots & b_{d-l-i+2} &  b_{d-l-i+1}\\
b_{d-2l-1} & b_{d-2l-2} & \cdots & b_{d-l-i+1} &  b_{d-l-i}\\
\vdots&       & \ddots & &\vdots\\
b_{d-2l-i+2} & b_{d-2l-i+1} & \cdots & b_{d-2i+4} &  b_{d-2i+3}\\
b_{d-2l-i+1} & b_{d-2l-i} & \cdots & b_{d-2i+3} &  b_{d-2i+2}
\end{pmatrix}
\begin{pmatrix}
v_{2l}^0\\
v_{2l+2}^1\\
\vdots\\
v_{2i-4,}^{i-2-l}\\
v_{2i-2}^{i-1-l}
\end{pmatrix}
=0$$
in the case $2i <d$, and a similar matrix for $2i>d$. What is important here is that in both cases the matrix will be constant along its anti-diagonals (it is a Toeplitz matrix) and the top right entry is nonzero.
Now, if we can prove that all the matrices above are invertible, we
will conclude that $f^*(v_{2i}) + \sum a^{i-j} v_{2j}$ is in the
kernel of $\lim B_{\e}\tilde{\go}^{n-2i}A_{\e}$  if and only if
$v_{2j}=0$ for all $j < i$ and $v_{2i} \in
\ker\{\go^{n-2i}\}$. So the next lemma finishes the proposition.
\noqed\end{proof}

\begin{lem}\label{T:oliver}
Let $b^n_j = \begin{pmatrix} n\\ j \end{pmatrix}, n, j \in \N$. Then for any $p \in \N$
$$\Delta^{n,p+1}_k = \det \begin{pmatrix}
b^n_{k+p} & b^n_{k+p-1} & \cdots & b^n_{k+1} &  b^n_{k}\\
b^n_{k+p-1} & b^n_{k+p-2} & \cdots & b^n_{k} &  b^n_{k-1}\\
\vdots&       & \ddots & &\vdots\\
b^n_{k+1} & b^n_{k} & \cdots & b^n_{k-p+2} &  b^n_{k-p+1}\\
b^n_{k} & b^n_{k-1} & \cdots & b^n_{k-p+1} &  b^n_{k-p}
\end{pmatrix}
\neq 0
$$
if $b^n_k \neq 0$. 
\end{lem}
\begin{proof}
Initially we observe that $b^n_k \neq 0$ if and only if $n\geq k \geq 0$ and for $n=k$ the matrix above has zeros above the anti-diagonal and ones on it, so the determinant is a power of $-1$. Further, by adding to each column the one to its right and using the binomial identity $b^n_k + b^n_{k-1} = b^{n+1}_k$ we get
$$\Delta^{n,p+1}_k = \det \begin{pmatrix}
b^{n+p}_{k+p} & b^{n+p-1}_{k+p-1} & \cdots & b^{n+1}_{k+1} &  b^n_{k}\\
b^{n+p}_{k+p-1} & b^{n+p-1}_{k+p-2} & \cdots & b^{n+1}_{k} &  b^n_{k-1}\\
\vdots&       & \ddots & &\vdots\\
b^{n+p}_{k+1} & b^{n+p-1}_{k} & \cdots & b^{n+1}_{k-p+2} &  b^{n}_{k-p+1}\\
b^{n+p}_{k} & b^{n+p-1}_{k-1} & \cdots & b^{n+1}_{k-p+1} &  b^n_{k-p}
\end{pmatrix}
$$
Now it is easy to check that 
$$\Delta^{n+1,p+1}_k  =\frac{(n+p+1)!(n-k)!}{n!(n+p-k+1)!}\Delta^{n,p+1}_k,$$
showing that $\Delta^{n+1,p+1}_k$ is nonzero if $\Delta^{n,p+1}_k$ is nonzero and we obtain the result by induction.
\end{proof}

These three propositions give us the following

\begin{theo}{lefschetz}
Let $i: (M^{2d},\sigma) \hookrightarrow (X^{2n},\go)$ be a symplectic
embedding with $M$ and $X$ compact and both satisfying the
Lefschetz property and $2d < n$. Let $(\tilde{X},\go+\e \ga)$ be the blow-up of $X$ along $M$ with the symplectic form from theorem \ref{T:cohomology}. Then, for \e\ small enough, $\tilde{X}$ also satisfies the Lefschetz property.
\end{theo}

\section{Massey Products and the Blow-up}\label{S:massey}

Having determined how the Lefschetz property behaves under blow-up, we turn our attention to formality.
As remarked in Chapter \ref{minimalmodels}, a very useful tool used to prove that manifolds are not formal are the Massey products.  The object this section is to prove that under mild codimension conditions, these products are preserved in the blow-up. This will allow us to find examples of nonformal symplectic manifolds in the next section.

\begin{theo}{massey products}
Let $i:M^{2(n-k)} \hookrightarrow X^{2n}$ be a symplectic embedding  with
$M$ compact
and let $\tilde{X}$ be the blown-up manifold, then:
\begin{itemize}
\item if $X$ has a nontrivial triple Massey product, so does $\tilde{X}$,
\item {\em (Babenko and Taimanov \cite{BT00a})} if $M$ has a nontrivial triple Massey product and $k > 3$, so does $\tilde{X}$.
\end{itemize}
\end{theo}
\begin{proof}
We start with the first claim and assume the Massey product $\langle v_1, v_2, v_3
\rangle$ is nonzero
in $X$. This means that there is $u$ representing such a product
with $[u] \not \in \mathcal{I}([v_1], [v_3])$, the ideal generated by $[v_1]$ and
$[v_3]$ in $H^*(X)$. If we consider the
product $\langle f^*v_1, f^*v_2,f^*v_3 \rangle$ we see that $f^*u$ is a
representative for it. The question then is whether $f^*[u]$ is in
the ideal $(f^*[v_1], f^*[v_3])$. Let us assume there was a relation
of the kind 
$$f^*[u] =  f^*[v_1](f^*\xi_1+a \zeta_1^1 + \cdots a^{k-1}
\zeta_1^{k-1}) +  f^*[v_3](f^*\xi_3+a \zeta_3^1 + \cdots a^{k-1}
\zeta_3^{k-1})$$  
Then, using the product rules \eqref{E:product rules},
$$f^*[u] = f^*([v_1]\xi_1 +[v_3] \xi_3) + a (i^*[v_1] \zeta_1^1 +  i^*[v_3]
\zeta_3^1)+ \cdots + a^{k-1}  (i^*[v_1] \zeta_1^{k-1} + i^*[v_3] \zeta_3^{k-1}).$$
Now, since the sum above is a direct one, all the
coefficients of the powers of $a$ must vanish and the following must hold:
$$f^*[u] = f^*( [v_1] \xi_1 + [v_3] \xi_3).$$ 
Since $f^*$ is an injection, we conclude that $[u] \in ([v_1], [v_3])$
which contradicts our initial assumption. 

Now we treat the second case. We
start by assuming that $v_1$, $v_2$ and $v_3 \in \gO(M)$ are closed forms
satisfying
$$ v_1 \wedge v_2 = dw_1 ~~~\mbox{ and }~~~ v_2\wedge v_3 = dw_2,$$
with $[w_1 v_3 -(-1)^{|v_1|} v_1 w_2] \not \in ([v_1],[v_3])$.
Letting $\gf:\tilde{V} \into V$ be the map of diagram \eqref{E:diagram}
and $\pi:V \into M$ the projection of the disc bundle, we have the
following relations in $H^*(\tilde{X})$
$$ \ga \gf^*\pi^* v_1 \wedge \ga \gf^*\pi^* v_2 = d (\ga^2
\gf^*\pi^*w_1) ~~~ \mbox{ and } ~~~  \ga \gf^*\pi^* v_2 \wedge \ga
\gf^*\pi^* v_3 = d (\ga^2 \gf^*\pi^*w_2)$$
The question then is again whether the cohomology class of the form
$$\ga \gf^*\pi^*v_1 \ga^2 \gf^*\pi^* w_1-(-1)^{|v_1|}\ga \gf^*
\pi^*v_1 \ga^2 \gf^*\pi^*w_2$$
is in the ideal generated by $a [v_1]$ and $a [v_3]$.

Suppose it was. Then there would be a relation of the type
\begin{align*}
a^3  [w_1 v_3 - (-1)^{|v_1|}v_1 w_2] &= a[v_1](f^*\xi_1+a \zeta_1^1+
\cdots + a^{k-1} \zeta_1^{k-1}) +\\& + a[v_3](f^*\xi_3+a \zeta_3^1+
\cdots + a^{k-1} \zeta_3^{k-1})\\
&= a([v_1] i^*\xi_1 [v_3] i^*\xi_3)+ a^2([v_1] \zeta_1^1 + [v_3]
\zeta_3^1) + \cdots +\\ &+ a^{k-1} ([v_1] \zeta_1^{k-2} + [v_3]
\zeta_1^{k-2}) + a^k ([v_1] \zeta_1^{k-1} + [v_3] \zeta_3^{k-1}).
\end{align*}
Expanding $a^k$ and using again that the result is a direct sum, we
look at the coefficient of $a^3$. Comparing both sides we see that it equals $[w_1 v_3 -
(-1)^{|v_1|}v_1 w_2]$, so
$$[w_1 v_3 - (-1)^{|v_1|}v_1 w_2] = [v_1] \zeta_1^2 + [v_3]
\zeta_3^2 - c_{k-3} \zeta_1^{k-1} [v_1] + c_{k-3} \zeta_3^{k-1} [v_3].$$
But this contradicts the fact that  $[w_1 v_3 -(-1)^{|v_1|} v_1 w_2]
\not\in ([v_1],[v_3])$.
\end{proof}

\section{Examples}\label{S:nonformal examples}

We want to find a 6-dimensional symplectic manifold with nontrivial
triple Massey product, $b_1$ even, not satisfying the Lefschetz
property and then use the blow-up procedure to eliminate the kernel of
$\go^i$. We have already encountered manifolds with these properties before, namely, nilmanifolds.

It is a result of Benson and Gordon \cite{BG88} that nontoroidal nilmanifolds never satisfy the Lefschetz property. Also, Nomizu's theorem \ref{T:Nomizu} implies that the Lie algebra of the corresponding Lie group with its differential $(\wedge^{\bullet}\Gg^*,d)$ furnishes a minimal model for the nilmanifold, therefore, no nontoroidal nilmanifold is formal. Indeed, it is a result of Cordero {\it et al} \cite{CFC87} that they always have nontrivial (maybe higher order) Massey products.

The simplest nilmanifold with the properties required is the one obtained from the product of two copies of the Heisenberg group.

\begin{ex}{heisenberg}
If $G$ is the 3 dimensional Heisenberg group then the Lie algebra has a basis 
formed by the left invariant vector fields whose values at the
identity are 
$$
\del_1 = \begin{pmatrix}
0&1&0\\
0&0&0\\
0&0&0\\
\end{pmatrix};~~~
\del_2 = \begin{pmatrix}
0&0&0\\
0&0&1\\
0&0&0\\
\end{pmatrix};~~~
\del_3 = \begin{pmatrix}
0&0&-1\\
0&0&0\\
0&0&0\\
\end{pmatrix}
$$

Then we check that $[\del_1,\del_2] = -\del_3$ and $[\del_1, \del_3] = [\del_2,\del_3] =0$.
Therefore, the quotient manifold, $\mathbb{H}$, of the 3--Heisenberg group
by the lattice generated by $\exp \del_i$ will have Chevalley-Eilenberg 
complex $(\wedge^{\bullet}\frak{g}^*,d)$  generated
by the invariant 1-forms $e_i$ dual to the $\del_i$ related by
$$ de_1 = de_2 =0;~~~de_3 = e_1\wedge e_2.$$
So, this manifold is the nilmanifold $\mathbb{H} = (0,0,12)$.

Hence $H^1(\mathbb{H})$ is generated by $\{e_1,e_2\}$,
$H^2$, by $\{e_{13}, e_{23}\}$ and $H^3$, by $\{ e_{123}\}$, where, as usual, $e_{ij}$ is the shorthand for $e_i\wedge e_j$. Therefore $b_1= b_2= 2$ and $b_0 = b_3 =1$.
\vskip6pt
\noindent 
{\bf Massey products.}
With the relations above for $e_1$, $e_2$ and $e_3$ we get that
$$ e_1 \wedge e_2 = de_3 ~~~~ \mbox{ and } ~~~~ e_2 \wedge e_1 =
d(-e_3).$$

Therefore we can form the Massey product $\langle e_1, e_2, e_1 \rangle = 
e_3 \wedge e_1
+ e_1 \wedge (-e_3) = - 2 e_{13} \neq 0$. Observe that in this
case, since $e_1 \wedge e_2$ is exact, Massey products have no
indeterminacy and the above is a nontrivial one in $\mathbb{H}$.
\end{ex}

\begin{ex}{gil}
Now consider the product $\mathbb{H} \times \mathbb{H} = (0,0,12,0,0,45)$.
The triple product $\langle e_1, e_2,e_1 \rangle$ is still nonzero. Further, the form
$$\go = e_{14} + e_{23} + e_{56}$$
is closed and has top power $6 e_{123456}$, which is
everywhere nonvanishing. Hence $\go$ is a symplectic form in
$\mathbb{H} \times \mathbb{H}$.

It is easy to see
that the kernel of $\go:H^2 \into H^4$ is $\mbox{span}\{e_{25}\}$ and the kernel of $\go^2$ in $H^1$ is $\mbox{span}\{e_2,e_5\}$.

In $\mathbb{H}$, consider the path
$$\ga(t) = \exp(t(\del_1+\del_2+\del_3)) =
\begin{pmatrix}
1 & t & t+ \dfrac{1}{2} t^2\\
0 & 1 & t \\
0 & 0 & 1
\end{pmatrix}, ~~t \in [0,2],$$
then $\ga' = \del_1 + \del_2 + \del_3$ and, besides this, $\ga(2) \approx \ga(0)$,
hence this is a circle.

In $\mathbb{H} \times \mathbb{H}$ there are two copies of $\ga$ (one in
each factor) making 
a torus $T^2$. A basis for the tangent space of this torus is given by
$\{\del_1+\del_2+\del_3, \del_4 + \del_5 + \del_6\}$. The symplectic form
evaluated on this basis equals 1 everywhere, hence this torus is a
symplectic submanifold. On the other hand, $e_{25}$
evaluated on this basis also equals 1 everywhere.

Therefore, by Theorems \ref{T:lefschetz in m2} and \ref{T:massey products}, the blow-up $M^6$ of $\mathbb{H} \times \mathbb{H}$ 
along this torus satisfies the Lefschetz property (for \e\ small enough)
and has a nontrivial triple product.
\end{ex}

\begin{ex}{munoz}
Still let $M$ be the manifold from the previous example. Using Fern\'andez and Mu\~noz's result on formality of Donaldson submanifolds \cite{FM02}, the Poincar\'e dual $N^4$ of the symplectic class $[\go]$ (after changing the symplectic form slightly and re-scaling to get integral periods) will be  a symplectic submanifold of $M$ which still satisfies the Lefschetz property and by Donaldson's theorem \cite{Do96} the inclusion $N \hookrightarrow M$ induces an isomorphism $H^1(M) \cong H^1(N)$ and an injection $H^2(M) \hookrightarrow H^2(N)$. Now, the Massey product in $M$ comes from three 1-forms and therefore still exists in $N$ and further, since we have an injection in $H^2$, this product is nonzero in $N$. So $N$ is a nonformal symplectic 4-manifold satisfying the Lefschetz property.
\end{ex}

\begin{ex}{gil2}
Let $(N^4,\sigma)$ be the manifold obtained in Example \ref{T:munoz}
that has a nontrivial triple product and satisfies the Lefschetz property.
By construction, $\sigma$ is an integral cohomology class, therefore, $(N, \sigma)$ can be symplectically embedded in
$\mathbb{C}P^{6}$, by Gromov's Embedding Theorem \cite{Gr71,Ti77}. By Theorem
\ref{T:lefschetz}, the blow-up of $\mathbb{C}P^{6}$ along $N$ will have the
Lefschetz Property. According to Theorem \ref{T:massey products}, it will have a
nonvanishing triple product (and thus is not formal) and from Theorem
\ref{T:cohomology} it is 
simply connected.
\end{ex}


\chapter{T-duality}\label{tduality}

T-duality in physics is a symmetry which relates IIA and IIB string theory and T-duality transformations act on spaces in which at least one direction has the topology of a circle.
In this chapter, we consider a mathematical version of T-duality introduced by Bouwknegt, Evslin and Mathai for principal circle bundles with nonzero twisting 3-form $H$ (also called $H$-flux) \cite{BEM03,BEM03b,BHM03}. Our main result is that any {\it invariant twisted \gcs} (or \gcy\ or \gk\ or generalized metric Calabi-Yau, for that matter) can be transported to the T-dual circle bundle with the appropriate twist. Hence it becomes probable that \gc\ geometry is a natural place to study mirror symmetry. Moreover, most of the structures derived from a \gcs\ correspond to their counterparts in the T-dual, for example, the invariant forms in $\mc{U}^{k}$, the decomposition $d = \del + \delbar$ and even some submanifolds (invariant in some sense) can be transported.

These properties of T-duality are very interesting, as they match properties on a manifold and its dual, therefore allowing one to prove difficult results on one side by studying them in the probably simpler mirror. As an example of this idea, we prove that no 6-nilmanifold can have an invariant \gk\ structure. T-duality  also allows us to determine \gk\ structures on semi-simple Lie groups different from those in Example \ref{T:liegroups}.

Of course other geometric structures can be transported via T-duality. Notably, the results in this chapter and in \cite{CG04b} inspired Witt to use T-duality to transport other generalized structures, such as generalized $G_2$ and $\Spin{7}$ and their weaker versions \cite{Wi04}, thereby obtaining new examples of those structures. Also, in \cite{CG04b}, we study how a generalized metric $g+b$ (see Chapter \ref{gcss}, Section \ref{linear algebra of gcs}) transforms and obtain, in a geometrical way, the Buscher rules \cite{Bu87,Bu88}.

This chapter is organized in the following way. In the first section we introduce T-duality for principal circle bundles as presented in \cite{BEM03} and prove the main result in that paper stating that T-dual manifolds have isomorphic twisted cohomologies. In Section \ref{tduality and gcss} we prove that \gcss\ can be transported via T-duality (our main result) and from there we draw a series of corollaries.

The material of this chapter and further results on this subject will be present in a collaborative paper with Gualtieri \cite{CG04b}, who I thank for useful conversations about this topic and who helped me to shape the results the way they appear in this chapter.

\section{Topological T-duality}

The construction of the T-dual we introduce in this section is due to Bouwknegt, Evslin and Mathai and is the subject of a series of papers by these authors \cite{BEM03,BEM03b,BHM03}. They set T-duality in the context of principal torus bundles and then work with twisted cohomology. One further requirement is that the twisting 3-form $H$ (also called an NS-flux) is invariant under the circle action. As a consequence of these restrictions, when we want to consider geometric structures later on, we will be obliged to work only with invariant structures.

The starting point is a principal circle bundle with an {\it integral} closed 3-form $H$ over a manifold $M$: $(E,H)$. We choose a connection $\theta$ on $E$, so that $\theta(\del/\del\theta) =1$, where $\del/\del\theta$ is the vector field generated by a fixed element in the Lie algebra of $S^1$. The curvature of this bundle is $d\theta = F$. As $H$ is closed and invariant,
$$ 0=\mc{L}_{\del/\del\theta} H = \frac{\del}{\del\theta}\lfloor dH + d (\frac{\del}{\del\theta}\lfloor H) = d (\frac{\del}{\del\theta}\lfloor H).$$
So, if we write $H = \tilde{F} \theta + h$, with $\tilde{F}$ and $h$ pull-backs from $M$, we see that $\tilde{F}$ is closed. Also, as $H$ is integral, and $\int_{S^1}\theta =2\pi$, we get that $2\pi \tilde{F}$ is integral. Hence we can construct a circle bundle $\tilde{E}$ over $M$ and choose a connection form $\tilde{\theta}$ whose curvature is $\tilde{F}$. We associate the 3-form $\tilde{H} =  F \tilde{\theta} + h$ to $\tilde{E}$:

\begin{definition}
The {\it $T$-dual} to $(E,H)$ is the space $(\tilde{E},\tilde{H})$ defined as the circle bundle over $M$ with a connection $\tilde{\theta}$  and curvature $d\tilde{\theta}=\tilde{F}$ with the 3-form $\tilde{H} = F \tilde{\theta} + h$.
\end{definition}
\setlength{\unitlength}{1mm}
\begin{picture}(60,30)(-95,0)
\thinlines
\put(-45,25){\makebox(0,0){$E$}}
\put(-57,20){\makebox(0,0){$H = \tilde{F} \theta + h$}}
\put(-57,15){\makebox(0,0){$d\theta = F$}}
\put(-30,5){\makebox(0,0){$M$}}
\qbezier(-42,23)(-37,15)(-33,8)
\qbezier(-35,9)(-34,8.5)(-33,8)
\qbezier(-33,10)(-33,9)(-33,8)

\put(-15,25){\makebox(0,0){$\tilde{E}$}}
\put(-3,20){\makebox(0,0){$\tilde{H} = F \tilde{\theta}+ h$}}
\put(-3,15){\makebox(0,0){$d\tilde{\theta} = \tilde{F}$}}
\qbezier(-18,23)(-23,15)(-27,8)
\qbezier(-25,9)(-26,8.5)(-27,8)
\qbezier(-27,10)(-27,9)(-27,8)
\end{picture}

\begin{remark}
In the definition above, the T-dual circle bundle $\tilde{E}$ is well defined by $(E,[H])$. However the corresponding connection $\tilde{\theta}$ is defined up to a closed form and the 3-form $\tilde{H}$ relies on the particular choice of connection. In the sequel we always assume such a choice was made when talking about `the' T-dual. 
\end{remark}

\begin{definition}
If $(E,H)$ and $(\tilde{E},\tilde{H})$ are T-dual to each other, the fiber product 
$$E\times_{M}\tilde{E} \hookrightarrow E \times \tilde{E}$$
is the {\it correspondence space}.
\end{definition}
Using the correspondence space, we can define the following map between the complexes of invariant differential forms:
\begin{equation}\label{E:tau}
\tau: \gO^{\bullet}_{S^1}(E) \into \gO^{\bullet}_{S^1}(\tilde{E}) \qquad \tau(\rho) = \frac{1}{2\pi}\int_{S_1} e^{-\theta \wedge \tilde{\theta}}\rho,
\end{equation}
where the $S^1$ where the integration takes place is the fiber of $E \times_{M} \tilde{E} \into \tilde{E}$, so the result is an invariant form in $\tilde{E}$. This formalism is suggestive, but, more concretely, any invariant form $\rho$ in $E$ can be written as $\rho = \theta \rho_1 + \rho_0$. In this case it is easy to check that
\begin{equation}\label{E:tau2}
\tau(\theta \rho_1 + \rho_0) = \rho_1 - \tilde{\theta}\rho_0.
\end{equation}

\begin{remark}
Needless to say, $\tau$ does not preserve degrees. Nevertheless it is well behaved under the $\Z_2$-grading of differential forms as $\tau$ reverses the parity of its argument:
$$\tau(\gO^{ev/od}_{S^1}) \subset \gO^{od/ev}_{S^1}.$$
\end{remark}

Also, it is clear from \eqref{E:tau2} and the setting that if we T-dualize twice and choose $\tilde{\tilde{\theta}} =\theta$, we get $(E,H)$ back and, on invariant forms, $\tau^2 =-\Id$. But $\gO^{\bullet}_{S^1}(E)$ is naturally a $\Z_2$-graded differential complex --- without the Leibniz rule --- with differential $d_H = d+ H$. Bouwknegt's main theorem \cite{BEM03} can be stated in the following way for forms:

\begin{theo}{BEM}
The map $\tau: (\gO^{\bullet}_{S^1}(E),d_H) \into (\gO^{\bullet}_{S^1}(\tilde{E}),-d_{\tilde{H}})$ is an isomorphism of differential complexes.
\end{theo}
\begin{proof}
Given that $\tau$ has an inverse, obtained by T-dualizing again, we only have to check that $\tau$ preserves the differentials, i.e., $-d_{\tilde{H}} \circ \tau = \tau\circ d_H $. To obtain this relation we use equation \eqref{E:tau}:
\begin{align*}
-d_{\tilde{H}}\tau(\rho) &= \frac{1}{2\pi}\int_{S^1} d_{\tilde{H}}(e^{-\theta\tilde{\theta}} \rho)\\
&= \frac{1}{2\pi}\int_{S^1} (H-\tilde{H})e^{-\theta\tilde{\theta}} \rho + e^{-\theta\tilde{\theta}} d \rho + \tilde{H}e^{-\theta\tilde{\theta}} \rho\\
&= \frac{1}{2\pi}\int_{S^1} H e^{-\theta\tilde{\theta}} \rho + e^{-\theta\tilde{\theta}} d \rho\\
& = \tau(d_H \rho)
\end{align*}
\end{proof}
\begin{remark}
If one considers $\tau$ as a map of the complexes of differential forms (no invariance required), it will not be invertible. Nonetheless, every $d_H$-cohomology class has an invariant representative, hence $\tau$ is a quasi-isomorphism.
\end{remark}

\begin{ex}{s3}
The Hopf fibration makes the 3-sphere, $S^3$, a principal $S^1$ bundle over $S^2$. The curvature of this bundle is a volume form of $S^2$, $\sigma$. So $S^3$ with zero twist is $T$-dual to $(S^2 \times S^1, \sigma \wedge \theta)$. The closed nonexact form $1 \in \wedge^{0}S^3$ is mapped to $\theta$, the volume form of $S^1$ in $S^2 \times S^1$ via the T-duality map while the volume form of $S^3$ is mapped into the volume form of $S^2$.
\end{ex}

\begin{ex}{2step nilmanifolds}
Using the notation of Chapter \ref{nilmanifolds}, consider a nilmanifold $E$ whose structure is given by $(0,0,0,c_1)$, where $c_1 \in \wedge^2\mbox{span}\{1,2,3\}$. $E$ is naturally a principal circle bundle over a 3-torus. If we take $H=0$ we obtain that the T-dual, $\tilde{E}$, is just the 4-torus with a 3-form $\tilde{H} = c_1 \wedge 4$.

If we had a 2-step nilmanifold $E$, i.e., $\nil(E)=2$, and $H=0$, we could write the structure of $E^{j+k}$ as $(0, \cdots,0, c_1, \cdots, c_k)$, with $c_i \in \wedge^2\mbox{span}\{1,\cdots,j\}$, indeed, as this nilmanifold is a principal $k$-torus bundle over a torus, we can even take $c_i$ to be the Chern classes of each circle bundle involved. This choice of circle bundles gives us a preferred way to T-dualize along the $k$ circles making the torus bundle. After T-dualizing, we obtain a $k+j$-torus with 3-form $\sum c_i \wedge i$.

This shows that every 2-step nilmanifold with vanishing 3-form can be T-dualized to a torus with nonvanishing 3-form.
\end{ex}


\begin{ex}{liegroups5} Let $(G,H)$ be a semi-simple Lie group with 3-form $H(X,Y,Z) = K([X,Y],Z)$, the Cartan form generating $H^3(G,\Z)$, where $K$ is the Killing form.

With a choice of maximal torus $T$ and a circle subgroup $S^1 <T$, we can think of $G$ as a principal circle bundle. For $X =\del/\del\theta \in \frak{g}$ tangent to $S^1 < G$ and of length $-1$ according to the Killing form, a natural connection on $G$ is given by $-K(X, \cdot)$. The curvature of this connection is given by
$$d(-K(X, \cdot))(Y,Z) = K(X, [Y,Z])= H(X,Y,Z),$$
hence the pull back of $F$ to $G$ is
$$F =  H(X,\cdot,\cdot) = X\lfloor H = \tilde{F}.$$
Which shows that semi-simple Lie groups with the Cartan 3-form are self T-dual.
Of course, one can repeat this with any other circle making up the maximal torus.
\end{ex}

\begin{remark}
Using the idea of the previous two examples we can T-dualize principal torus bundles by decomposing them into circle bundles and then T-dualizing each circle bundle in turn. One can show that the T-dual does not depend on the particular decomposition of the torus into circles, as different decompositions in the left side just induce different decompositions of the T-dual tori, but keep the principal torus bundle unchanged.

However, it is not always possible to T-dualize principal torus bundles, as we have done. To do so, some restrictions have to be placed on $H$. Namely, in order for a torus bundle $(E,H)$ to be T-dualizable, $H(X,Y, \cdot)$ has to vanish for every pair of vertical vectors $X$ and $Y$ \cite{BEM03}. In Example \ref{T:2step nilmanifolds}, this was trivially true, while in Example \ref{T:liegroups5} we had $H(X,Y,Z) =K([X,Y],Z)$, which vanishes for any $X,Y$ in the Lie algebra of the torus.
\end{remark}

\section{T-duality and Generalized Complex Structures}\label{tduality and gcss}

We have seen that T-duality comes with a map of differential complexes $\tau$ which is an isomorphism of the complexes of invariant differential forms. Now we prove that this map also preserves invariant twisted \gcss.

But first we introduce a map on $T+T^*$. Let $(E,H)$ and $(\tilde{E},\tilde{H})$ be a T-dual pair. Any invariant section of $TE\oplus TE^*$ can be written as $X + f \del/\del\theta + \xi + g \theta$, where $X$ is a horizontal vector and $\xi$ is a pull-back from $M$. We define $\gf: (TE \oplus T^*E)_{S^1} \into (T\tilde{E} \oplus T^*\tilde{E})_{S^1}$ by:
\begin{equation}\label{E:gf}
\gf(X + f \frac{\del}{\del\theta} + \xi + g \theta) = -X - g \frac{\del}{\del\tilde{\theta}} - \xi -f \tilde{\theta}.
\end{equation}

The core of our results is encoded in this theorem:

\begin{theo}{gf} The following hold:
\begin{enumerate}
\item the map \gf\ is orthogonal \wrt\ the natural pairing  on $TE \oplus T^*E$ and $T\tilde{E} \oplus T^*\tilde{E}$, hence it induces an isomorphism of Clifford algebras, $Cl((TE \oplus T^*E)_{S^1}) \cong Cl((T\tilde{E} \oplus T^*\tilde{E})_{S^1})$;
\item for $V \in C^{\infty}((TE \oplus TE^*)_{S^1})$ we have
$$\tau( V \cdot \rho) = \gf(V) \cdot \tau(\rho),$$
for any invariant form $\rho$. Therefore $\tau$ induces an isomorphism of Clifford modules $\wedge^{\bullet}T_{S^1}^*E \cong \wedge^{\bullet}T_{S^1}^*\tilde{E}$;
\item the map \gf\ is an isomorphism of Courant algebroids, $(C^{\infty}((TE \oplus T^*E)_{S^1}),[\cdot,\cdot]_H) \cong (C^{\infty}((T\tilde{E} \oplus T^*\tilde{E})_{S^1}),-[\cdot,\cdot]_{\tilde{H}})$, i.e., for $V_i \in C^{\infty}((TE \oplus T^*E)_{S^1})$
$$\gf([V_1,V_2]_H) = - [\gf(V_1), \gf(V_2)]_{\tilde{H}}.$$
\end{enumerate}
\end{theo}
\begin{proof}
\begin{enumerate}
\item It is obvious from equation \eqref{E:gf} that \gf\ is orthogonal with respect to the natural pairing and hence the first claim holds.

\item Splitting $\rho = \theta \rho_1 + \rho_0$ and $V = X+ f \del/\del\theta + \xi + g \theta$, equation \eqref{E:tau2} for $\tau$ gives the result:
\begin{align*}\tau(V\cdot \rho) & = \tau(\theta(-X \lfloor \rho_1 - \xi \rho_1 + g \rho_0) + X \lfloor \rho_0 + f \rho_1 + \xi \rho_0)\\
& = -X \lfloor \rho_1 - \xi \rho_1 + g \rho_0 +\tilde{\theta}( -X \lfloor \rho_0 - f \rho_1 - \xi \rho_0).
\end{align*}
While
\begin{align*}
\gf(V) \tau (\rho) &=  (-X - g \del/\del\theta - \xi - f \theta) (\rho_1 - \tilde{\theta} \rho_0)\\
& = -X \lfloor \rho_1 -\xi \rho_1 + g \rho_0 + \tilde{\theta}(-X \lfloor \rho_0 - \xi \rho_0 - f \rho_1).
\end{align*}

\item Finally, we have established that under the isomorphisms $\gf$ of Clifford algebras and $\tau$ of Clifford modules, $d_H$ corresponds to $-d_{\tilde{H}}$, hence the induced brackets (according to equation \eqref{E:twisted bracket}) are the same.
\end{enumerate}
\end{proof}

From this theorem and Theorem \ref{T:BEM}, we get that any structure defined on $E$ in terms of the natural pairing, Courant bracket and closed forms will correspond to one on $\tilde{E}$.

\begin{theo}{gcss and t-duality}
Any twisted invariant \gcs, \gcy\ structure or \gk\ structure on $E$ is transformed into a similar one via \gf.
\end{theo}
\begin{proof}
If $L < (TE\oplus T^*E)\tensor \C $ is the $+i$-eigenspace of an invariant $H$-\tgcs\ on $E$, then, by Theorem \ref{T:gf} (3), $\gf(L)$ is closed under the $\tilde{H}$-twisted Courant bracket on $(T\tilde{E} \oplus T^*\tilde{E})\tensor \C$. As \gf\ is orthogonal, $\gf(L)$ is still maximal isotropic, hence is a \gcs\ on $\tilde{E}$.

If $E$ has an $H$-twisted \gcy\ structure defined by a $d_H$-closed form $\rho$, with Clifford annihilator $L$, then the Clifford annihilator of $\tau(\rho)$ is $\gf(L)$, showing that $\tau(\rho)$ is pure, i.e., its annihilator has maximal dimension. By Theorem \ref{T:BEM}, $\tau(\rho)$ is $d_{\tilde{H}}$-closed, hence it induces an $\tilde{H}$-twisted \gcy\ structure on $\tilde{E}$.

If $L_1$, $L_2$ are twisted structures furnishing a \gk\ structure on $E$, then for $v \in L_1\cap L_2\backslash \{0\}$, 
$$\langle v, \overline v \rangle < 0,$$
and similarly for $v \in  L_1\cap \overline{L_2}$. As \gf\ is orthogonal \wrt\ the natural pairing, the same is true for $\gf(v) \in \gf(L_1)\cap \gf(L_2)\backslash \{0\}$. Therefore $\gf(L_1)$ and $\gf(L_2)$ induce a twisted \gk\ structure on $\tilde{E}$.
\end{proof}

As an application of this result we prove that no 2-step nilmanifold can be \gk.

\begin{theo}{nilmanifold and gk condition}
No 2-step nilmanifold admits a left invariant \gk\ structure. In particular, no 6-nilmanifold admits such structure.
\end{theo}
\begin{proof}
Recall that according to Theorem \ref{T:proposition5.17twisted} and the remark following it, any twisted \gk\ manifold admits a SKT structure. If a 2-step nilmanifold admits a \gk\ structure, according to Example \ref{T:2step nilmanifolds}, this nilmanifold can be T-dualized to a torus with nonzero 3-form, therefore furnishing the torus with an invariant SKT structure. But every invariant form in the torus is closed. In particular $d^c \go =0$ for the K\"ahler form induced by the metric and the complex structure, which can not happen in a SKT structure.

For the 6-dimensional case, we remark that Fino {\it et al} \cite{FPS02} have classified which 6-nilmanifolds admit invariant SKT structures (which would be the case for any admitting \gk\ structures) and those are all 2-step.
\end{proof}

As another easy application of Theorem \ref{T:gcss and t-duality} one can find further \tgk\ structures on Lie groups, following Examples \ref{T:liegroups} and \ref{T:liegroups5}.

Generalized Calabi-Yau metrics are also preserved by T-duality, but to prove that we need one extra lemma.
\begin{lem}
Given two invariant forms $\ga$ and $\gb$, define a form $\gamma$ in the base manifold by $(\ga,\gb) = \theta\wedge \gamma$. Then $(\tau(\ga),\tau(\gb)) = -\tilde{\theta}\wedge \gamma$.
\end{lem}
\begin{proof}
Write
$$\ga = \sum_{k=0}^{2n} \ga_0^{k} + u \ga_1^{k},$$
with $deg(\ga_i^k) = k$ and similarly for $\gb$. Using equation \eqref{E:mukai pairing} and expression \eqref{E:tau2} for $\tau$ a simple computation gives the result.
\end{proof}

\begin{remark}
If one thinks of $\tau$ as the Clifford action of $-\theta + \del/\del\theta$, and then relabeling $\theta$ to $\tilde{\theta}$, this lemma is just a consequence of the fact that
$$(v \cdot \ga, v \cdot \gb) = \langle v, v \rangle (\ga,\gb).$$
\end{remark}

\begin{theo}{ddbar and t-duality}
If two twisted \gcm s $(E,\J_1)$ and $(\tilde{E},\J_2)$ correspond via T-duality, then $\tau(\mc{U}^k_E) = \mc{U}^k_{\tilde{E}}$. In particular,
$$ \tau(\del_{E} \psi) = \del_{\tilde{E}}\tau(\psi) \qquad  \tau(\delbar_{E} \psi) = \delbar_{\tilde{E}}\tau(\psi)$$
and $E$ satisfies the $d_{H}d^{\J}$-lemma if and only if $\tilde{E}$ does (with $\tilde{H}$ in place of $H$).
\end{theo}
\begin{proof}
It is clear that $\tau(\mc{U}^n_E) = \mc{U}^n_{\tilde{E}}$ and that $\gf(\overline{L})$ is the $-i$-eigenspace of the induced \gcs\ on $\tilde{E}$. The result follows from Theorem \ref{T:gf}, item (2).
\end{proof}

Unfortunately this theorem is limited in its applications, at least as a way to generate new manifolds satisfying the $dd^{\J}$-lemma. For example, we have the following simple fact.

\begin{theo}{lefschetz and circle bundles}
The total space of a nonflat circle bundle over a compact manifold does not admit any symplectic structure for which the Lefschetz property holds.
\end{theo}
\begin{proof}
Although no invariance on the symplectic structure is required, we can always average the symplectic form to obtain an invariant form $\go = a \theta + b$ cohomologous to the original symplectic form, where $a$ and $b$ are pull-backs from the base and $\theta$ a connection form.

Since $\go$ is closed, we obtain that $da =0$. If $a$ is exact, say, $a = df$, then
$$\go^n = (n-1) a \theta b^{n-1} = (n-1) d(f \theta b^{n-1}),$$
as $d(\theta b^{n-1}) = 0$. Therefore $\go^{n}$ represents the trivial class, which can not happen in a compact symplectic manifold.

Hence it must be the case that $a$ represents a nontrivial cohomology class on the base and therefore its pull back is a nontrivial class on the total space and then we can consider the image of $[a]$ under the Lefschetz map $\go^{n-1}$:
$$[\go^{n-1} a] = [a b^{n-1}],$$
which is the pull-back of a top degree cohomology class from the base. As the circle bundle is not flat, this cohomology class vanishes in the total space and the Lefschetz property does not hold (at level 1).
\end{proof}

\begin{ex}{t-duality and the B-field}
One intriguing fact about T-duality in the present setting is that it does not preserve $B$-field actions. If $\rho = \theta \rho_1 + \rho_0$ is an invariant form and we let $B = \theta\wedge b_1 + b_2$ act on $\rho$ we get
\begin{equation}\label{E:t-duality and the B-field}
\tau(e^B \rho) = e^{b_2}\tau(\theta(\rho_1 + b_1 \wedge \rho_0) + \rho_0) = e^{b_2}((\tilde{\theta}-b_1) \wedge \rho_0 -\rho_1),
\end{equation}
which corresponds to the action of the not necessarily closed 2-form ${b_2}$ and the automorphism $\exp(b_1 \del/\del\tilde{\theta})$, as in Example \ref{T:endomorphisms}.

Another way to account for the appearance of $\tilde{\theta} - b_1$ instead of just $\tilde{\theta}$  in \eqref{E:t-duality and the B-field} is by saying that the action of a $B$-field forces a change in the connection in the T-dual from $\tilde{\theta}$ to $\tilde{\theta} - b_1$, but this however still does not explain satisfactorily the action of the nonclosed $b_2$.

In \cite{Gu03}, Chapter 8, Gualtieri proposes an alternative definition of T-duality using \gc\ submanifolds (for more details see \cite{CG04b}). That version of T-duality is such that a B-field action does not change structures on the T-dual, but it is not clear yet whether that interpretation agrees with the physics behind this theory.
\end{ex}

\begin{ex}{type change}
As even and odd forms get swapped with T-duality, the type of a \gcs\ is not preserved. However, it can only change, at a point by $\pm1$. Indeed, if $\rho = e^{B+i \go} \gO$ is an invariant form determining a \gcs\ there are two possibilities: If \gO\ is a pull-back from the base, the type will increase by 1, otherwise will decrease by $1$.

One would-be application of this fact is related to type changes in \gcy\ manifolds. This is because even if we start with a \gcy\ structure without type change on a circle bundle, the T-dual will still be a \gcy, but may have type change depending on whether $\gO$ is a pull back or not. For example, take $\R{5} \times S^1$ and consider the \gcs\ given by
$$\rho = e^{idx_5 \wedge d\theta}(dx_1 + i dx_2)(dx_3+ i dx_4 + (x_1+ix_2)d\theta).$$
The T-dual structure still in $\R{5} \times S^1$ is given by
$$\tau(\rho) =-e^{\frac{(dx_3 + i dx_4)(d\tilde{\theta} + i dx_5)}{x_1 + ix_2}}(x_1 + i x_2)(dx_1 + i dx_2),$$
which has type change along $x_1=x_2 =0$.

We could have replaced the coefficient of $d\theta$ in the expression for $\rho$ above by any function $f$, but the \gcy\ condition would only hold if $f = f(x_1,x_2)$ was a holomorphic function on $x_1 + i x_2$ and the type change in the T-dual would occur along the zeros of $f$. On the one hand, this suggests on what kind of locus we can expect to have type change, on the other hand this also suggests that this procedure won't work in compact examples, as then any holomorphic $f$ is constant.
Currently no example of compact \gcy\ with type change is known.
\end{ex}

\chapter{Formality of $k$-connected manifolds of dimension $4k+3$ and $4k+4$}\label{7manifolds}

In this chapter we leave the realm of generalized geometry and focus on formality of 7- and 8-manifolds. A theorem of Miller \cite{Mi79} states that any compact orientable $k$-connected manifold of dimension $d \leq (4k+2)$ is formal. In particular, a compact simply-connected $n$-manifold is formal if $n \leq 6$. Recently, Fern\'andez and Mu\~noz gave examples of simply-connected nonformal 7- and 8-manifolds \cite{FM02b} and Dranishnikov and Rudyak gave examples of $k$-connected nonformal $4k+3$- and $4k+4$-manifolds \cite{DrRu03}, therefore proving that Miller's theorem can not be improved without further hypotheses.

Here, we adopt the point of view that for a $k$-connected manifold, the smaller the $k+1^{th}$ Betti number, $b_{k+1}$, is, the simpler the topology. Then we establish the biggest value of $b_{k+1}$ for which one can still assure formality of an $n$-manifold, $n= 4k+3, 4k+4$:

\begin{theo}{miller improved}
A compact orientable $k$-connected manifold of dimension $4k+3$ or $4k+4$ with $b_{k+1} =1$ is formal.
\end{theo}


One motivation for the study of such manifolds comes from the existence of special geometric structures in 7- and 8-manifolds, e.g., $G_2$, $Spin(7)$ and symplectic structures. A compact irreducible Riemannian manifold with holonomy group $G_2$ or $Spin(7)$ has finite fundamental group and hence its universal cover is compact, simply-connected and has special holonomy. Other examples of manifolds with special holonomy are K\"ahler manifolds which are formal (see Theorem \ref{T:formality1}), hence it is conceivable that there is a connexion between the existence of special holonomy metrics and formality.

Another reason to study formality of 8-manifolds comes from symplectic geometry, as this is the only dimension where the question of existence of compact 1-connected nonformal examples is still open.

One property shared by irreducible $G_2$- and $Spin(7)$-manifolds is that they have a hard Lefschetz-like property. If $M^n$ is one such manifold, there is a closed $(n-4)$-form \gf\ for which:
\begin{equation}\label{E:motivation}
[\gf] \cup :H^2(M,\RR) \stackrel{\cong}{\into} H^{n-2}(M,\RR),
\end{equation}
is an isomorphism \cite{Jo00}. Of course, this last property is also shared by 2-Lefschetz symplectic 8-manifolds, where, by definition, \eqref{E:motivation} holds with $\gf = \go^2$.

We are interested in how \eqref{E:motivation} can be used to improve on Theorem \ref{T:miller improved}. We prove that if $M^n$, $n =4k+3,4k+4$, is $k$-connected and there is $\gf \in H^{n-2k-2}(M)$ for which
\begin{equation}\label{E:gf and h2}
\gf \cup :H^{k+1}(M,\RR) \stackrel{\cong}{\into} H^{n-k-1}(M,\RR)
\end{equation}
is an isomorphism and $b_{k+1}(M)=2$, then, $M$ is formal. If $n= 4k+3$ and $b_{k+1}(M)=3$, then all Massey products on $M$ vanish uniformly. We construct examples showing that our bounds are optimal.

This chapter is organized as follows. In Section \ref{miller}, we prove Theorem \ref{T:miller improved} by explicitly constructing a minimal model. In Section \ref{lefschetz}, we prove that the existence of the isomorphism \eqref{E:gf and h2} has a `formalizing tendency' in a $k$-connected orientable $4k+3$- or $4k+4$-manifold, as it implies formality for $b_2 =2$ and, in the $4k+3$-dimensional case, vanishing of Massey products for $b_2 =3$. In the last section we study examples from symplectic and Riemannian geometry where our results can be applied.

This chapter is based on \cite{Ca04c}.

\section{Extending Miller's bounds}\label{miller}

In this section we prove Theorem \ref{T:miller improved}. 
The way to prove formality in this case is by constructing the beginning of the minimal model for the manifold and then using the theorem of Fern\'andez and Mu\~noz (Theorem \ref{T:fernandezmunoz} above).

\vskip6pt
\noindent
{\sc Proof of Theorem \ref{T:miller improved}}. We will only prove the theorem for $4k+4$-manifolds, as the other case is analogous. The proof is accomplished by constructing the minimal model.

As $M$ is $k$-connected, $\mc{M}_{2k}$ is a free algebra generated by $H^{i}(M)$ in degree $i \leq 2k$ with vanishing differential and $\rho$ maps linearly each cohomology class to a representative form. The first time we may have to use a nonzero differential (and hence introduce one of the $B^j$ spaces) is in degree $2k+1$. This will be the case only if the generator $a \in \mc{M}_{k+1}$ satisfies $a^2 \neq 0$ and $\rho(a^2)$ is an exact form. Hence if either $\rho(a^2)$ is not exact or $a^2=0$,  all the spaces $B^j$ are trivial for $j \leq 2k+1$, showing that $M$ is $2k+1$-formal and hence, by Theorem \ref{T:fernandezmunoz}, formal.

So we only have to consider the case where the cohomology class $\ga \in H^{k+1}(M)$ satisfies $\ga^2 =0$ and $k+1$ is even. In this case, the Hirsch extension in degree $2k+1$ is given by
$$V^{2k+1} = H^{2k+1}(M)\oplus \mbox{span}\{b\},$$
where $d$ vanishes in $H^{2k+1}(M)$, $db = a^2$ and $\rho(b)$ is a form such that $d\rho(b) = \rho(a^2)$. With this splitting, $\mc{I}(\oplus_{j \leq 2k+1}(B^j)) < \mc{M}_{2k+1}$ is just the ideal generated by $b$ and to prove formality using Theorem \ref{T:fernandezmunoz} we have to show that any closed form in this ideal is exact in $\mc{M}$.

A closed form in this ideal, being the product of $b$ and an element of degree at least $k+1$, will have degree at least $3k+2$. Since $M$ is $k$-connected, Poincar\'e duality gives $H^{j}(M)=\{0\}$ for $3k+3 < j < 4k+4$. If an element in $\mc{I}(b)$ of degree $3k+2$ or $3k+3$ is closed and nonexact in $\mc{M}$, Poincar\'e duality implies that its dual is in either $\mc{M}_{k+2}^{k+2} \cong H^{k+2}(M)$ (in the former case) or in $\mc{M}_{k+1}^{k+1} \cong H^{k+1}(M)$ (in the latter). Either way, from there we can produce a degree $4k+4$ closed nonexact element in $\mc{I}(b)$. So we only have to check that any $4k+4$ closed element in $\mc{I}(b)$ is exact.

The only elements in $\mc{I}(b)$ in degree $4k+4$ are of the form $a b v$, with $v \in H^{k+2}(M) \cong V^{k+2}$, which have derivative $d(a b v) = a^3v$ and hence are not closed if $v \neq 0$. This shows that $M$ is $(2k+1)$-formal and therefore formal.
\qed
\vskip6pt

\section{Formality of hard Lefschetz manifolds}\label{lefschetz}

In this section we study compact orientable $k$-connected $n$-manifolds, $n=4k+3, 4k+4$, for which there is a cohomology class $\gf \in H^{n-2k-2}(M)$ inducing an isomorphism
\begin{equation}\tag{\ref{E:gf and h2}}
\gf \cup :H^{k+1}(M) \stackrel{\cong}{\into} H^{n-k-1}(M).
\end{equation}
We prove that this property has a `formalizing tendency' in the following sense. 

\begin{theo}{g2 with b2 leq 3} 
A compact orientable $k$-connected manifold $M^n$, $n=4k+3, 4k+4$, satisfying \eqref{E:gf and h2} with $b_{k+1}=2$ is formal. If $n=4k+3$ and $b_{k+1}=3$ all the Massey products vanish uniformly.
\end{theo}
\begin{proof}
We start considering the case $b_{k+1}=2$. The cases $n=4k+3$ and $n=4k+4$ are similar, so we only deal with the latter. Due to \eqref{E:gf and h2}, the class $\gf$ induces a nondegenerate bilinear form on $H^{k+1}(M)$. If $k$ is even, this bilinear form is skew, and if $k$ is odd, it is symmetric.

For $k$ odd, the signature of the bilinear form, i.e., the difference between the number of positive and negative eigenvalues, is either $2$, $0$ or $-2$. By changing $\gf$ to $-\gf$, the case of signature $-2$ can be transformed into signature 2, hence there are two possibilities to consider. As they are similar, we will only treat the signature 2 case.

Let $a_1,a_2$ ($da_1=da_2=0$) be generators of $\mathcal{M}_{k+1} = Sym^{\bullet}H^{k+1}(M)$, the
 first nontrivial stage of the minimal model for $M$. We may further assume that the bilinear form induced by \gf\ is diagonal in the basis $\{a_1,a_2\}$ and
$$ \int \gf a_i^2 =1.$$
As with Theorem \ref{T:miller improved}, $\mc{M}^{2k}$ is just the free algebra generated by $H^j(M)$ in degree $j \leq k$ with vanishing differential and the first time we may have to introduce one of the $B^j$ spaces is in degree $2k+1$, where
$$B^{2k+1} \cong \ker Sym^2H^{k+1}(M) \stackrel{\cup}{\into} H^{2k+2}(M).$$

We know that both
 $a_1^2$ and $a_2^2$ are nonzero in $H^{2k+2}(M)$, since by \eqref{E:gf and h2} they pair nontrivially with \gf. So
$$\dim\,B^{2k+1} = \dim(\ker Sym^2H^{k+1}(M) \stackrel{\cup}{\into} H^{2k+2}(M)) \leq \dim(Sym^2 H^{k+1}(M)) - 1 = 2,$$
and hence there may be at most
 two generators in degree $2k+1$ in $\mathcal{M}_{2k+1}$ to kill cohomology
 classes in $\mathcal{M}_{2k}$ that are not present in $H^{2k+2}(M)$. The case
 when only one generator is added has a proof very similar to the one
 of Theorem \ref{T:miller improved}, so we move on to the case when there are two generators $b_1$ and $b_2$ added to kill cohomology in degree $2k+2$.

Then
$$db_i = k^i_{11}a_1^2 + k^i_{12}a_1 a_2 + k^i_{22}a_2^2,\qquad i=1,2.$$
Multiplying by \gf\ and integrating we get $k^i_{11} =
-k^i_{22}$. Hence, by re-scaling and taking linear combinations we may assume that $db_1= a_1^2 - a_2^2$ and $db_2 = a_1 a_2$.

Now, if $c \in \mc{I}(b_1, b_2)_{\leq 2k+1}$ is a closed element (again, by Poincar\'e duality we may assume it has
degree $n$) we can write it as
$$ c = (a_1 c_{11} + a_2c_{12}) b_1 + (a_1 c_{21} + a_2c_{22}) b_2 + (k_1a_1 + k_2 a_2)b_1b_2,$$
where $c_{ij}$ are closed elements of degree $k+2$ and $k_i \in \RR$. The condition $dc=0$ implies that $c=0$, thus every closed form in $(b_1,
b_2)_{\leq 2k+1}$ is exact and $M$ is $2k+1$-formal and therefore formal.

The case when $k$ is even is easier. Indeed, the argument above can be used again, but with $Sym^2 H^{k+1}(M)$ replaced by $\wedge^2 H^{k+1}(M)$ and hence the nondegeneracy of the pairing implies
$$\dim\,B^{2k+1} = \dim(\ker \wedge^2 H^{k+1}(M) \stackrel{\cup}{\into} H^{2k+2}(M)) \leq \dim(\wedge^2 H^{k+1}(M)) - 1 = 0,$$
hence $M$ is trivially $2k+1$-formal and therefore formal.
\vskip6pt

To finish the proof, we consider a $4k+3$-manifold with $b_{k+1}=3$ and prove that all triple Massey products vanish if \eqref{E:gf and h2} holds. Initially we remark that $b_{k+1}=3$ and \eqref{E:gf and h2} can not happen if $k$ is even, as there is no nondegenerate skew bilinear form in an odd dimensional vector space.

We also observe that the Massey product $\langle a,b,c \rangle$ has
degree at least $3k+2$ and since $H^j(M)=\{0\}$, for $3k+2 < j < 4k+3$, this
product will vanish, whenever defined, if its degree is neither $3k+2$ nor $4k+3$. If $\langle a,b,c \rangle \in H^{4k+3}(M)$, it lies in the ideal
generated by $([a],[c])$, so the product also vanishes, therefore the only case left is when $a$, $b$ and $c$ have degree $k+1$ and the product lies in $H^{3k+2}(M)$.

This product vanishes trivially if $c= \lambda a$, hence we can assume $a$ and $c$ linearly independent. Since $\gf a b$ and $\gf b c$ vanish, but \gf\ induces a nondegenerate bilinear form, there is $\beta \in H^{k+1}(M)$ such that $\gf b \beta \neq 0$ and we can express the Massey product in the basis $\{\gf [a],\gf [\beta],\gf [c]\}$:
$$\langle a,b,c \rangle = -v_1 c + a v_2= k_1 \gf[a] + k_2\gf[\beta]+k_3\gf[c].$$
where $ab = dv_1$ and $bc= dv_2$. Multiplying the equation above by $-[b]$ and integrating over $M$ we get
$$ - k_2 \int \gf b \beta  = \int  - v_1 c b + a v_2 b = \int - v_1dv_2 + v_2dv_1 = \int d(v_1v_2) =0.$$
Thus $k_2=0$ and $\langle a,b,c \rangle = k_1 \gf[a] + k_3\gf[c] \in ([a],[c])$. So the Massey product vanishes.

Observe that if $ac = dv_3$, it is still possible to choose $v_i$, $i=1,2,3$ in such a way that the Massey products $\langle a,b,c \rangle$, $\langle b,c,a \rangle$ and $\langle c,a,b  \rangle$ vanish simultaneously. Indeed, let us assume $a$, $b$ and $c$ are linearly independent (the linearly dependent case is similar) and a set of choices of $v_i$ was made:
$$
\begin{cases}
\langle a,b,c \rangle = [v_1 c - a v_2] = k_1 \gf[a] + k_2\gf[c],&\\
\langle b,c,a \rangle = [v_2 a - b v_3] = l_1 \gf[b] + l_2\gf[a],&\\
\langle c,a,b \rangle = [v_3 b - c v_1] = m_1 \gf[c] + m_2\gf[b].&
\end{cases}
$$
Then we can set $v_1 = \tilde{v}_1 + k_2 \gf$, $v_2 = \tilde{v}_2 + l_2 \gf$ and $v_3 = \tilde{v}_3+ m_2 \gf$. With these choices we get
$$
\begin{cases}
\langle a,b,c \rangle = [\tilde{v}_1 c - a \tilde{v}_2] = \tilde{k} \gf[a],&\\
\langle b,c,a \rangle = [\tilde{v}_2 a - b \tilde{v}_3] = \tilde{l} \gf[b],&\\
\langle c,a,b \rangle = [\tilde{v}_3 b - c \tilde{v}_1] = \tilde{m} \gf[c].&
\end{cases}
$$
Adding them up, the left hand side vanishes, giving
$$0 = \tilde{k}\gf [a] + \tilde{l} \gf [b] + \tilde{m} \gf [c].$$
Since $\gf[a]$, $\gf[b]$ and $\gf[c]$ are linearly independent,
$\tilde{k}, \tilde{l}$ and $\tilde{m}$ vanish and, with these choices,
all the Massey products vanish {\it simultaneously}, hinting
at formality.
\end{proof}

\section{Examples}

In this section give simple examples showing that the bounds established in Theorem \ref{T:miller improved} are sharp. We also apply our results to the blow-up of $\C P^8$ along a symplectic submanifold and to Kovalev's examples of $G_2$-manifolds. We finish with an example of a compact 1-connected 7-manifold which satisfies all known topological restrictions imposed by a $G_2$ structure, has $b_2 =4$ and nonvanishing Massey products. This last example shows that the results of Theorem \ref{T:g2 with b2 leq 3} are sharp in 7 dimensions and that one can not answer the question of formality of $G_2$-manifolds using only the currently known topological properties of those.

\begin{ex}{sharp}
Now we construct $k$-connected $(4k+3)$-manifolds with $b_{k+1}=2$ which have nonvanishing Massey products and hence are not formal. This shows that the bounds obtained in Theorem \ref{T:miller improved} are sharp in these dimensions.

Let $V$ be a $(2k+2)$-vector bundle over $S^{k+1}\times S^{k+1}$ with nonzero Euler class $\chi$. Then, the total space, $X$, of the sphere bundle associated to $V$ is a compact $k$-connected $4k+3$-manifold. Using the Gysin sequence for this sphere bundle
$$ \cdots \stackrel{\chi\cup}{\into} H^j(S^{k+1}\times S^{k+1}) \stackrel{\pi^*}{\into} H^j(X) \stackrel{\pi_*}{\into}  H^{j-2k-1}(S^{k+1}\times S^{k+1}) \stackrel{\chi\cup}{\into} H^{j+1}(S^{k+1}\times S^{k+1}) \into \cdots$$
we see that $H^{k+1}(X) = \mbox{span}\{v_1, v_2\}$ and $H^{2k+2}(X) = \{0\}$, where $v_i$ are generators for the top degree cohomology of each sphere. Therefore, $v_1 \cup v_2 =0$. If $\go_i$ are volume forms pulled back from each $S^{k+1}$, then $\go_1 \wedge \go_2 = d\xi$, where
$$\int_{S^{2k+1}}\xi = \int_{S^{k+1} \times S^{k+1}}\chi.$$
Therefore we can compute the Massey product
$$ \langle v_1, v_2, v_2 \rangle = - [\xi \wedge \go_2].$$
This is not exact, as it pairs nontrivially with $v_1$,
$$\int_X \xi \wedge \go_1 \wedge \go_2 =  \int_{S^{k+1} \times S^{k+1}} \chi \cdot \int_{S^{k+1}\times S^{k+1}} \go_1 \wedge \go_2  \neq 0$$
and has no indeterminacy, hence $X$ is not formal.
\end{ex}

\begin{ex}{8symplectic}
Blowing up $\C P^n$ along a suitable submanifold, Babenko and Taimanov proved  in \cite{BT00b} that there are compact 1-connected nonformal symplectic manifolds in any dimension $2k \geq 10$. Due to Miller's theorem, such examples do not exist in dimensions 6 or less and the question is still open for 8-manifolds.

Using the techniques of chapter \ref{blowup}, one can show that the blow-up, $X$, of $\C P^n$ along any submanifold is always 2-Lefschetz and, if the submanifold is connected, $b_2(X) = 2$. Therefore, our results imply  the blow-up of $\C P^8$ along a connected symplectic submanifold is always formal. One can also try to blow up $\C P^n$ and then take a sequence of Donaldson submanifolds until the result is an 8-manifold, but the manifold obtained this way is also 2-Lefschetz \cite{FM02} and has $b_2 =2$ \cite{Do96}.
\end{ex}

\begin{ex}{kovalev}
In \cite{Kov00}, Kovalev produces a series of examples of 1-connected compact $G_2$-manifolds from pairs of Fano 3-folds via twisted connected sums. If $M^7$ is obtained from the Fanos $F_1$ and $F_2$, he proves that $b_2(M) \leq \min\{b_2(F_1),b_2(F_2)\} -1$. Our results imply that if either of the Fanos involved have $b_2(F_i) \leq 3$, $M^7$ is formal, while if $b_2(F_1)=4$ and $b_2(F_2) \geq 4$, the Massey products vanish. Hence the only possibility for one of his examples to have nontrivial Massey products is if it is constructed from 2 Fanos with $b_2(F_i) \geq 5$.

According to the classification of Fano 3-folds \cite{MoMu86}, these have $b_2 \leq 10$ and if $b_2 \geq 6$, the Fano is just the blow-up of $\C P^3$ in an appropriate number of points. It is easy to follow Kovalev's construction to prove that if one of the summands is $\C P^3$ with some points blown-up, $M^7$ is formal. The case where both the Fanos have $b_2=5$ is more difficult, but one still has formality.
\end{ex}

In \cite{Jo00}, Joyce proves the following
\begin{theo}{topological g2}
A compact Riemannian manifold with holonomy $G_2$ has finite $\pi_1$, nonvanishing first Pontryagin class $p_1$ and, if $\gf$ is the closed 3-form determining the structure,
\begin{equation}\label{E:g2 conditions}
\int a^2 \gf < 0 \mbox{ for } a \in H^{2}(M)\backslash \{0\} \qquad \mbox{ and } \qquad \int p_1 \wedge [\gf ] < 0.
\end{equation}
\end{theo}

\begin{ex}{m4crosss3}
The conditions in Theorem \ref{T:topological g2} are stronger than just \eqref{E:gf and h2}. For example, let $M^4$ be a simply-connected 4-manifold and take $\gf$ to be the volume form of $S^3$ in $X^7 = M^4 \times S^3$. Then the pairing induced by \gf\ in $H^2(X)$ is just the nondegenerate intersection form in $H^2(M)$ and hence \eqref{E:gf and h2} holds. However, no manifold obtained this way satisfies the conditions of Theorem \ref{T:topological g2}. Indeed, $H^3(S^3  \times M^4) = H^3(S^3)$. So the cohomology class of the 3-form \gf\ is the pull back of a multiple of the top   degree class in $S^3$ and as before the bilinear form induced by  \gf\ is just the intersection form of the 4-manifold. Therefore if $X$ satisfies the conditions of Theorem \ref{T:topological g2},   $M$ has (negative) definite intersection form. But by Donaldson's  Theorem (\cite{DK98}, Theorem 1.3.1) such an $M$ can not be spin unless $b_2(M)=0$. If $b_2(M)=0$, the Hirzebruch Signature Theorem in dimension 4 implies that $p_1(M) = 3 \tau(M)=0$, where $\tau$ is the signature of $M$, and hence $M$ (and $S^3 \times M$) has  vanishing first Pontryagin class.
\end{ex}

Although the topological properties of $G_2$ manifolds are stronger than \eqref{E:gf and h2}, they bring no extra information about formality of $G_2$-manifolds. Indeed, using circle bundles we can construct a nonformal manifold satisfying all the topological properties from Theorem \ref{T:topological g2}. The key is Wall's classification of 1-connected spin 6-manifolds.

\begin{theo}{Wall}
{\em (Wall \cite{Wa66})} Diffeomorphism classes of oriented 6-manifolds with
torsion-free homology and vanishing second Stiefel-Witney class
correspond bijectively to isomorphism classes of systems of
invariants:
\begin{itemize}
\item Two finitely generated free abelian groups $H^2$, $H^3$, the latter of even rank;
\item A symmetric trilinear map $\mu:H^2\times H^2\times H^2 \into Z$;
\item A homomorphism $p_1:H^2 \into Z$;
\item Subject to: for $x, y \in H^2$,
$$\mu(x,x,y) = \mu(x,y,y) ~~~\mod(2),$$
\noindent
\phantom{Subject to:} for $x \in H^2$,
$$p_1(x) = 4 \mu(x,x,x) ~~~ \mod(24). $$ 
\end{itemize}
\end{theo}

With appropriate choices for the pairing $\mu$ and for the Chern class of the principal circle bundle, we can obtain nonformal 7-manifolds satisfying \eqref{E:g2 conditions}. As the base manifold is spin, so will be the total space of the circle bundle. We finish this paper with one example constructed this way.

\begin{ex}{second 7manifold}
We let $H^2 = \langle\go, \ga_1,\ga_2,\ga_3,\gamma\rangle$ and
define the cup product on $H^2$ so as to have the following relations
\begin{center}
\vskip6pt
\noindent
\begin{tabular}{c c c c c}
$\go \gamma =0,$& $\go \ga_i \ga_j
= 2\gd_{ij},$& $\go^2 \ga_i=0,$& $\go^3=2,$ & $\go \ga_1 = \gamma \ga_3,$\\
&$\ga_1\ga_2=0,$& $\ga_2 \ga_3 =\gamma\ga_1$& and& $\ga_3 \go = \gamma \ga_1$.
\end{tabular}
\vskip6pt
\end{center}
One set of choices that gives the desired result is the following
\vskip6pt
\noindent
\begin{alignat*}{6}
\go \ga_i \ga_j &= 2\delta_{ij}&\qquad \ga_1\ga_2^2 &=0 &\qquad \ga_1^2 \ga_2 &= 0 & \qquad \ga_2 \ga_3^2 & =2 &\qquad \ga_3 \gamma^2 &=0 \\
\go \gamma \ga_i &=0&\qquad \ga_1 \ga_2 \ga_3 &=0&\ga_1^2 \ga_3 &= 0&\qquad \ga_2 \ga_3 \gamma &=0& \ga_3^2 \gamma &=0\\
\go \gamma^2 &=0& \ga_1 \ga_2 \gamma&=0 & \ga_1^2 \gamma &=0 & \ga_2 \gamma^2 & =2 & \ga_3^3&=0 \\
\go^2 \ga_i &=0& \ga_1 \ga_3^2 &= 0 &\ga_1^3 &=0 & \ga_2^2 \ga_3 &= 0 & \gamma^3&=0\\
\go^2 \gamma &=0& \ga_1 \ga_3 \gamma &= 2& & & \ga_2^2 \gamma &= 2& &\\
\go^3 &=2 & \ga_1 \gamma^2&=0 & & & \ga_2^3 &= 0& & \\
\end{alignat*}
With these choices, $\go \gamma=0$ and
this is the only 2-cohomology class that pairs trivially with $\gamma$.

Now we let $M$ be a simply connected spin 6-manifold with $H^2(M)=H^2$,
cup product as described above, arbitrary $H^3$ and first Pontryagin
class $p_1 = 4 \go^2$
. Let $w, a_i, c_1$ be  a set of  closed forms representing $\go, \ga_i, \gamma$ and let $X$ be a circle bundle
over $M$ with connection form $\theta$ and first Chern
class $\gamma$ with $d\theta = c_1$. Then $X$ is spin, has first Pontryagin class $p_1 =
4\go^2$ and has degree 2 and 3 cohomology $H^2(X)=\langle \go, \ga_1,
\ga_2, \ga_3 \rangle$, $H^3(X) 
= H^3(M)\oplus \langle [\theta\tensor \go] \rangle$. The term $[\theta \tensor \go]$ is the one we are concerned about: As $\gamma\go =0$, in the form level we have $c_1\wedge w  = d\xi$, for some 3-form $\xi$ pulled back from $M$, hence $\theta \wedge w - \xi$ is a closed form. This form represents the cohomology class $[\theta \tensor \go]$.

Letting $\gf = -\theta \wedge w +\xi$ we see
that \gf\ induces a negative definite bilinear form on $H^2(X)$ as
$$\int_X \gf \wedge a_i\wedge a_j = \int_M - w \wedge a_i \wedge a_j = -2\delta_{ij},$$
and similarly
$$\int_X \gf \wedge w\wedge a_j = 0 \qquad\mbox{and}\qquad \int_X \gf \wedge w\wedge w =  -2.$$

Also,
$$\int_X \gf p_1 = \int_X -\theta w 4w^2 = -\int_M4w^3 = - 8.$$

Finally, since $\gamma$ pulls back to zero in $H^2(X)$, we can define the
Massey products $\langle \go, \ga_1, \ga_2 \rangle$, $\langle \ga_1,
\ga_2,\ga_3\rangle$, $\langle \ga_2, \ga_3, \go \rangle$ and $\langle \ga_3,
\go, \ga_1\rangle$. To prove that $X$ is not
formal we compute $\langle \ga_1, \ga_2,\ga_3 \rangle$:
$$ a_1 a_2 = da_{12},\qquad a_2 a_3 = d(\theta a_1 - a_{23}),$$
where $a_{ij}$ are pull-backs from the base. So
$$\int_X \langle \ga_1, \ga_2,\ga_3 \rangle \go = \int_X a_{12} a_3 w - a_1
(\theta a_1-a_{23})w = - \int_X \theta  w a_1^2 = - \int_M w a_1^2 = - 2.$$
Which means that, for these choices, the Massey product pairs
nontrivially with the cohomology class $\go$ and therefore is a closed
nonexact form. One can also check that different choices {\it
  keep the integral above unchanged} so the Massey product does not
vanish.
\end{ex}

\begin{remark}
If we deal with the $Spin(7)$ case in a similar fashion, i.e., stripping off the Riemannian structure and working only with the implied topological properties, circle bundles will not provide possible examples of nonformal $Spin(7)$-manifolds. This is because $Spin(7)$-manifolds have $\hat{A}$-genus 1 and by a result of Atiyah and Hirzebruch \cite{AH70} if a compact connected Lie group acts differentiably and non-trivially on a compact  orientable spin manifold $X$, then $\hat A(X)=0$.
\end{remark}

\bibliographystyle{abbrv}
\bibliography{references}


\end{document}